\documentclass[12pt]{article} 

\usepackage{amsfonts,amsmath,amssymb,epsf}
\usepackage{mathrsfs}
\usepackage{graphics} 
\usepackage{ifpdf}
\ifpdf
\usepackage{epstopdf}
\usepackage{hyperref}
\else
\usepackage[hypertex]{hyperref}
\fi
\usepackage{xypic}

\advance\textheight by \topskip
\newcounter{thm}

\newcommand\sect[1]{\section{#1}\setcounter{equation}0\setcounter{thm}0} 
                   
\newtheorem{thm}{Theorem}[section]
\newtheorem{defn}[thm]{Definition}
\newtheorem{prop}[thm]{Proposition}
\newtheorem{cor}[thm]{Corollary}
\newtheorem{rema}[thm]{Remark}
\newtheorem{lemma}[thm]{Lemma}

\newcommand\void[1]       {}

\newcommand\be            {\begin{equation}}
\newcommand\bea           {\begin{eqnarray}}
\newcommand\bearll        {\begin{array}{ll}\displaystyle}

\newcommand\ee            {\end{equation}}
\newcommand{\eea}         {\end{eqnarray}}
\newcommand\eear          {\end{array}}
\newcommand\enl           {\\[.5em]\displaystyle}
\newcommand\etb           {&\!\! \displaystyle}
\newcommand\erf[1]        {(\ref{#1})}
\newcommand\labl[1]       {\label{#1}\ee}
\newcommand{\nn}          {\nonumber \\}

\newcommand\arxiv[2]      {\href{http://arXiv.org/abs/#1}{#2}}
\newcommand\doi[2]        {\href{http://dx.doi.org/#1}{#2}}

\newcommand\Acl           {A_{\text{\rm{cl}}}}
\newcommand\Aop           {A_{\text{\rm{op}}}}
\newcommand\CcC           {\Cc^2_{\pm}}
\newcommand\itco          {\tilde{\iota}_{\text{{\rm cl-op}}}}
\newcommand\ico           {\iota_{\text{{\rm cl-op}}}}

\renewcommand\cir         {\,{\circ}\,}
\newcommand\Dim           {\mathrm{Dim}}
\newcommand\eps           {\varepsilon}
\newcommand\Hom           {\mathrm{Hom}}
\newcommand\id            {{\rm id}}
\newcommand\In            {\,{\in}\,}

\newcommand\one           {{\bf1}}
\newcommand\oti           {\,{\otimes}\,}
\newcommand\ti            {\,{\times}\,}

\newcommand{\halmos}{\rule{1ex}{1.4ex}}
\newcommand{\pf}{\noindent{\it Proof.}\hspace{2ex}}
\newcommand{\epf}{\hspace*{\fill}\mbox{$\halmos$}}

\newcommand\Cb            {\mathbb{C}}
\newcommand\Nb            {\mathbb{N}}

\newcommand\Zb            {\mathbb{Z}}

\newcommand\Cc            {\mathcal{C}}
\newcommand\Ic            {\mathcal{I}}

\begin{document}
\thispagestyle{empty}
\def\thefootnote{\fnsymbol{footnote}}
\begin{flushright}
KCL-MTH-08-06\\
0807.3356 [math.QA]
\end{flushright}
\vskip 1.0em
\begin{center}\LARGE
Cardy algebras and sewing constraints, I
\end{center}\vskip 1.5em
\begin{center}\large
  Liang Kong%
  $\,{}^{a,b}$\footnote{Email: {\tt kong@mpim-bonn.mpg.de}}
  and
  Ingo Runkel%
  $\,{}^{c}$\footnote{Email: {\tt ingo.runkel@kcl.ac.uk}}%
\end{center}
\begin{center}\it$^a$
Max-Planck-Institut f\"ur Mathematik \\
Vivatsgasse 7, 53111 Bonn, Germany
\end{center}
\begin{center}\it$^b$
Hausdorff Research Institute for Mathematics \\
Poppelsdorfer Allee 45, 53115 Bonn, Germany
\end{center}
\begin{center}\it$^c$
  Department of Mathematics, King's College London \\
  Strand, London WC2R 2LS, United Kingdom  
\end{center}

\vskip 1em

\begin{abstract} 
This is part one of a two-part work that relates two 
different approaches to two-dimensional open-closed rational conformal 
field theory. In part one we review the definition of a Cardy algebra,
which captures the necessary consistency conditions of the
theory at genus 0 and 1. 
We investigate the properties of these algebras
and prove uniqueness and existence theorems.
One implication is that under certain natural assumptions, every rational closed CFT is extendable to an open-closed CFT.
The relation of Cardy algebras to the solutions of 
the sewing constraints is the topic of part two.
\end{abstract}

{\small
\tableofcontents
}

\newpage

\setcounter{footnote}{0}
\def\thefootnote{\arabic{footnote}}

\sect{Introduction and summary}

This is part I of a two-part work which relates two different approaches 
to two-dimensional open-closed rational conformal field theory (CFT). 

\medskip

The first approach uses a three-dimensional topological field theory to express correlators of the  open-closed
CFT \cite{tftcardycase,tft1,unique}. 
Here one starts from a modular tensor category, which defines a 
three-dimensional topological field theory
\cite{retu,turaev-bk}, and from a special symmetric Frobenius 
algebra in this modular tensor category. To each open-closed world
sheet $X$ one assigns a 3-bordism $M_X$ with embedded ribbon graph
constructed from this Frobenius algebra. To the boundary of $M_X$ the 
topological field theory assigns a vector
space $B\ell(X)$ and to $M_X$ itself a vector
$C_X \in B\ell(X)$. One proves that this collection of vectors $C_X$  
provides a so-called solution to the sewing constraints \cite{unique}. 
If the modular tensor category is the category of representations of a 
suitable vertex operator algebra, the spaces $B\ell(X)$ are spaces of 
conformal blocks, and the $C_X$ are the correlators of  
an open-closed CFT. In this approach one thus makes an ansatz for the 
correlators on all world sheets simultaneously and then proves that 
they obey the necessary consistency conditions. The relation to CFT 
rests on convergence and factorisation properties of higher genus 
conformal blocks, and the precise list of conditions the vertex 
operator algebra has to fulfil for these properties to hold is not 
known. However, from a physical perspective one expects that 
interesting classes of models \cite{Witten:1988hf,Frohlich:1989gr} 
will have all the necessary properties.

The second approach uses the theory of vertex operator algebras to construct directly the correlators of the
genus 0 and genus 1 open-closed CFT \cite{osvoa,ffa,cardy}. 
More precisely, in this approach one uses a notion of CFT defined in 
\cite[sect.\,1]{cardy} (and called partial CFT\footnote{
The qualifier `partial' refers to the fact that the gluing of
punctures
is only defined if the coordinates $\zeta_1$, $\zeta_2$ around two punctures
can be analytically extended to a large enough region containing no other punctures, so that the identification $\zeta_1 \sim 1/\zeta_2$ is well-defined. That is, if $\zeta_1$ can be extended to a disc of radius $r$, then $\zeta_2$ must be defined on a disc of radius greater than $1/r$. Both discs must not contain further punctures.
}
there), where one glues Riemann surfaces
around punctures with local coordinates as in \cite{vafa,h-book} 
instead of gluing around parametrised circles as in \cite{Segal}. 
This approach is based on the precise relation between genus-0 CFT and 
vertex operator algebras \cite{h-book}, and on the fact that the category of 
modules over a rational vertex operator algebra is a modular tensor 
category \cite{hl-vtc,h-mtc}. Let us call a vertex operator algebra 
rational if it satisfies the conditions in \cite[sect.\,1]{h-mtc}. If one 
analyses the consistency conditions of a genus-0,1 open-closed
CFT, one arrives at a structure called Cardy $\Cc_V|\Cc_{V\otimes V}$-
algebra in \cite{cardy}. It is formulated in purely categorical terms in 
the categories $\Cc_V$ and $\Cc_{V\otimes V}$ of modules over the rational 
vertex operator algebras $V$ and $V \oti V$, respectively. Cardy algebras 
(defintion \ref{def-I}) are the central objects in 
part I of this work, and we will describe their relation to CFT in slightly 
more detail below. The data in a Cardy algebra amounts to an open-closed
CFT on a generating set of world sheets, from which the entire CFT can be 
obtained by repeated gluing. The conditions on this data are necessary for 
this procedure to give a consistent genus-0,1 open-closed CFT.

The two approaches just outlined start at opposite ends of the same problem. 
In both cases the difficulty to obtain a complete answer lies in the lack of control over the properties of higher genus conformal blocks. 
Nonetheless, both approaches give rise to notions formulated in entirely categorical terms, and we can compare the structures one finds. In part II we will come to the satisfying conclusion that giving a solution to the sewing constraints is essentially equivalent, in a sense made precise in part II, to giving a Cardy algebra.

\medskip

To motivate the notion of a Cardy algebra and our interest in it, we would 
like to outline how it emerges when formulating closed CFT and open-closed 
CFT at genus-0,1 in the language of vertex operator algebras. 
The next one and a half pages, together with a few remarks in the main 
text, are the only places where we make reference to vertex operator 
algebras. The reader who is not familiar with this structure is invited 
to skip ahead.

\medskip

All types of field algebras occurring below are called self-dual if they are endowed with non-degenerate invariant bilinear forms. 

A genus-0 closed CFT is equivalent to an algebra over a partial dioperad consisting of spheres with arbitrary in-coming and out-going punctures. The dioperad structure allows to compose one in-going and one out-going puncture of distinct spheres, so that the result is again a sphere.
Such an algebra with additional natural properties is canonically equivalent to a so-called self-dual conformal full field algebra \cite{ffa,ko-ffa}. A conformal full field algebra contains 
chiral and anti-chiral parts, 
the easiest nontrivial example is given by $V\otimes V$, where $V$ is a vertex operator algebra. A conformal full field algebra containing $V\otimes V$ as a subalgebra is called a conformal full field algebra over $V\otimes V$. 
When $V$ is rational, the category of self-dual conformal full field algebras over $V\otimes V$ is isomorphic to the category of commutative symmetric Frobenius algebras in $\Cc_{V\otimes V}$ 
\cite[thm.\,4.15]{ko-ffa}.

Similarly, a genus-0 open CFT is an algebra over a partial dioperad consisting of disks with an arbitrary number of in-coming and out-going boundary punctures. Such an algebra with additional natural properties is canonically equivalent to a self-dual open-string vertex operator algebra as defined in \cite{osvoa}. A vertex operator algebra $V$ is naturally an open-string vertex operator algebra. An open-string vertex operator algebra containing $V$ as a subalgebra in its meromorphic centre is called open-string vertex operator algebra over $V$. When $V$ is rational, the category of self-dual open-string vertex operator algebras over $V$ is isomorphic to the category of symmetric Frobenius algebras in $\Cc_V$, 
see \cite[thm.\,4.3]{osvoa} and \cite[thm.\,6.10]{cardy}. 

Finally, a genus-0 open-closed CFT is an algebra over the Swiss-cheese partial dioperad, which consists 
of disks with both interior punctures and boundary punctures, and is
equipped with an action of the partial spheres dioperad. 
Such an algebra can be constructed from a so-called self-dual open-closed field algebra \cite{ocfa}. It consists of a self-dual conformal full field algebra $\Acl$, a self-dual open-string vertex operator algebra $\Aop$, and interactions between $\Acl$ and $\Aop$ satisfying certain compatibility conditions. 
Namely, if $\Acl$ is defined over $V\otimes V$ and $\Aop$ over $V$, one requires that the boundary condition on a disc is $V$-invariant in the sense that both the chiral copy $V\otimes \one$ and the anti-chiral copy $\one \otimes V$ of $V$ in $\Acl$ give the copy of $V$ in $\Aop$ in the limit of the insertion point approaching a point on the boundary of the disc \cite[def.\,1.25]{ocfa}.
An open-closed field algebra with $V$-invariant boundary condition is called an open-closed field algebra over $V$. When $V$ is rational, the category of self-dual open-closed field algebras over $V$ is isomorphic to the category of triples $(\Aop|\Acl, \itco)$, 
where $\Acl$ is a commutative symmetric Frobenius $\Cc_{V\otimes V}$-algebra, $\Aop$ a symmetric Frobenius $\Cc_V$-algebra, and $\itco$ an algebra homomorphism $T(\Acl) \rightarrow \Aop$ 
satisfying a centre condition (given in \eqref{eq:comm-C} below),
see \cite[thm.\,3.14]{ocfa} and \cite[sect.\,6.2]{cardy}. 
Here $T: \Cc_{V\otimes V}\rightarrow \Cc_V$ is the Huang-Lepowsky tensor product functor \cite{hl-vtc}. 

The genus-1 theory does not provide new data as it is determined by taking traces of genus-0 correlators, but it does provide two additional consistency conditions: the modular invariance condition for one-point correlators on the torus 
\cite{Sonoda:1988fq}, and the Cardy condition for boundary-two-point correlators on the annulus 
   \cite{Cardy:1989ir,Lew}. 
Their categorical formulations have been worked out in \cite{hk-modinv, cardy}. Adding them to the axioms of a self-dual open-closed field algebra over $V$ finally results in the notion of a Cardy $\Cc_V|\Cc_{V\otimes V}$-algebra. One can prove that the category of self-dual open-closed field algebras over a rational vertex operator algebra $V$ satisfying the two genus-1 consistency conditions is isomorphic to the category of  Cardy $\Cc_V|\Cc_{V\otimes V}$-algebras 
\cite[thm.\,6.15]{cardy}.

If $V$ is rational, then so is $V \oti V$ \cite{DMZ,ffa}. Thus both $\Cc_V$ and $\Cc_{V \otimes V}$ are modular tensor categories. In fact, $\Cc_{V \otimes V} \cong \Cc_V \boxtimes (\Cc_V)_-$ (see \cite[thm.\,4.7.4]{FHL} and \cite[thm\,2.7]{DMZ}), where the minus sign relates to the particular braiding used for $\Cc_{V \otimes V}$. Namely, for a given modular tensor category $\mathcal{D}$, we denote by $\mathcal{D}_-$ the modular tensor category obtained from $\mathcal{D}$ by inverting braiding and twist. We will also sometimes write $\mathcal{D}_+$ for $\mathcal{D}$. The product $\boxtimes$ amounts to taking direct sums of pairs of objects and tensor products of morphisms spaces. The definition of a Cardy algebra can be stated in a way that no longer makes reference to the vertex operator algebra $V$, and therefore makes sense in an arbitrary modular tensor category $\Cc$. Abbreviating $\CcC \equiv \Cc_+ \boxtimes \Cc_-$, this leads to the definition of a Cardy $\Cc|\CcC$-algebra. 

\medskip

The relation to genus-0,1 open-closed CFT outlined above is the main motivation for our interest in Cardy $\Cc|\CcC$-algebras. In part I of this work we investigate how much one can learn about Cardy algebras in the categorical setting, and without the assumption that the modular tensor category $\Cc$ is given by $\Cc_V$ for some $V$. We briefly summarise our approach and results below.

\medskip

In section \ref{sec:colax-fun}--\ref{sec:MTC-def}, we recall some basic notions we will need, such as (co)lax tensor functors, Frobenius functors, and modular tensor categories. In Section \ref{sec:TR-def}, we study the functor $T: \CcC \rightarrow \Cc$, which is defined by the tensor product on $\Cc$ via $T( \oplus_i A_i \times B_i) = \oplus_i A_i \otimes B_i$ for $A_i, B_i \in \Cc$. 
Using the braiding of $\Cc$ one can turn $T$ into a tensor functor.
A tensor functor is automatically also a Frobenius functor, and so takes a Frobenius algebra $A$ in its domain category to a Frobenius algebra $F(A)$ in its target category.

An important object in this work is the functor $R : \Cc\rightarrow \CcC$, also defined in section \ref{sec:TR-def}. We show that $R$ is 
left and right adjoint to $T$. 
As a consequence, $R$ is automatically a lax and colax tensor 
functor, but it is in general not a tensor functor. However, we will 
show that it is still a Frobenius functor, and so takes Frobenius algebras 
in $\Cc$ to Frobenius algebras in $\CcC$. In fact, it also preserves the 
properties simple, special and symmetric of a Frobenius algebra. In the 
case $\Cc=\Cc_V$ the functor 
$R: \mathcal{C}_V\rightarrow \mathcal{C}_{V\otimes V}$ was first 
constructed in \cite{Li-regrep-1,Li-regrep-2} using techniques from vertex 
operator algebras. This motivated the present construction and notation.
    The functor $R$ was also considered in a slightly different context in
    \cite{Etingof:2004}.

The above results imply that $R$ and $T$ form an ambidextrous adjunction, and we will use this adjunction to transport algebraic structures between $\Cc$ and $\CcC$. For example, the algebra homomorphism $\itco : T(\Acl) \rightarrow \Aop$ in $\Cc$ is transported to an algebra homomorphism $\ico : \Acl \rightarrow R(\Aop)$ in $\CcC$.
This gives rise to an alternative definition of a Cardy $\Cc|\CcC$-algebra as 
a triple $(\Aop|\Acl, \ico)$. 

To prepare the definition of a Cardy algebra, in section \ref{sec:mod-inv} we discuss the so-called modular invariance condition for algebras in $\CcC$ (definition \ref{def:mod-inv} below). We show that when $\Acl$ is simple,
the modular invariance condition can be replaced by an easier
condition on the quantum dimension of $\Acl$ (namely, the dimension of $\Acl$ has to be that of the modular tensor category $\Cc$), see theorem \ref{thm:modinv-dim}. 

In section \ref{sec:cardy-def} we give the two definitions of 
a Cardy algebra and prove their equivalence.
Section \ref{sec:cardy-reconstr} contains our main results. We first show 
that for each special symmetric Frobenius algebra $A$ in $\Cc$ 
(see section \ref{sec:algebra} for the definition of special) 
one obtains a Cardy algebra $(A|Z(A),e)$, where $Z(A)$ is the full centre  
of $A$ (theorem \ref{thm:reconst}). The full centre \cite[def.\,4.9]{unique}
is a subobject of $R(A)$ and 
$e : Z(A) \rightarrow R(A)$ is the canonical embedding. Next we prove a uniqueness theorem (theorem \ref{thm:unique}), 
which states that if $(\Aop|\Acl,\ico)$ is a Cardy algebra such that 
$\dim\Aop\neq0$ and $\Acl$ is simple, then $\Aop$ is special and 
$(\Aop|\Acl,\ico)$ is isomorphic to $(\Aop|Z(\Aop), e)$. When combined 
with part II of this work, this result amounts to 
\cite[thm.\,4.26]{unique} and provides an alternative (and shorter) proof.
Finally we show that for every simple
modular invariant commutative symmetric Frobenius algebra $\Acl$ in $\CcC$ there exists a simple special symmetric Frobenius algebra $\Aop$ and an algebra homomorphism $\ico : \Acl \rightarrow R(\Aop)$ such that $(\Aop,\Acl,\ico)$ is a Cardy algebra (theorem \ref{thm:open-exists}). 
This theorem is closely related to a result announced in \cite{mug-conf}
and provides an independent proof in the framework of Cardy algebras.

In physical terms these two theorems mean that a rational open-closed CFT with a unique closed vacuum state can be uniquely reconstructed from its correlators involving only discs with boundary punctures, and that every closed CFT with unique vacuum and left/right rational chiral algebra $V \oti V$ occurs as part of such an open-closed CFT. 

\bigskip
\noindent
{\bf Acknowledgements:}\\
We would like to thank the organisers of the Oberwolfach Arbeitsgemeinschaft ``Algebraic structures in conformal field theories'' (April 2007), where this work was started, for an inspiring meeting. We would further like to thank the Hausdorff Institute for Mathematics in Bonn and the organisers of the stimulating meeting ``Geometry and Physics'' (May 2008). We are indebted to  
Alexei Davydov, Jens Fjelstad, J\"urgen Fuchs, Yi-Zhi Huang, 
Alexei Kitaev, Urs Schreiber, Christoph Schweigert, 
Stephan Stolz and Peter Teichner for helpful discussions and/or comments 
on a draft of this paper. 
The research of IR was partially supported by
the EPSRC First Grant EP/E005047/1, the PPARC rolling grant
PP/C507145/1 and the Marie Curie network `Superstring Theory'
(MRTN-CT-2004-512194).

\newpage

\sect{Preliminaries on tensor categories} \label{sec:prel}

In this section, we review some basic facts of tensor categories 
and fix our conventions and notations along the way. 

\subsection{Tensor categories and (co)lax tensor functors} \label{sec:colax-fun}

In a tensor (or monoidal) category $\Cc$ with 
tensor product bifunctor $\otimes$ and unit object $\one$, 
for $U,V,W \in \Cc$,
we denote the associator $U\otimes (V\otimes W) \xrightarrow{\cong} (U\otimes V)\otimes W$ by $\alpha_{U,V,W}$,  
the left unit isomorphism $\one \otimes U \xrightarrow{\cong} U$ by $l_U$, and
the right unit isomorphism $U\otimes \one \xrightarrow{\cong} U$ by $r_U$. 
If $\Cc$ is braided, for $U,V\in \Cc$ we write the braiding isomorphism as $c_{U,V} : U\otimes V \rightarrow V\otimes U$. 

Let $\Cc_1$ and $\Cc_2$ be two tensor categories with 
units $\one_1$ and $\one_2$ respectively. For simplicity, we will often
write $\otimes, \alpha, l, r$ for the data of 
both $\Cc_1$ and $\Cc_2$. Lax and colax tensor functors
are defined as follows, see e.g.\ \cite[ch.\,I.3]{Yetter} 
or \cite[ch.\,I.1.2]{Leinster}.

\begin{defn}  A {\em lax tensor functor}
$G: \Cc_1 \rightarrow \Cc_2$ is a functor equipped 
with a morphism $\phi_0^G: \one_2 \rightarrow G(\one_1)$ in $\Cc_2$
and a natural transformation 
$\phi_2^G: \otimes \circ (G\ti G) \rightarrow G \circ \otimes$
such that the following three diagrams commute, 
\be
\raisebox{4em}{\xymatrix{
G(A)\otimes (G(B)\otimes G(C)) \ar[r]^{\alpha}
\ar[d]_{\id_{G(A)}\otimes \phi_2^G}  &   
(G(A)\otimes G(B))\otimes G(C) \ar[d]^{\phi_2^G \otimes \id_{G(C)} }  \\
G(A)\otimes G(B\otimes C) \ar[d]_{\phi_2^G} &
G(A\otimes B)\otimes G(C) \ar[d]^{\phi_2^G} \\
G(A\otimes (B\otimes C)) \ar[r]^{G(\alpha)}  & G((A\otimes B)\otimes C)}}   
\quad ,
\label{R-ten-fun-1}
\ee
\be
\raisebox{2em}{\xymatrix{
\one_2\otimes G(A) \ar[d]_{\phi_0^G\otimes \id_{G(A)}}
\ar[r]^{l_{G(A)}} & G(A) \ar[d]^{G(l_A^{-1})} \\
G(\one_1)\otimes G(A) \ar[r]^{\phi_2^G} &  
G(\one_1 \otimes A) }}
\quad , \quad
\raisebox{2em}{\xymatrix{
G(A)\otimes \one_2 \ar[d]_{ \id_{G(A)}\otimes \phi_0^G}
\ar[r]^{r_{G(A)}} & G(A)  \ar[d]^{G(r_A^{-1})} \\
G(A)\otimes G(\one_1) \ar[r]^{\phi_2^G} &
G(A\otimes \one_1)}} \,\, . 
\label{R-ten-fun-2}
\ee
\end{defn}

\begin{defn} 
A {\em colax tensor functor}
is a functor $F: \Cc_1\rightarrow \Cc_2$ 
equipped with a morphism
$\psi_0^F: F(\one_1) \rightarrow \one_2$  in $\Cc_2$, and a natural transformation
$\psi_2^F:  F\circ \otimes \rightarrow  \otimes \circ (F\ti F)$
such that the following three diagrams commute,
\be  \label{L-ten-fun-1}
\raisebox{4em}{\xymatrix{
F(A)\otimes (F(B)\otimes F(C)) \ar[r]^{\alpha} &   
(F(A)\otimes F(B))\otimes F(C)  \\
F(A)\otimes F(B\otimes C) \ar[u]^{\id_{F(A)}\otimes \psi_2^F} &
F(A\otimes B)\otimes F(C)  \ar[u]_{\psi_2^F \otimes \id_{F(C)} }\\
F(A\otimes (B\otimes C)) \ar[r]^{F(\alpha)} \ar[u]^{\psi_2^F} 
& F((A\otimes B)\otimes C)\ar[u]_{\psi_2^F}}}
\quad ,
\ee
\be
\raisebox{2em}{\xymatrix{
\one_2\otimes F(A) 
 & F(A) \ar[l]_{l_{F(A)}^{-1}} \\
F(\one_1)\otimes F(A)  
\ar[u]^{\psi_0^F\otimes \id_{F(A)}}&  
F(\one_1 \otimes A) \ar[l]_{\psi_2^F}\ar[u]_{F(l_A)}}}
\quad , \quad
\raisebox{2em}{\xymatrix{
F(A)\otimes \one_2 
& F(A) \ar[l]_{r_{F(A)}^{-1}}  \\
F(A)\otimes F(\one_1) \ar[u]^{ \id_{F(A)}\otimes \psi_0^F}&
F(A\otimes \one_1) \ar[u]_{F(r_A)} \ar[l]_{\psi_2^F}}}\,\, .
\labl{L-ten-fun-2}
\end{defn}

We denote a lax tensor functor by $(G,\phi_2^G, \phi_0^G)$ 
or just $G$, and a colax tensor functor by $(F, \psi_2^F, \psi_0^F)$  
or $F$.

\begin{defn}
A {\em tensor functor} $T: \Cc_1 \rightarrow \Cc_2$ 
is a lax tensor functor $(T, \phi_2^T, \phi_0^T)$ such that 
$\phi_0^T, \phi_2^T$ are both isomorphisms. 
\end{defn}

A tensor functor $(T, \phi_2^T, \phi_0^T)$ is 
automatically a colax tensor functor $(T, \psi_2^T, \psi_0^T)$ with
$\psi_0^T = (\phi_0^T)^{-1}$ and $\psi_2^T =(\phi_2^T)^{-1}$.

In the next section we will discuss algebras in tensor categories.
The defining properties (\ref{R-ten-fun-1}) and (\ref{R-ten-fun-2})
of a lax tensor functor are analogues of the associativity, the left-unit, 
and the right-unit properties of an algebra. Indeed, a lax tensor functor 
$G:\Cc_1\rightarrow \Cc_2$ maps a $\Cc_1$-algebra to a $\Cc_2$-algebra.  
Similarly, (\ref{L-ten-fun-1}) and (\ref{L-ten-fun-2}) are analogues of the 
coassociativity, the left-counit and the right-counit properties of a 
coalgebra, and a colax tensor functor $F:\Cc_1\rightarrow \Cc_2$ maps a 
$\Cc_1$-coalgebra to a $\Cc_2$-coalgebra. We will later make use of 
functors that take Frobenius algebras to Frobenius algebras. This requires 
  a stronger condition than being lax and colax and leads to the notion of a 
  `functor with Frobenius structure' or `Frobenius monoidal functor'
  \cite{Sz,day-pastro,Pfe}, which we will simply refer to as Frobenius functor.

\begin{defn} \label{def:frob-tensor}
A {\em Frobenius functor} $F: \Cc_1\rightarrow \Cc_2$ is a tuple $F \equiv (F, \phi_2^F, \phi_0^F, \psi_2^F, \psi_0^F)$ such that
$(F,\phi_2^F, \phi_0^F)$ is a lax tensor functor, 
$(F,\psi_2^F, \psi_0^F)$ is a colax tensor functor, and such that
the following two diagrams commute: 
\be
\raisebox{4em}{\xymatrix{
F(A)\otimes (F(B)\otimes F(C)) \ar[r]^{\alpha}  &   
(F(A)\otimes F(B))\otimes F(C) \ar[d]^{\phi_2^F \otimes \id_{F(C)} }   \\
F(A)\otimes F(B\otimes C) \ar[u]^{ \id_{F(A)} \otimes \psi_2^F }  \ar[d]_{\phi_2^F} &
F(A\otimes B)\otimes F(C)  \\
F(A\otimes (B\otimes C)) \ar[r]^{F(\alpha)}  & F((A\otimes B)\otimes C) 
\ar[u]_{\psi_2^F}}}  \label{Frob-ten-cat-1}
\ee
\be
\raisebox{4em}{\xymatrix{
F(A)\otimes (F(B)\otimes F(C)) \ar[d]_{ \id_{F(A)} \otimes \phi_2^F } &   
(F(A)\otimes F(B))\otimes F(C)  \ar[l]_{\alpha^{-1}}   \\
F(A)\otimes F(B\otimes C)   &
F(A\otimes B)\otimes F(C)  \ar[u]_{\psi_2^F \otimes \id_{F(C)} } \ar[d]^{\phi_2^F} \\
F(A\otimes (B\otimes C))   \ar[u]^{\psi_2^F} & F((A\otimes B)\otimes C)  
\ar[l]_{ F(\alpha^{-1}) }}}  \label{Frob-ten-cat-2}
\ee
\end{defn}

\begin{prop}\label{prop:tensor-implies-frobenius}
If $(F, \phi_2^F, \phi_0^F)$ is a tensor functor, then $F$ is a Frobenius functor 
with $\psi_0^F = (\phi_0^F)^{-1}$ and $\psi_2^F =(\phi_2^F)^{-1}$.
\end{prop}
\pf
Since $F$ is a tensor functor, it is lax and colax. If we replace $\psi_2^F$ by $(\phi_2^F)^{-1}$ in \eqref{Frob-ten-cat-1} and \eqref{Frob-ten-cat-2}, both commuting diagrams are equivalent to \eqref{R-ten-fun-1}, which holds because $F$ is lax. Thus $F$ is a Frobenius functor.
\epf

\medskip

The converse statement does not
hold. For example, the functor $R$ which we define in section 
\ref{sec:TR-def} is Frobenius but not tensor.

\medskip

Let us recall the notion of adjunctions 
and adjoint functors \cite[ch.\,IV.1]{Ma}. 

\begin{defn}
An {\em adjunction} from $\Cc_1$ to $\Cc_2$ is
a triple $\langle F, G, \chi\rangle$, where $F$ and $G$ are 
functors
$$
F: \Cc_1 \rightarrow \Cc_2~~, \quad G: \Cc_2 \rightarrow 
\Cc_1~, 
$$
and $\chi$ is a natural isomorphism which assigns to each pair of objects
$A_1 \in \Cc_1$, $A_2 \in \Cc_2$ a bijective map
$$
\chi_{A_1, A_2}:  \Hom_{\Cc_2}(F(A_1), A_2) 
\xrightarrow{~~\cong~~} 
\Hom_{\Cc_1}(A_1, G(A_2))~,
$$
which is natural in both $A_1$ and $A_2$. 
$F$ is called a {\em left-adjoint} of $G$ and $G$  is called 
a {\em right-adjoint} of $F$. 
\end{defn}

For simplicity, we will often abbreviate $\chi_{A_1, A_2}$ as $\chi$. 
Associated to each 
adjunction $\langle F, G, \chi\rangle$, there are two natural transformations 
$\id_{\Cc_1} \xrightarrow{\delta} GF$
and $FG \xrightarrow{\rho} \id_{\Cc_2}$, where
$\id_{\Cc_1}$ and $\id_{\Cc_2}$ are identity functors, 
given by 
\be
\delta_{A_1}= \chi(\id_{F(A_1)}), \quad\quad \rho_{A_2}= \chi^{-1}(\id_{G(A_2)})
\ee
for $A_i\in \Cc_i, i=1,2$. 
They satisfy the following two identities: 
\be  \label{eq:adj-unit-counit}
G \xrightarrow{\delta G} GFG \xrightarrow{G\rho} G
= G \xrightarrow{\id_G} G, 
\quad\quad 
F \xrightarrow{F\delta} FGF \xrightarrow{\rho F} F
= F \xrightarrow{\id_F} F~.
\ee
We have, 
for $g: F(A_1) \rightarrow A_2$ and $f: A_1\rightarrow G(A_2)$, 
\be  \label{eq:chi-eta-epsilon}
\chi(g) = G(g) \circ \delta_{A_1},  \quad\quad 
\chi^{-1}(f) =  \rho_{A_2} \circ F(f)~.  
\ee
For simplicity, $\delta_{A_1}$ and $\rho_{A_2}$ are often 
abbreviated as $\delta$ and $\rho$, respectively.  

Let $\langle F, G, \chi\rangle$ be an adjunction from a tensor category 
$\Cc_1$ to a tensor category $\Cc_2$
and $(F, \psi_2^F, \psi_0^F)$ a colax tensor functor from $\Cc_1$ to $\Cc_2$. 
We can define a morphism $\phi_0^G: \one_1 \rightarrow G(\one_2)$ 
and a natural transformation 
$\phi_2^G : \otimes \circ (G \ti G) \rightarrow G\circ \otimes$ by,
for $A,B\in \Cc_2$,
\be\begin{array}{l}\displaystyle
\phi_0^G =\chi(\psi_0^F) = \,\,
\one_1 \xrightarrow{\delta_{\one_1}} GF(\one_1) \xrightarrow{G(\psi_0^F)}
G(\one_2),  \enl
\phi_2^G = \chi( (\rho_A\otimes \rho_B) \circ \psi_2^F) = \,\,
 G(A)\otimes G(B)\xrightarrow{\delta} GF\big(G(A)\otimes G(B)\big)  
 \\ \displaystyle
\hspace{4cm}
\xrightarrow{G(\psi_2^F)} G\big(FG(A)\otimes FG(B)\big)
\xrightarrow{G(\rho_A\otimes \rho_B)} G(A\otimes B). 
\end{array}
\labl{eq:phi2F-def-psi}
where we have used the first identity in (\ref{eq:chi-eta-epsilon}). 
Notice that $\phi_2^G$ is natural because it is a composition of 
natural transformations. One can easily show that $\psi_0^F$ and
$\psi_2^F$ can be re-obtained from $\phi_0^G$ and $\phi_2^G$ as follows:
\be\begin{array}{l}\displaystyle
\psi_0^F = \chi^{-1}(\phi_0^G) = \,\, F(\one_1) \xrightarrow{F\phi_0^G} FG(\one_2) 
\xrightarrow{\rho} \one_2   \enl
\psi_2^F = \chi^{-1}(\phi_2^G \circ (\delta\otimes \delta)) = \,\,
F(U\otimes V) \xrightarrow{F(\delta\otimes \delta)} F\big(GF(U)\otimes GF(V)\big)   \\ \displaystyle
\hspace{4cm}
\xrightarrow{F\phi_2^G} FG\big(F(U) \otimes F(V)\big)  \xrightarrow{\rho} F(U)\otimes F(V). 
\end{array}
\labl{eq:psi-by-phi}
for $U,V \in \Cc_1$. 
The following result is standard; 
for the sake of completeness, we give a proof in appendix \ref{app:lax-colax}. 

\begin{lemma}  \label{lem:lax-colax-adj}
$(F, \psi_2^F, \psi_0^F)$ is a colax tensor functor iff 
$(G, \phi_2^G, \phi_0^G)$ is a lax tensor functor. 
\end{lemma}

\subsection{Algebras in tensor categories}\label{sec:algebra}

An algebra in a tensor category $\Cc$, or a $\Cc$-algebra, is a triple $A=(A,m,\eta)$ where $A$ is an object of $\Cc$, $m$ (the multiplication) is a morphism $A \oti A \rightarrow A$ such that $m \cir (m \oti \id_A) \cir \alpha_{A,A,A} = m \cir (\id_A \oti m)$, and $\eta$ (the unit) is a morphism $\one \rightarrow A$ such that $m \cir (\id_A \oti \eta) = \id_A  \cir r_A$ and $m \cir (\eta\oti \id_A) = \id_A \cir l_A$. If $\Cc$ is braided and $m \cir c_{A,A} = m$, then $A$ is called commutative.  

A left $A$-module is a pair $(M, m_M)$, where $M\in \Cc$ 
and $m_M$ is a morphism $A\otimes M \rightarrow M$ such that 
$m_M\circ (\id_A\otimes m_M) = m_M \circ (m_A \otimes \id_M)\circ \alpha_{A,A,M}$ 
and $m_M\circ (\eta_A \otimes \id_M) = \id_M \circ l_M$. 
Right $A$-modules and $A$-bimodules are defined similarly. 

\begin{defn}\label{def:simp-asimp-hapl}
Let $\Cc$ be a tensor category and let $A$ be an algebra in $\Cc$.
\\[.3em]
(i) $A$ is called {\em simple} iff it is simple as a bimodule over itself.    
\\[.3em]
Let $\Cc$ be in addition $\Bbbk$-linear, for $\Bbbk$ a field.
\\[.3em]
(ii) $A$ is called {\em absolutely simple} iff the space of $A$-bimodule maps from $A$ to itself is one-dimensional, $\dim_{\Bbbk} \Hom_{A|A}(A,A)=1$. 
\\
(iii) $A$ is called {\em haploid} iff
$\dim_{\Bbbk} \Hom(\one,A)=1$ {\rm \cite[def.\,4.3]{fs-cat}}.
\end{defn}
In the following we will assume that all tensor categories are strict to avoid spelling out associators and unit constraints. 

\medskip

A $\Cc$-coalgebra $A = (A,\Delta,\eps)$ is defined analogously to a $\Cc$-algebra, i.e.\ $\Delta : A \rightarrow A \oti A$ and $\eps : A \rightarrow \one$ obey coassociativity and counit conditions. 

If $\Cc$ is braided and if $A$ and $B$ are $\Cc$-algebras, there are two in general non-isomorphic algebra structures on $A\otimes B$. We choose $A \oti B$ to be the $\Cc$-algebra with multiplication $m_{A\otimes B}= (m_A\otimes m_B) \circ (\id_A \otimes c_{A,B}^{-1} \otimes \id_B)$ and unit $\eta_{A\otimes B}=\eta_A\otimes \eta_B$. 
Similarly, if $A$ and $B$ are $\Cc$-coalgebras, then $A \oti B$
becomes a $\Cc$-coalgebra if we choose
the comultiplication $\Delta_{A\otimes B}= (\id_A\otimes c_{A,B} \otimes \id_B) \circ (\Delta_A \otimes \Delta_B)$ and 
the counit $\eps_{A\otimes B}= \eps_A\otimes \eps_B$. 

\begin{defn}
A {\em Frobenius algebra} $A = (A,m,\eta,\Delta,\eps)$ is an algebra and a coalgebra such that the coproduct is an intertwiner of $A$-bimodules,
$$(\id_A \oti m) \cir (\Delta \oti \id_A) = \Delta \oti m = (m \oti \id_A) \cir (\id_A \oti \Delta)~.$$ 
\end{defn}

We will use the following graphical representation for the morphisms of a Frobenius algebra,
\be
  m = \raisebox{-20pt}{
  \begin{picture}(30,45)
   \put(0,6){\scalebox{.75}{\includegraphics{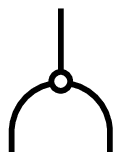}}}
   \put(0,6){
     \setlength{\unitlength}{.75pt}\put(-146,-155){
     \put(143,145)  {\scriptsize $ A $}
     \put(169,145)  {\scriptsize $ A $}
     \put(157,202)  {\scriptsize $ A $}
     }\setlength{\unitlength}{1pt}}
  \end{picture}}  
  ~~,\quad
  \eta = \raisebox{-15pt}{
  \begin{picture}(10,30)
   \put(0,6){\scalebox{.75}{\includegraphics{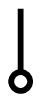}}}
   \put(0,6){
     \setlength{\unitlength}{.75pt}\put(-146,-155){
     \put(146,185)  {\scriptsize $ A $}
     }\setlength{\unitlength}{1pt}}
  \end{picture}}
  ~~,\quad
  \Delta = \raisebox{-20pt}{
  \begin{picture}(30,45)
   \put(0,6){\scalebox{.75}{\includegraphics{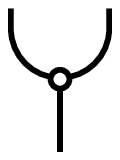}}}
   \put(0,6){
     \setlength{\unitlength}{.75pt}\put(-146,-155){
     \put(143,202)  {\scriptsize $ A $}
     \put(169,202)  {\scriptsize $ A $}
     \put(157,145)  {\scriptsize $ A $}
     }\setlength{\unitlength}{1pt}}
  \end{picture}}
  ~~,\quad
  \eps = \raisebox{-15pt}{
  \begin{picture}(10,30)
   \put(0,10){\scalebox{.75}{\includegraphics{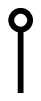}}}
   \put(0,10){
     \setlength{\unitlength}{.75pt}\put(-146,-155){
     \put(146,145)  {\scriptsize $ A $}
     }\setlength{\unitlength}{1pt}}
  \end{picture}}
  ~~.
\ee
A Frobenius algebra $A$ in a $\Bbbk$-linear tensor category, for $\Bbbk$ a field, is called {\em special} iff $m \circ \Delta = \zeta \, \id_A$ and $\eps \circ \eta = \xi \, \id_\one$ for nonzero constants $\zeta$, $\xi \In \Bbbk$. If $\zeta=1$ we call $A$ {\em normalised-special}. A Frobenius algebra homomorphism between two 
Frobenius algebras is both an algebra homomorphism and 
a coalgebra homomorphism.

A (strictly) sovereign tensor category is a tensor category equipped with a left and a right duality which agree on objects and morphisms (see e.g.\ \cite{bichon,fs-cat} for more details). We will write the dualities as
\be
\begin{array}{llll}
  \raisebox{-8pt}{
  \begin{picture}(26,22)
   \put(0,6){\scalebox{.75}{\includegraphics{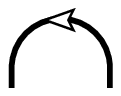}}}
   \put(0,6){
     \setlength{\unitlength}{.75pt}\put(-146,-155){
     \put(143,145)  {\scriptsize $ U^\vee $}
     \put(169,145)  {\scriptsize $ U $}
     }\setlength{\unitlength}{1pt}}
  \end{picture}}  
  \etb= d_U : U^\vee \oti U \rightarrow \one
  ~~,\qquad &
  \raisebox{-8pt}{
  \begin{picture}(26,22)
   \put(0,6){\scalebox{.75}{\includegraphics{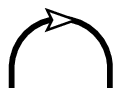}}}
   \put(0,6){
     \setlength{\unitlength}{.75pt}\put(-146,-155){
     \put(143,145)  {\scriptsize $ U $}
     \put(169,145)  {\scriptsize $ U^\vee $}
     }\setlength{\unitlength}{1pt}}
  \end{picture}}  
  \etb= \tilde d_U : U \oti U^\vee \rightarrow \one
  ~~,
\\[2em]
  \raisebox{-8pt}{
  \begin{picture}(26,22)
   \put(0,0){\scalebox{.75}{\includegraphics{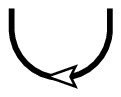}}}
   \put(0,0){
     \setlength{\unitlength}{.75pt}\put(-146,-155){
     \put(143,183)  {\scriptsize $ U $}
     \put(169,183)  {\scriptsize $ U^\vee $}
     }\setlength{\unitlength}{1pt}}
  \end{picture}}  
  \etb= b_U : \one \rightarrow U \oti U^\vee
  ~~,
  &
  \raisebox{-8pt}{
  \begin{picture}(26,22)
   \put(0,0){\scalebox{.75}{\includegraphics{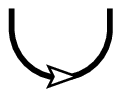}}}
   \put(0,0){
     \setlength{\unitlength}{.75pt}\put(-146,-155){
     \put(143,183)  {\scriptsize $ U^\vee $}
     \put(169,183)  {\scriptsize $ U $}
     }\setlength{\unitlength}{1pt}}
  \end{picture}}  
  \etb= \tilde b_U : \one \rightarrow U^\vee \oti U
  ~.
\eear\labl{eq:duality-mor}
In terms of these we define the left and right dimension of an
object $U$ as
\be
  \dim_l U = d_U \circ \tilde b_U \qquad , \quad
  \dim_r U = \tilde d_U \circ b_U ~,
\labl{eq:dim-def}
both of which are elements of $\Hom(\one,\one)$.

Let now $\Cc$ be a sovereign tensor category. For a Frobenius algebra $A$ in $\Cc$, we define two morphisms: 
\be
\Phi_A =\,\,
\raisebox{-35pt}{
  \begin{picture}(50,75)
   \put(0,8){\scalebox{.75}{\includegraphics{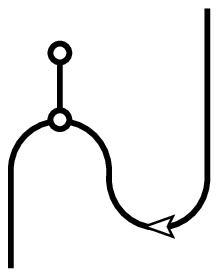}}}
   \put(0,8){
     \setlength{\unitlength}{.75pt}\put(-34,-37){
     \put(31, 28)  {\scriptsize $ A $}
     \put(87,117)  {\scriptsize $ A^\vee $}
     }\setlength{\unitlength}{1pt}}
  \end{picture}},  \quad\quad\quad
\Phi_A' =\,\,
  \raisebox{-35pt}{
  \begin{picture}(50,75)
   \put(0,8){\scalebox{.75}{\includegraphics{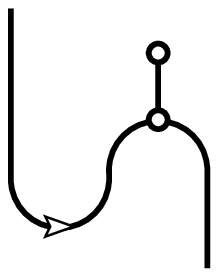}}}
   \put(0,8){
     \setlength{\unitlength}{.75pt}\put(-34,-37){
     \put(87, 28)  {\scriptsize $ A $}
     \put(31,117)  {\scriptsize $ A^\vee $}
     }\setlength{\unitlength}{1pt}}
  \end{picture}}
  ~~.
\labl{eq:Frob-sym-cond}

\begin{defn}
A Frobenius algebra $A$ is {\em symmetric} iff $\Phi_A = \Phi_A'$. 
\end{defn}

The following lemma shows that under certain conditions we do not need to distinguish the various notions of simplicity in definition \ref{def:simp-asimp-hapl}.

\begin{lemma}\label{lem:modinv-simp-hapl}
Let $A$ be a commutative symmetric Frobenius algebra in a 
$\mathbb{C}$-linear semi-simple sovereign braided tensor category $\Cc$ 
and suppose that $\dim_l A \neq 0$. Then the following are equivalent.
\\
(i)\phantom{ii} $A$ is simple.
\\
(ii)\phantom{i} $A$ is absolutely simple.
\\
(iii) $A$ is haploid.
\end{lemma}

\pf
(ii)$\Leftrightarrow$(iii): 
$A$ is haploid iff it is absolutely simple as a left module over itself \cite[eqn.\,(4.17)]{fs-cat}. Furthermore, for a commutative algebra we have $\Hom_{A}(A,A) = \Hom_{A|A}(A,A)$, and so $A$ is haploid iff it is absolutely simple.
\\[.3em]
(i)$\Rightarrow$(ii): If $A$ is simple, then
every nonzero element of $\Hom_{A|A}(A,A)$ is
invertible. Hence this space forms a division algebra over $\Cb$,
and is therefore isomorphic to $\Cb$.
\\[.3em]
(iii)$\Rightarrow$(i): 
Since $\Cc$ is semi-simple and $A$ is haploid, also 
   $\Hom(A,\one)$ 
is one-dimensional. 
The counit $\eps$
is a nonzero element in this space, and so gives a basis. 
This implies firstly, 
that $\eps \circ \eta \neq 0$, and
secondly, that
there is a constant $\beta \in \Cb$ such that
\be
  \beta \cdot \eps
  = d_{A} \circ 
  (\id_{A^\vee} \otimes m) \circ 
  (\tilde b_{A}\otimes \id_{A}) ~.
\ee
Composing with $\eta$ from the right yields
$\beta \, \eps \circ \eta = \dim_l A$.
The right hand side is nonzero, and so $\beta \neq 0$. 
By \cite[lem.\,3.11]{tft1}, $A$ is special.
We have already proved (ii)$\Leftrightarrow$(iii), and so
$A$ is absolutely simple.
A special
Frobenius algebra in a semi-simple category has a semi-simple
category of bimodules
(apply \cite[prop.\,5.24]{fs-cat} to the algebra 
tensored with its opposite algebra).
For semi-simple $\Cb$-linear
categories, simple and absolutely
simple are equivalent\footnote{
  \label{fn:HomUU=C}
  To see this note that if $U$ is simple, then
  the $\Cb$-vector space 
  $\Hom(U,U)$ is a division algebra, and hence 
  $\Hom(U,U) = \Cb \, \id_U$.
  Conversely, if $U$ is not simple, then $U = U_1 \oplus U_2$ and
  $\Hom(U,U)$ contains at least two linearly independent elements,
  namely $\id_{U_1}$ and $\id_{U_2}$.}.
Thus $A$ is simple.
\epf

\begin{rema}{\rm
For a Frobenius algebra $A$ the morphisms \eqref{eq:Frob-sym-cond}
are invertible, and hence $A \cong A^\vee$. In this case
one has $\dim_l A = \dim_r A$  \cite[rem.\,3.6.3]{fs-cat}
and so we could have stated the above lemma equivalently
with the condition $\dim_r A\neq 0$.
}
\end{rema}

Let $F: \Cc_1\rightarrow \Cc_2$ be 
a lax tensor functor between two tensor categories $\Cc_1,\Cc_2$ and
let $(A, m_A, \eta_A)$ be an algebra 
in $\Cc_1$. Define morphisms 
$F(A)\otimes F(A)\xrightarrow{m_{F(A)}} F(A)$ and
$\one_2 \xrightarrow{\eta_{F(A)}} F(A)$ as  
\be  
m_{F(A)} = F(m_A) \circ \phi_2^F, \quad\quad 
\eta_{F(A)}= F(\eta_A) \circ \phi_0^F ~.
\labl{eq:FA-alg-st}
Then $(F(A), m_{F(A)}, \eta_{F(A)})$ is an algebra in $\Cc_2$ 
\cite[prop.\,5.5]{JS}.
If $f: A\rightarrow B$ is an algebra homomorphism between two algebras $A, B\in \Cc_1$, then $F(f): F(A) \rightarrow F(B)$ is also an algebra homomorphism. If $(M, m_M)$ is a left (or right) $A$-module in $\Cc_1$, then $(F(M), F(m_M)\circ \phi_2^F)$ is a left (or right) $F(A)$-module; 
if $M$ has a $A$-bimodule structure, 
then $F(M)$ naturally has a $F(A)$-bimodule structure. 

Similarly, if $(A, \Delta_A, \eps_A)$ is a  
coalgebra in $\Cc_1$ and $F: \Cc_1\rightarrow \Cc_2$ is 
a colax tensor functor, then $F(A)$ with coproduct $F(A) \xrightarrow{\Delta_{F(A)}} F(A)\otimes F(A)$ and counit 
$F(A)\xrightarrow{\eps_{F(A)}} \one_2$ given by 
\be  \label{eq:FA-coalg-st}  
\Delta_{F(A)} = \psi_2^F \circ F(\Delta_A), \qquad 
\eps_{F(A)} = \psi_0^F \circ F(\eps_A),
\ee
is a coalgebra in $\Cc_2$. 
If $f: A\rightarrow B$ is a coalgebra homomorphism between two coalgebras $A, B\in \Cc_1$, then $F(f): F(A) \rightarrow F(B)$ is also a coalgebra homomorphism.

\begin{prop} 
\footnote{After the preprint of the present paper appeared we 
noticed that this proposition is also proved in 
\cite[cor.\,5]{day-pastro}.}
\label{prop:t-frob}
If $F: \, \Cc_1\rightarrow \Cc_2$ is a Frobenius functor and $(A, m_A, \eta_A, \Delta_A, \eps_A)$ a Frobenius algebra in $\Cc_1$, then $(F(A), m_{F(A)}, \eta_{F(A)}, \Delta_{F(A)}, \eps_{F(A)})$ is a Frobenius algebra in $\Cc_2$. 
\end{prop}
\pf
One Frobenius property,
$(m_{F(A)} \otimes \id_{F(A)}) \circ (\id_{F(A)} \otimes \Delta_{F(A)}) 
= \Delta_{F(A)}\circ m_{F(A)}$, follows 
from the commutativity of the following diagram (we spell out the associativity
isomorphisms): 
\be
\raisebox{5em}{  
\small{ 
\xymatrix{
F(A)\otimes F(A) \ar[rr]^{\hspace{-0.3cm}\id_{F(A)}\otimes F(\Delta_A)}
\ar[d]^{\phi_2^F} & & F(A) \otimes F(A\otimes A) 
\ar[rr]^{\hspace{-0.3cm}\id_{F(A)} \otimes \psi_2^F}\ar[d]^{\phi_2^F} 
& & F(A)\otimes (F(A)\otimes F(A)) \ar[d]^{\alpha_{F(A), F(A), F(A)}}   \\
F(A\otimes A) \ar[rr]^{F(\id_A\otimes \Delta_A)}\ar[dd]^{F(m_A)} 
& & F(A\otimes (A\otimes A)) \ar[d]^{F(\alpha_{A,A,A})}
& &(F(A)\otimes F(A))\otimes F(A) \ar[d]^{\phi_2^F \otimes \id_{F(A)}}\\
 & & F((A\otimes A) \otimes A) \ar[d]^{F(m_A\otimes \id_A)}\ar[rr]^{\psi_2^F} 
& & F(A\otimes A) \otimes F(A) \ar[d]^{F(m_A)\otimes \id_{F(A)}} \\
F(A) \ar[rr]^{F(\Delta_A)} & & F(A\otimes A) \ar[rr]^{\psi_2^F} 
& & F(A)\otimes F(A) 
}}}
\ee
The commutativity of the upper-left subdiagram follows from the 
naturalness of $\phi_2^F$, that of the upper-right subdiagram 
follows from (\ref{Frob-ten-cat-1}), that of the lower-left 
subdiagram follows from the Frobenius properties of $A$,
and that of the lower-right subdiagram follows from 
the naturalness of $\psi_2^F$. The proof of the other Frobenius
property is similar. 
\epf

\begin{prop} \label{prop:T-sp-Frob}
If $F: \Cc_1 \rightarrow \Cc_2$ is a tensor functor and $A$ a Frobenius algebra in $\Cc_1$, then:  \\
(i) $F(A)$ has a natural structure of Frobenius algebra as given in proposition \ref{prop:t-frob};  \\
(ii) If $A$ is (normalised-)special, so is $F(A)$. 
\end{prop}
\pf
Part (i) follows from propositions \ref{prop:tensor-implies-frobenius} and \ref{prop:t-frob}. Part (ii) is a straightforward verification of the definition, using 
$\psi^F_2 = (\phi^F_2)^{-1}$ and $\psi^F_0 = (\phi^F_0)^{-1}$.
\epf

\medskip

Let $\Cc_1, \Cc_2$ be sovereign tensor categories and $F: \Cc_1\rightarrow \Cc_2$ a Frobenius functor. We define two morphisms
 $I_{F(A^{\vee})}, I_{F(A^{\vee})}':  F(A^{\vee}) \rightarrow F(A)^{\vee}$, 
for a Frobenius algebra $A$ in $\Cc_1$, as follows:  
\be\begin{array}{l}\displaystyle
I_{F(A^{\vee})} =
((\psi_0^F\circ F(d_A) \circ \phi_2^F) 
\otimes \id_{F(A)^{\vee}}) \circ
(\id_{F(A^{\vee})} \otimes b_{F(A)})
~,\enl
I_{F(A^{\vee})}' =
(\id_{F(A)^{\vee}} \otimes (\psi_0^F\circ F(\tilde{d}_A) \circ \phi_2^F))
\circ (\tilde{b}_{F(A)} \otimes \id_{F(A^{\vee})}) 
~.  
\end{array}
\labl{equ:II'-def}
It is easy to see that these are isomorphisms. 

\begin{lemma} \label{lem:Phi-T(A)}
If $F:\Cc_1\rightarrow \Cc_2$ is a Frobenius functor and $A$ a Frobenius algebra in $\Cc_1$, then
\be \label{eq:Phi-I}
\Phi_{F(A)}  = I_{F(A^{\vee})} \circ F(\Phi_A) \quad , \quad
\Phi_{F(A)}' = I_{F(A^{\vee})}' \circ F(\Phi_A').
\ee
\end{lemma}
\pf
We only prove the first equality, the second one can be seen in 
the same way. By definition, we have
\bea
&&I_{F(A^{\vee})} \circ F(\Phi_A)  \nn
&& \hspace{0.2cm}= 
((\psi_0^F \circ F(d_A) \circ \phi_2^F) 
\otimes \id_{F(A)^{\vee}}) \circ
(\id_{F(A^{\vee})} \otimes b_{F(A)}) \circ 
F(\Phi_A) \nn
&& \hspace{0.2cm}= 
\big\{ \big[ (\psi_0^F\circ F(d_A) \circ \phi_2^F) \circ (F(\Phi_A) 
\otimes \id_{F(A)}) \big] \otimes \id_{F(A)^{\vee}} \big\} \circ (\id_{F(A)} 
\otimes b_{F(A)}) ~.  \nonumber
\eea
For the term inside the square brackets we find
\be\begin{array}{l}\displaystyle
\psi_0^F\circ F(d_A) \circ \phi_2^F \circ (F(\Phi_A) 
\otimes \id_{F(A)} ) =
\psi_0^F \circ F(d_A) \circ 
F(\Phi_A\otimes \id_A)\circ \phi_2^F \enl
= \psi_0^F \circ F(d_A \circ (\Phi_A\otimes \id_A)) \circ \phi_2^F 
= \psi_0^F \circ F(\eps_A \circ m_A) \circ \phi_2^F~.
\end{array}
\ee
On the other hand, by definition, $\Phi_{F(A)} = 
[((\psi_0^F \circ F(\eps_A)  \circ (F(m_A) \circ \phi_2^F)) 
\otimes \id_{F(A)^{\vee}}] \circ 
(\id_{F(A)} \otimes b_{F(A)})$. 
This demonstrates the first equality in (\ref{eq:Phi-I}).
\epf

\begin{prop}  \label{prop:f-til-f}
Let $F: \Cc_1 \rightarrow \Cc_2$ be a tensor functor, 
$G : \Cc_2 \rightarrow \Cc_1$ a functor,
$\langle F,G,\chi \rangle$ 
an adjunction, $A$ a $\Cc_1$-algebra, and $B$ a $\Cc_2$-algebra. 
Then $f: A\rightarrow G(B)$ is an algebra homomorphism if and only if 
$\tilde{f}= \chi^{-1}(f):  F(A) \rightarrow B$ is an algebra homomorphism. 
\end{prop}
\pf
We need to show that
\be
m_{G(B)} \circ (f \otimes f)  =  f \circ m_A 
\qquad \text{and} \qquad
f \circ \eta_A   = \eta_{G(B)},
\labl{eq:f-alg-map+unit}
is equivalent to 
\be  
m_B \circ (\tilde{f} \otimes \tilde{f})  =  \tilde{f} \circ m_{F(A)}
\qquad \text{and} \qquad
\tilde{f} \circ \eta_{F(A)} = \eta_B ~.
\labl{eq:f-til-alg-map+unit}
We first prove that the first identity in 
\eqref{eq:f-alg-map+unit} is equivalent to
the first identity in  \eqref{eq:f-til-alg-map+unit}.
For the left hand side of the first identity in \eqref{eq:f-alg-map+unit}
we have the following equalities,
\bea
m_{G(B)} \circ (f\otimes f) 
&\overset{(1)}{=}& 
G(m_B) \circ \phi_2^G \circ (f\otimes f)  \nn
&\overset{(2)}{=}& 
G(m_B) \circ G(\rho \otimes \rho) \circ G(\psi_2^F) \circ \delta \circ (f\otimes f) \nn
&\overset{(3)}{=}& 
G(m_B) \circ G(\rho \otimes \rho) \circ G(\psi_2^F) \circ GF(f\otimes f) \circ \delta
\nn
&\overset{(4)}{=}& 
G(m_B) \circ G(\rho \otimes \rho) \circ G(F(f)\otimes F(f)) \circ G(\psi_2^F) \circ \delta  \nn
&\overset{(5)}{=}& 
G\big( m_B \circ (\rho \otimes \rho) \circ (F(f) \otimes F(f)) \circ \psi_2^F \big) \circ \delta \nn
&\overset{(6)}{=}& 
\chi \big( m_B \circ (\rho \otimes \rho) \circ (F(f) \otimes F(f)) \circ \psi_2^F 
\big)~,   \label{eq:m-FB-ff}
\eea
where (1) is the definition of $m_{G(B)}$ in \eqref{eq:FA-alg-st}, (2) is the second identity in \eqref{eq:phi2F-def-psi}, (3) and (4) are naturality of $\delta$ and $\psi_2^F$, respectively, step (5) is functoriality of $G$ and finally step (6) is \eqref{eq:chi-eta-epsilon}.
For the right hand side of the first identity in \eqref{eq:f-alg-map+unit}
we get
\bea
f\circ m_A 
&\overset{(1)}{=}& 
G\rho \circ \delta G \circ (f \circ m_A)  \nn
&\overset{(2)}{=}& 
G\rho \circ GF(f \circ m_A)  \circ \delta   \nn
&\overset{(3)}{=}& 
G( \rho \circ F(f \circ m_A)) \circ \delta \nn
&\overset{(4)}{=}& 
\chi( \rho \circ F(f) \circ F(m_A))~,   \label{eq:f-mA}
\eea
where (1) is the adjunction property \eqref{eq:adj-unit-counit}, (2) is naturality of $\delta$, (3) functoriality of $G$, and (4) amounts to \eqref{eq:chi-eta-epsilon} and functoriality of $F$.

On the other hand, we see that the first equality in 
\eqref{eq:f-til-alg-map+unit}
is equivalent to 
\be 
m_B \circ ( \rho \otimes \rho) \circ \big( F(f) \otimes F(f) \big) 
= \rho \circ F(f) \circ F(m_A) \circ \phi_2^F~.
\labl{eq:mff=fm}
Using that $\phi_2^F$ is invertible with inverse $(\phi_2^F)^{-1} = \psi_2^F$ and that $\chi$ is an isomorphism, it follows that the statement that \erf{eq:m-FB-ff} is equal to \erf{eq:f-mA} is equivalent to the identity \erf{eq:mff=fm}.

Now we prove that the second identity in 
\eqref{eq:f-alg-map+unit} is equivalent to
the second identity in  \eqref{eq:f-til-alg-map+unit}.
Using \eqref{eq:FA-alg-st} and \eqref{eq:phi2F-def-psi} we can write
$\eta_{G(B)} = G(\eta_B) \circ \phi_0^F =
G(\eta_B) \circ G(\psi_0^F) \circ \delta_\one$. 
Together with \eqref{eq:chi-eta-epsilon} this shows that the second identity in 
\eqref{eq:f-alg-map+unit} is equivalent to
\be  \label{eq:f-unit-1}
f\circ \eta_A = \chi( \eta_B \circ \psi_0^F)~. 
\ee
On the other hand, the second identity in
\eqref{eq:f-til-alg-map+unit} is equivalent to   
\be
\rho \circ F(f) \circ F(\eta_A) \circ \phi_0^F = \eta_B~, 
\ee
which, by $\phi_0^F = (\psi_0^F)^{-1}$ and (\ref{eq:chi-eta-epsilon}), 
is further equivalent to (\ref{eq:f-unit-1}). 
\epf

\medskip

\begin{defn}  {\rm
Let $(A, m_A, \eta_A, \Delta_A, \eps_A)$ 
and $(B, m_B, \eta_B, \Delta_B, \eps_B)$ 
be two Frobenius algebras in a 
tensor category $\mathcal{C}$. 
For $f: A\rightarrow B$, we define $f^*: B\rightarrow A$ 
by 
\be  \label{f-star-def}
f^* =  ( (\eps_B \circ m_B) \otimes \id_A) 
\circ (\id_B \otimes f \otimes \id_A) 
\circ (\id_B \otimes (\Delta_A \circ \eta_A)).
\ee
}
\end{defn}

The following lemma is immediate from the definition
of $(\cdot)^*$ and the properties of Frobenius algebras.
We omit the proof.

\begin{lemma}\label{lem:f-star-prop}
Let $\Cc$ be a tensor category, 
let $A,B,C$ be Frobenius algebras in $\Cc$, and
let $f : A \rightarrow B$ and $g : B \rightarrow C$ be morphisms.
\\
(i)\phantom{ii}~$(g \circ f)^* = f^* \circ g^*$. \\
(ii)\phantom{i}~$f$ is a monomorphism iff $f^*$ is an epimorphism. \\
(iii)~$f$ is an algebra map iff $f^*$ is a coalgebra map. \\
(iv)\,~If $f$ is a homomorphism of Frobenius algebras, then 
$f^* \circ f = \id_A$ and $f \circ f^* = \id_B$.\\
(v)\phantom{i}\,~If $\Cc$ is sovereign and if $A$ and $B$ 
are symmetric, then $f^{**} = f$. \\
\end{lemma}

Let $\Cc$ and $\mathcal{D}$ be tensor categories and let
$F : \Cc \rightarrow \mathcal{D}$ be a Frobenius functor.
Given Frobenius algebras $A,B$ in $\Cc$ and a morphism
$f : A \rightarrow B$,
the next lemma shows how $(\cdot)^*$ behaves under $F$.

\begin{lemma}\label{lem:Fstar-starF}
$F(f^*) = F(f)^*$.
\end{lemma}

\pf
The definition of the structure morphisms of the Frobenius
algebra $F(A)$ is given in \eqref{eq:FA-alg-st} and 
\eqref{eq:FA-coalg-st}. Substituting these definitions gives
\be\begin{array}{l}\displaystyle
F(f)^*
=
\big[
\big( \psi_0^F \circ F(\eps_B) \circ F(m_B) \circ \phi_2^F \big)
\otimes \id_{F(A)} \big] 
\circ
\big[ \id_{F(B)} \otimes F(f) \otimes \id_{F(A)} \big]
\\[.2em] \displaystyle
\hspace*{8em}
\circ
\big[ \id_{F(B)} \otimes \big( \psi_2^F \circ F(\Delta_A)
\circ F(\eta_A) \circ \phi_0^F \big)\big]
\enl
=
(\psi_0^F \otimes \id_{F(A)}) \circ
\big[
F\big(\eps_B \circ m_B \circ (\id_B \otimes f)\big)
\otimes \id_{F(A)} \big] 
\\[.2em] \displaystyle
\hspace*{8em}
\circ
( \phi_2^F \otimes \id_{F(A)} )
\circ
( \id_{F(B)} \otimes \psi_2^F )
\\[.2em] \displaystyle
\hspace*{8em}
\circ
\big[ \id_{F(B)} \otimes F(\Delta_A \circ \eta_A) \big]
\circ (\id_{F(B)} \otimes \phi_0^F) ~.
\eear
\ee
In the middle line of the last expression we can use the
defining property (\ref{Frob-ten-cat-1}) of $F$, namely we substitute
$( \phi_2^F \otimes \id_{F(A)} )
\circ
( \id_{F(B)} \otimes \psi_2^F )
= \psi_2^F \circ \phi_2^F$.
Then $\psi_2^F$ can be moved to the left, 
and $\phi_2^F$ to the right, until they can be omitted 
against $\psi_0^F$ and $\phi_0^F$, respectively, using 
\eqref{R-ten-fun-2} and \eqref{L-ten-fun-2}. This results in 
\be
F(f)^* = 
F\big( (\eps_B \circ m_B \circ (\id_B \otimes f)) \otimes
\id_A \big) \circ
F\big(\id_B \otimes (\Delta_A \circ \eta_A)\big)~,
\ee
which is nothing but $F(f^*)$.
\epf

\subsection{Modular tensor categories}\label{sec:MTC-def}

Let $\Cc$ be a modular tensor category \cite{turaev-bk,baki},
i.e.\ an abelian semi-simple finite $\Cb$-linear ribbon category 
with simple tensor unit $\one$ and a non-degeneracy condition on 
the braiding (to be stated in a moment).
We denote the set of equivalence classes of simple objects in $\Cc$
by $\Ic$, elements in $\Ic$ 
by $i,j,k\in \Ic$ and their representatives by $U_i, U_j, U_k$. 
We also set $U_0=\one$ and for an index $k\in \Ic$ we define $\bar{k}$ by
$U_{\bar{k}} \cong U_k^{\vee}$. 

   Since the tensor unit is simple, we shall for modular tensor categories
   identify $\Hom(\one,\one) \cong \Cb$ (cf.\ footnote \ref{fn:HomUU=C}).
Define numbers $s_{i,j} \in \Cb$ 
by\footnote{
  In the graphical notation used below, we have given an orientation to
  the ribbons indicated by the arrows. For example, it is understood 
  that this orientation determines which of the duality morphisms in
  \eqref{eq:duality-mor} to use.
}
\be  
 s_{i,j}   ~=~  \quad
\raisebox{-30pt}{
  \begin{picture}(110,65)
   \put(0,8){\scalebox{.75}{\includegraphics{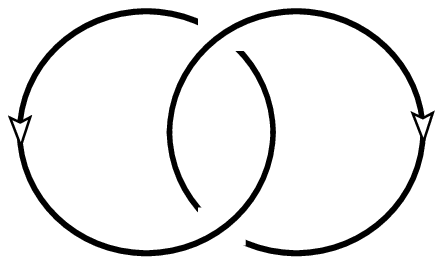}}}
   \put(0,8){
     \setlength{\unitlength}{.75pt}\put(-18,-19){
     \put( 98, 48)       {\scriptsize $ U_i $}
     \put( 50, 48)      {\scriptsize $ U_j $}
     }\setlength{\unitlength}{1pt}}
  \end{picture}}
~~.
\ee
They obey $s_{i,j} = s_{j,i}$ and $s_{0,i} = \dim U_i$, see 
e.g.\ \cite[sect.\,3.1]{baki}. 
(In a ribbon category the left and right dimension 
\eqref{eq:dim-def} of $U_i$ coincide and are denoted by $\dim U_i$.) 
The non-degeneracy condition on the braiding of a modular tensor category is that the $|\Ic|{\times}|\Ic|$-matrix $s$ should be invertible. In fact \cite[thm.\,3.1.7]{baki},
\be
  \sum_{k \in \Ic} s_{ik} \, s_{kj} = 
  \Dim\,\Cc \, \delta_{i,\bar\jmath} ~,
\ee 
where $\Dim\,\Cc = \sum_{i\in\Ic} (\dim U_i)^2$. One can show (even in the weaker context of fusion categories 
over $\Cb$) that $\Dim\,\Cc \ge 1$ \cite[thm.\,2.3]{eno}. 
In particular, $\Dim\,\Cc \neq 0$.
We fix once and for all a square root
$\sqrt{\Dim\,\Cc}$ of $\Dim\,\Cc$.

Let us fix a basis 
$\{ \lambda_{(i,j)k}^{\alpha} \}_{\alpha=1}^{N_{ij}^k}$ 
in $\Hom_{\Cc} (U_i\otimes U_j, U_k)$ and the dual basis
$\{ \Upsilon^{(i,j)k}_{\alpha} \}_{\alpha=1}^{N_{ij}^k}$ 
in $\Hom_{\Cc} (U_k, U_i\otimes U_j)$. 
The duality of the bases means that 
$\lambda_{(i,j)k}^{\alpha} \cir \Upsilon^{(i,j)k}_{\beta}
 = \delta_{\alpha,\beta}\, \id_{U_k}$. We also fix
$\lambda_{(0,i)i} = \lambda_{(i,0)i} = \id_{U_i}$.
We denote the basis vectors graphically as follows: 
\be
\lambda_{(i,j)k}^{\alpha} = 
  \raisebox{-23pt}{
  \begin{picture}(30,52)
   \put(0,8){\scalebox{.75}{\includegraphics{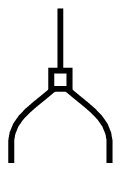}}}
   \put(0,8){
     \setlength{\unitlength}{.75pt}\put(-18,-11){
     \put(39,36)  {\scriptsize $ \alpha $}
     \put(30,61)  {\scriptsize $ U_k $}
     \put(15, 2)  {\scriptsize $ U_i $}
     \put(43, 2)  {\scriptsize $ U_j $}
     }\setlength{\unitlength}{1pt}}
  \end{picture}}
\quad , \qquad
\Upsilon_{\alpha}^{(i,j)k} =
  \raisebox{-23pt}{
  \begin{picture}(30,52)
   \put(0,8){\scalebox{.75}{\includegraphics{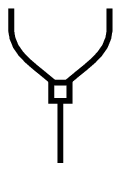}}}
   \put(0,8){
     \setlength{\unitlength}{.75pt}\put(-18,-11){
     \put(39,28)  {\scriptsize $ \alpha $}
     \put(30, 2)  {\scriptsize $ U_k $}
     \put(15,61)  {\scriptsize $ U_i $}
     \put(43,61)  {\scriptsize $ U_j $}
     }\setlength{\unitlength}{1pt}}
  \end{picture}} 
  ~~.
\ee

For $V \in \Cc$ 
we also choose a basis $\{ b_{V}^{(i;\alpha)} \}$ 
of $\Hom_{\Cc}(V, U_i)$
and the dual basis $\{ b_{(i;\beta)}^{V} \}$ of $\Hom_{\Cc}(U_i, V)$ 
for $i\in \Ic$ such that
$b^{(i;\alpha)}_{V} \circ b_{(i;\beta)}^{V} = 
\delta_{\alpha\beta}\, \id_{U_i}$. We use the graphical notation
\be   \label{eq:b-dual-b}
b_V^{(i; \alpha)} = 
  \raisebox{-23pt}{
  \begin{picture}(30,52)
   \put(0,8){\scalebox{.75}{\includegraphics{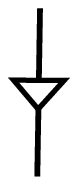}}}
   \put(0,8){
     \setlength{\unitlength}{.75pt}\put(-18,-11){
     \put(39,36)  {\scriptsize $ \alpha $}
     \put(23,65)  {\scriptsize $ U_i $}
     \put(23, 2)  {\scriptsize $ V $}
     }\setlength{\unitlength}{1pt}}
  \end{picture}}
\quad  , \qquad
b_{(i;\alpha)}^V =
  \raisebox{-23pt}{
  \begin{picture}(30,52)
   \put(0,8){\scalebox{.75}{\includegraphics{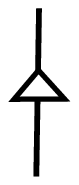}}}
   \put(0,8){
     \setlength{\unitlength}{.75pt}\put(-18,-11){
     \put(37,34)  {\scriptsize $ \alpha $}
     \put(20, 2)  {\scriptsize $ U_i $}
     \put(20,65)  {\scriptsize $ V $}
     }\setlength{\unitlength}{1pt}}
  \end{picture}} 
  ~~.
\ee

Given two modular tensor categories $\Cc$ and $\mathcal{D}$, by $\Cc \boxtimes \mathcal{D}$ we mean the tensor product of additive categories over $\Cb$ \cite[def.\,1.1.15]{baki}, i.e.\ the category whose objects are direct sums of pairs $V \ti W$ of objects $V \in \Cc$ and $W \in \mathcal{D}$ and whose morphism spaces are
\be
 \Hom_{\Cc \boxtimes \mathcal{D}}(V \ti W,V' \ti W') = 
 \Hom_{\Cc}(V,V')\otimes_\Cb \Hom_{\Cc}(W,W')
\ee
for pairs, and direct sums of these if the objects are direct sums of pairs.

If we replace the braiding and the twist in $\Cc$ by the
antibraiding $c^{-1}$ and the antitwist $\theta^{-1}$
respectively, we obtain another ribbon category 
structure on $\Cc$. 
In order to distinguish these two distinct structures, 
we denote $(\Cc, c, \theta)$
and $(\Cc, c^{-1}, \theta^{-1})$ by
$\Cc_+$ and $\Cc_-$ respectively. As in the introduction, we will abbreviate
\be
  \CcC = \Cc_+ \boxtimes \Cc_- ~.
\ee
Note that a set of representatives of the simple objects in $\CcC$ is given by $U_i\ti U_j$ for $i,j\in \Ic$.

\medskip

For the remainder of section \ref{sec:prel} we fix a modular tensor category $\Cc$.

\subsection{The functors $T$ and $R$}\label{sec:TR-def}

The tensor product bifunctor $\otimes$ can be naturally extended to a functor 
$T: \CcC \rightarrow \Cc$. 
Namely, 
$T( \oplus_{i=1}^{N} V_i \ti W_i) = \oplus_{i=1}^N V_i \otimes W_i$
for all $V_i, W_i \in \Cc$ and $N\in \Nb$. The functor $T$ becomes
a tensor functor as follows. For 
$\phi_0^T: \one \rightarrow T(\one \ti \one)$ 
take $\phi_0^T = \id_\one$ (or $l_{\one}^{-1}$ in the non-strict case).
Next notice that, for $U, V, W, X\in \Cc$,
\be\bearll
T(U\ti V) \otimes T(W\ti X) 
\etb= (U\oti V) \otimes (W\oti X), \enl 
T\big( (U\ti V) \oti (W\ti X) \big) \etb= 
(U\oti W) \otimes (V\oti X).
\eear\ee
We define $\phi_2^T: T(U\ti V) \otimes T(W\ti X)
\rightarrow T\big( (U\ti V) \otimes (W\ti X)\big)$ by
\be
\phi_2^T = 
\id_U \otimes c_{WV}^{-1} \otimes \id_X ~.
\ee
(In the non-strict case the appropriate associators have to be added.)
The above definition of $\phi_2^T$ can be naturally extended to a morphism 
$\phi_2^T: T(M_1 \otimes M_2) \rightarrow T(M_1)\otimes T(M_2)$ for 
any pair of objects $M_1, M_2$ in $\CcC$. 
The following result can be checked by direct calculation
\cite[prop.\,5.2]{JS}.

\begin{lemma} \label{lem:T-is-tensor}
The triple $(T, \phi_2^T, \phi_0^T)$ gives a tensor functor. 
\end{lemma}

In particular, $(T, \phi_2^T, \phi_0^T, \psi_2^T, \psi_0^T)$, where 
$\psi_2^T = (\phi_2^T)^{-1}$ and $\psi_0^T = (\phi_0^T)^{-1}$, 
gives a Frobenius functor. 

\medskip

Define the functor $R:\Cc \rightarrow \CcC$ as follows: 
for $A \in \Cc$ and $f\in \Hom_{\Cc}(A,B)$,
\be  
R(A) = \bigoplus_{i\in \Ic}\, (A\otimes U_i^{\vee}) \times U_i
\quad , \quad
R(f) = \bigoplus_{i\in \Ic}\, (f\otimes \id_{U_i^{\vee}}) \times \id_{U_i}
~.
\labl{eq:def-R}
    This functor was also considered in a slightly different context in
    \cite[prop.\,2.3]{Etingof:2004}.
The family of isomorphisms 
$\gamma_A^R= \oplus_{i\in \Ic} \, \tfrac{\Dim \Cc}{\dim U_i} 
\, \id_{(A\otimes U_i^{\vee}) \times U_i} \in \text{Aut}(R(A))$ 
defines a natural isomorphism $\gamma^R: R \rightarrow R$.  

Our next aim is to show that $R$ is left and right adjoint to $T$,
in other words $R$ and $T$ form an 
ambidextrous adjunction 
(see e.g.\ \cite{amadj} for a discussion of ambidextrous adjunctions).
To this end we introduce two linear isomorphisms, for $A \in \Cc$ and $M \in \CcC$,
\be\begin{array}{l}\displaystyle
  \hat\chi : 
  \Hom_{\Cc}(T(M), A) \longrightarrow \Hom_{\CcC}(M, R(A)) ~,
\enl
  \check\chi : 
  \Hom_{\Cc}(A, T(M)) \longrightarrow \Hom_{\CcC}(R(A), M) ~.
\eear
\labl{eq:chi-2-def}
If we decompose $M$ as $M=\oplus_{n=1}^N M_n^l \ti M_n^r$, 
then $\hat{\chi}$ and $\check{\chi}$ are given by
\be   \label{eq:hat-chi}
\hspace{-1.1cm}\hat{\chi}:  \bigoplus_{n=1}^N 
  \raisebox{-45pt}{
  \begin{picture}(50,90)
   \put(0,8){\scalebox{.75}{\includegraphics{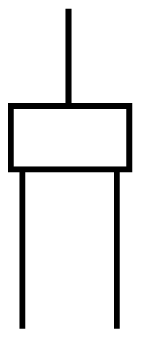}}}
   \put(0,8){
     \setlength{\unitlength}{.75pt}\put(-18,-11){
     \put(37,65)  {\scriptsize $ f_n $}
     \put(23, 2)  {\scriptsize $ M_n^l $}
     \put(50, 2)  {\scriptsize $ M_n^r $}
     \put(37, 110){\scriptsize $A$}
     }\setlength{\unitlength}{1pt}}
  \end{picture}}
\mapsto  \,\, \, \bigoplus_{n=1}^N \bigoplus_{i\in \Ic} \sum_{\alpha}\,\, 
\raisebox{-45pt}{
  \begin{picture}(60,90)
   \put(0,8){\scalebox{.75}{\includegraphics{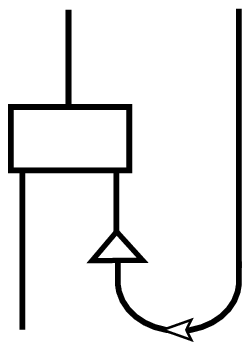}}}
   \put(0,8){
     \setlength{\unitlength}{.75pt}\put(-18,-11){
     \put(38, 69)  {\scriptsize $ f_n $}
     \put(23, 2)    {\scriptsize $ M_n^l $}
     \put(62, 49)  {\scriptsize $ M_n^r $}
     \put(37, 110){\scriptsize $A$}
     \put(38, 35)  {\scriptsize $\alpha$} 
     \put(87,110) {\scriptsize $U_i^{\vee}$}
     }\setlength{\unitlength}{1pt}}
  \end{picture}}
\,\, \times \,\, 
\raisebox{-45pt}{
  \begin{picture}(35,90)
  \put(0,8){\scalebox{.75}{\includegraphics{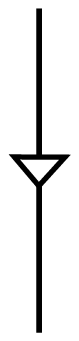}}}
   \put(0,8){
     \setlength{\unitlength}{.75pt}\put(-18,-11){
     \put(37,57)     {\scriptsize $ \alpha $}
     \put(23,2)       {\scriptsize $ M_n^r $}
     \put(20,110)   {\scriptsize $ U_i $}
  }\setlength{\unitlength}{1pt}}
  \end{picture}} 
\hspace{0.2cm}
\ee
and
\be   \label{eq:check-chi}
\check{\chi}:  \bigoplus_{n=1}^N 
  \raisebox{-45pt}{
  \begin{picture}(50,90)
   \put(0,8){\scalebox{.75}{\includegraphics{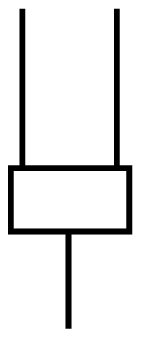}}}
   \put(0,8){
     \setlength{\unitlength}{.75pt}\put(-18,-11){
     \put(38,48)  {\scriptsize $ g_n $}
     \put(23,110)  {\scriptsize $ M_n^l $}
     \put(50,110)  {\scriptsize $ M_n^r $}
     \put(37,2){\scriptsize $A$}
     }\setlength{\unitlength}{1pt}}
  \end{picture}}
\mapsto  \,\, \, \bigoplus_{n=1}^N \bigoplus_{i\in \Ic} \sum_{\alpha}\,\, 
\raisebox{-45pt}{
  \begin{picture}(60,90)
   \put(0,8){\scalebox{.75}{\includegraphics{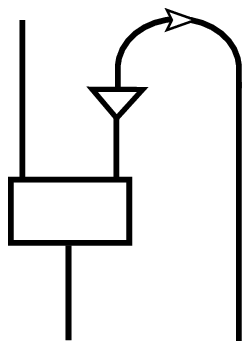}}}
   \put(0,8){
     \setlength{\unitlength}{.75pt}\put(-18,-11){
     \put(38, 49)  {\scriptsize $ g_n $}
     \put(23, 110)    {\scriptsize $ M_n^l $}
     \put(62, 67)  {\scriptsize $ M_n^r $}
     \put(37, 2){\scriptsize $A$}
     \put(39, 81)  {\scriptsize $\alpha$} 
     \put(87, 2) {\scriptsize $U_i^{\vee}$}
     }\setlength{\unitlength}{1pt}}
  \end{picture}}
\,\, \times \,\, 
\raisebox{-45pt}{
  \begin{picture}(55,90)
  \put(0,8){\scalebox{.75}{\includegraphics{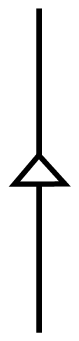}}}
   \put(0,8){
     \setlength{\unitlength}{.75pt}\put(-18,-11){
     \put(37,57)     {\scriptsize $ \alpha $}
     \put(23,2)       {\scriptsize $ U_i $}
     \put(20,110)   {\scriptsize $ M_n^r $}
     \put(55, 55){$\displaystyle\frac{\Dim \Cc}{\dim U_i}$}
  }\setlength{\unitlength}{1pt}}
\end{picture}} 
  \hspace{0.5cm}.
\ee
Notice that $\hat\chi$ and $\check\chi$ are independent 
of the choice of basis.

\begin{thm}\label{thm:RT-adjoint}
$\langle T,R,\hat\chi \rangle$ and 
$\langle R,T,\check\chi{}^{-1} \rangle$ are adjunctions, 
i.e.\ $R$ is both left and right adjoint of $T$. 
\end{thm}
\pf
Write $M$ as $M= \oplus_{n=1}^N M_{n}^l \ti M_{n}^r$. 
The isomorphism $\hat\chi$ amounts to the following composition
of natural isomorphisms,
\be\begin{array}{l}\displaystyle
\Hom_{\Cc}(T(M), A) ~=~ 
  \oplus_{n} \Hom_{\Cc}(M_{n}^l\otimes M_{n}^r, A)  
\enl \qquad \cong
  \oplus_{n, i} 
  \Hom_{\Cc}(M_{n}^l\otimes U_i, A) \otimes \Hom_{\Cc}(M_{n}^r, U_i)  
\enl \qquad \cong
  \oplus_{n, i} \Hom_{\Cc}(M_{n}^l, A\otimes U_i^{\vee}) \otimes 
  \Hom_{\Cc}(M_{n}^r, U_i)  
  ~=~ \Hom_{\CcC}(M, R(A))  ~.
\eear
\labl{chi-1}
Thus $\hat\chi$ is natural. 
Let $(\gamma_A^R)^*: \Hom_{\CcC}
(R(A), M)\rightarrow \Hom_{\CcC}
(R(A), M)$ denote the pull-back of $\gamma_A^R$. 
The isomorphism $\check{\chi}$ is equal to 
the composition of $(\gamma_A^R)^*$
and the following sequence of natural isomorphisms,
\be\begin{array}{l}\displaystyle
\Hom_{\Cc}(A, T(M)) ~=~
  \oplus_{n} \Hom_{\Cc}(A, M_{n}^l\otimes M_{n}^r)
\enl \qquad \cong
  \oplus_{n, i} \Hom_{\Cc}(A, M_{n}^l\otimes U_i) \otimes
  \Hom_{\Cc}(U_i, M_{n}^r) 
\enl \qquad \cong
  \oplus_{n, i} \Hom_{\Cc}(A\otimes U_i^{\vee}, M_{n}^l) \otimes 
  \Hom_{\Cc}(U_i, M_{n}^r)   ~=~\Hom_{\CcC}(R(A), M) ~.
\eear
\labl{chi-2}
We have proved that both $\hat\chi$ and $\check\chi$ are natural isomorphisms.
\epf

\medskip

There are four natural transformations associated to $\hat\chi$ and $\check\chi$, namely
\be
\id_{\CcC} 
\xrightarrow{\hat{\delta}} RT \xrightarrow{\check{\rho}} 
\id_{\CcC} \hspace{1cm} \text{and} \hspace{1cm} 
\id_{\Cc} 
\xrightarrow{\check{\delta}} TR \xrightarrow{\hat{\rho}} 
\id_{\Cc},
\ee
defined by, for $A\in \Cc$, $M \in \CcC$,
\be\begin{array}{ll}\displaystyle
\hat{\delta}_M = \hat{\chi}(\id_{T(M)})~, 
\hspace{2cm} \etb \hat{\rho}_A = \hat{\chi}^{-1}(\id_{R(A)})~, \hspace{2cm}   
\enl
\check{\rho}_M = \check{\chi}(\id_{T(M)})~,
\etb
\check{\delta}_A = \check{\chi}^{-1}(\id_{R(A)})~.
\eear
\ee
They can be expressed graphically as follows, 
with $M=\oplus_{n=1}^N M_n^l \ti M_n^r$,
\be \begin{array}{l}\displaystyle
\hat{\delta}_M = \bigoplus_{n,i}\sum_{\alpha}  
  \raisebox{-25pt}{
  \begin{picture}(60,65)
   \put(0,8){\scalebox{.75}{\includegraphics{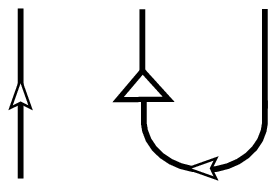}}}
   \put(0,8){
     \setlength{\unitlength}{.75pt}\put(-18,-11){
     \put(15, 2)  {\scriptsize $ M_n^l $}
     \put(15, 68)  {\scriptsize $ M_n^l $}
     \put(50, 68)  {\scriptsize $ M_n^r $}
     \put(84, 68){\scriptsize $U_i^{\vee}$}
     \put(38, 37){\scriptsize $\alpha$}
     }\setlength{\unitlength}{1pt}}
  \end{picture}}
 ~ \times ~
\raisebox{-25pt}{
  \begin{picture}(20,65)
   \put(0,8){\scalebox{.75}{\includegraphics{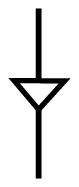}}}
   \put(0,8){
     \setlength{\unitlength}{.75pt}\put(-18,-11){
     \put(21, 68)  {\scriptsize $ U_i $}
     \put(21, 2)    {\scriptsize $ M_n^r $}
     \put(37, 36)  {\scriptsize $\alpha$} 
     }\setlength{\unitlength}{1pt}}
  \end{picture}}
\quad,\hspace{1.5cm}
\hat{\rho}_A = \bigoplus_{i\in \Ic} \,\,\,\,
  \raisebox{-25pt}{
  \begin{picture}(60,65)
   \put(0,8){\scalebox{.75}{\includegraphics{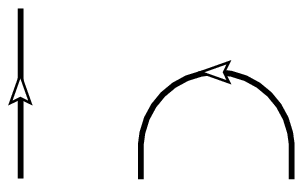}}}
   \put(0,8){
     \setlength{\unitlength}{.75pt}\put(-18,-11){
     \put(14, 2)  {\scriptsize $ A $}
     \put(14, 68)  {\scriptsize $ A $}
     \put(48, 2)  {\scriptsize $ U_i^{\vee} $}
     \put(93, 2){\scriptsize $ U_i $}
     }\setlength{\unitlength}{1pt}}
  \end{picture}}
  ~~,
\enl 
\check{\rho}_M = \bigoplus_{n,i} \sum_{\alpha}
  \raisebox{-25pt}{
  \begin{picture}(60,70)
   \put(0,8){\scalebox{.75}{\includegraphics{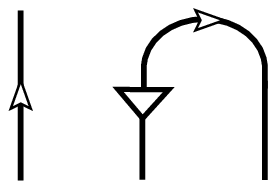}}}
   \put(0,8){
     \setlength{\unitlength}{.75pt}\put(-18,-11){
     \put(15, 2)  {\scriptsize $ M_n^l $}
     \put(15, 68)  {\scriptsize $ M_n^l $}
     \put(50, 2)  {\scriptsize $ M_n^r $}
     \put(84, 2){\scriptsize $U_i^{\vee}$}
     \put(38, 37){\scriptsize $\alpha$}
     }\setlength{\unitlength}{1pt}}
  \end{picture}}
~ \times ~
\raisebox{-25pt}{
  \begin{picture}(40,70)
   \put(0,8){\scalebox{.75}{\includegraphics{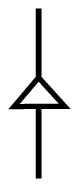}}}
   \put(0,8){
     \setlength{\unitlength}{.75pt}\put(-18,-11){
     \put(21, 2)  {\scriptsize $ U_i $}
     \put(21, 68)    {\scriptsize $ M_n^r $}
     \put(37, 36)  {\scriptsize $\alpha$}      
     }\setlength{\unitlength}{1pt}}
  \end{picture}}
\hspace{-0.5cm}\mbox{\small $\frac{\text{Dim} \Cc}{\dim U_i}$}
~~ , ~~~ 
\check{\delta}_A =  \bigoplus_{i\in \Ic} 
   ~~
  \raisebox{-25pt}{
  \begin{picture}(60,70)
   \put(0,8){\scalebox{.75}{\includegraphics{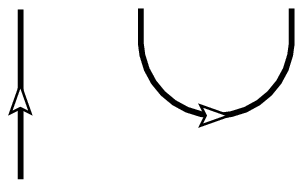}}}
   \put(0,8){
     \setlength{\unitlength}{.75pt}\put(-18,-11){
     \put(15, 2)  {\scriptsize $ A $}
     \put(15, 68)  {\scriptsize $ A $}
     \put(48, 68)  {\scriptsize $ U_i^{\vee} $}
     \put(93, 68){\scriptsize $ U_i $}
   }\setlength{\unitlength}{1pt}}
  \end{picture}}
~~~\mbox{\small $\frac{\dim U_i}{\text{Dim} \Cc}$ .}
\eear
\labl{eq:delta-rho-hat-check}
Note that
\be
  \check\rho_M \circ \hat\delta_M = \Dim\,\Cc \cdot \id_M 
  \qquad \text{and} \qquad
  \hat\rho_A \circ \check\delta_A = \id_A ~.
\labl{eq:delta-rho-id}

\begin{lemma} \label{lem:RT-injective}
The functors $T$ and $R$ as maps on the sets of morphisms have left inverses, and thus are injective.
\end{lemma}
\pf
Let $f : A\rightarrow B$ be a morphism in $\Cc$.
We define a map $Q_R:\Hom_{\CcC}(R(A), R(B)) \rightarrow \Hom_{\Cc}(A, B)$ by $f' \mapsto \hat\rho_B \circ T(f') \circ \check\delta_A$. Then we have  
\be
Q_R \circ R(f) = \hat\rho_B \circ TR(f) \circ \check\delta_A 
=
\hat\rho_B \circ \check\delta_B  \circ f
= f ~,
\ee
where we used naturality of $\check\delta$ and \eqref{eq:delta-rho-id}
in the second and third equalities, respectively. So $Q_R$ is a left inverse of $R$ on morphisms. 
Thus $R$ is injective on morphisms.
Similarly, let $g : M \rightarrow N$ be a morphism in $\CcC$.
We define a map $Q_T: \Hom_{\Cc}(T(M), T(N)) \rightarrow \Hom_{\CcC}(M, N)$ by 
$g'\mapsto (\Dim\,\Cc)^{-1} \cdot \check\rho_N \circ R(g') \circ \hat\delta_M$. Then we have 
\be
Q_T \circ T(g) 
=
(\Dim\,\Cc)^{-1} \cdot \check\rho_N \circ RT(g) \circ \hat\delta_M 
=
(\Dim\,\Cc)^{-1} \cdot \check\rho_N \circ \hat\delta_N  \circ g
= g ~.
\ee
So $Q_T$ is a left inverse of $T$ on morphisms. Thus $T$ is injective on morphisms. 
\epf

\medskip

Using \eqref{eq:chi-eta-epsilon} and
\eqref{eq:delta-rho-hat-check},
one can express the two inverse maps $\hat{\chi}^{-1}, \check{\chi}^{-1}$ as follows, for $f\in \Hom_{\CcC}(M, R(A))$ and $g\in \Hom_{\CcC}(R(A), M)$,
\be
\hat{\chi}^{-1}(f) = \hat{\rho} \circ T(f) 
\qquad , \qquad
\check{\chi}^{-1}(g) = T(g) \circ \check{\delta}~.
\labl{eq:chi-inv}
By proposition \ref{prop:tensor-implies-frobenius}
and lemma \ref{lem:lax-colax-adj},
$R$ is both a lax and colax tensor functor.
In particular, 
$\phi_0^R: \one \ti \one \rightarrow R(\one)$ is given by
\be  \label{eq:phi-0-R}
\phi_0^R =\hat{\chi}(\psi_0^T)= R(\psi_0^T)\circ \hat{\delta}_{\one\ti \one} = 
\id_{\one \times \one}
\ee
and $\phi_2^R : R(A)\otimes R(B)\rightarrow R(A\otimes B)$ by 
$\phi_2^R = R(\hat{\rho}_A\otimes \hat{\rho}_B) \circ R(\psi_2^T) 
\circ \hat{\delta}$,
which can be expressed graphically as 
\be   \label{eq:phi-2-R}
\phi_2^R ~=
\bigoplus_{i,j,k\in \Ic} \sum_{\alpha} ~~
  \raisebox{-45pt}{
  \begin{picture}(170,90)
   \put(0,8){\scalebox{.75}{\includegraphics{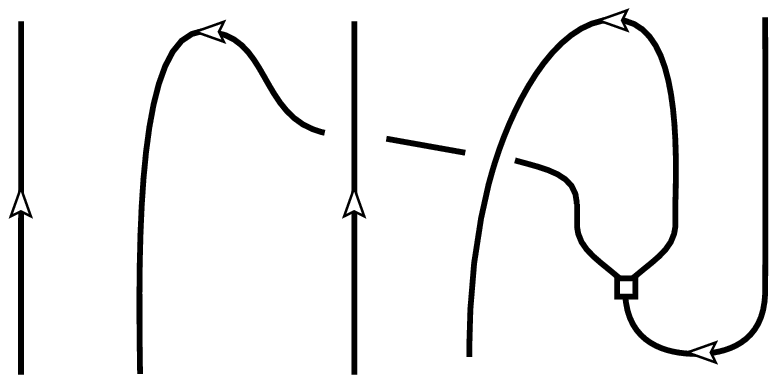}}}
   \put(0,8){
     \setlength{\unitlength}{.75pt}\put(-18,-11){
     \put(15, 2)        {\scriptsize $ A $}
     \put(15, 120)    {\scriptsize $ A $}
     \put(52, 2)        {\scriptsize $ U_i^{\vee} $}
     \put(113,2)       {\scriptsize $ B $}
     \put(113, 120)  {\scriptsize $B$  }
     \put(147, 2)      {\scriptsize $U_j^{\vee}$}
     \put(234, 120)  {\scriptsize $U_k^{\vee}$}
     \put(193, 46)      {\scriptsize $\alpha$}
     }\setlength{\unitlength}{1pt}}
  \end{picture}}
~~  \times ~~
\raisebox{-45pt}{
  \begin{picture}(20,90)
   \put(0,8){\scalebox{.75}{\includegraphics{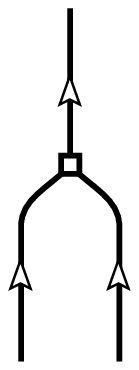}}}
   \put(0,8){
     \setlength{\unitlength}{.75pt}\put(-18,-11){
     \put(18, 2)     {\scriptsize $ U_i $}
     \put(46, 2)     {\scriptsize $ U_j $}
     \put(32, 56)   {\scriptsize $\alpha$} 
     \put(32, 120) {\scriptsize $ U_k $ }
     }\setlength{\unitlength}{1pt}}
  \end{picture}}
  \quad ~~.
  \phantom{\frac{\dim U_i \,\dim U_j}{\dim U_k \,\Dim\,\Cc}} 
\ee
Similarly, $\psi_0^R: R(\one) \rightarrow \one\ti \one$ is given by 
\be  \label{eq:psi-0-R}
\psi_0^R = \check{\rho}_{\one} \circ R(\phi_0^T) = 
\Dim\,\Cc ~\id_{\one\times\one} 
\ee
and $\psi_2^R : R(A\oti B) \rightarrow R(A) \otimes R(B)$ by 
$\psi_2^R = \check{\rho} \circ R(\phi_2^T) \circ R(\check{\delta}_A \otimes
\check{\delta}_B)$,
which in graphical notation reads
\be   \label{eq:psi-2-R}
\psi_2^R ~= \bigoplus_{i,j,k\in \Ic} \sum_{\alpha} ~~
  \raisebox{-45pt}{
  \begin{picture}(170,90)
   \put(0,8){\scalebox{.75}{\includegraphics{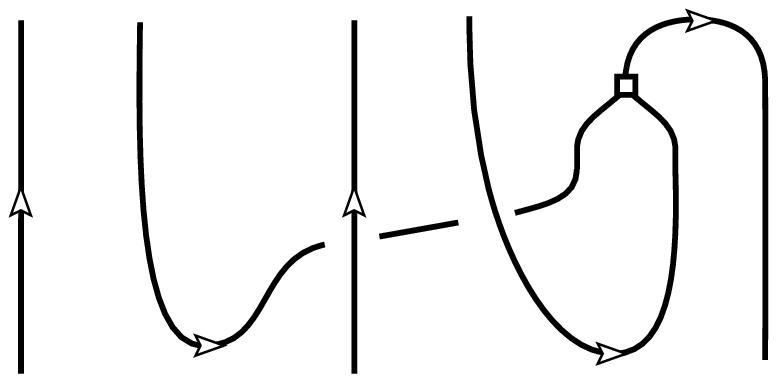}}}
   \put(0,8){
     \setlength{\unitlength}{.75pt}\put(-18,-11){
     \put(15, 2)        {\scriptsize $ A $}
     \put(15, 120)    {\scriptsize $ A $}
     \put(52, 120)        {\scriptsize $ U_i^{\vee} $}
     \put(113,2)       {\scriptsize $ B $}
     \put(113, 120)  {\scriptsize $B$  }
     \put(147, 120)      {\scriptsize $U_j^{\vee}$}
     \put(234, 2)  {\scriptsize $U_k^{\vee}$}
     \put(193, 82)      {\scriptsize $\alpha$}
     }\setlength{\unitlength}{1pt}}
  \end{picture}}
~~  \times ~~
\raisebox{-45pt}{
  \begin{picture}(20,90)
   \put(0,8){\scalebox{.75}{\includegraphics{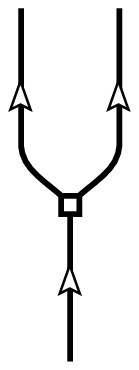}}}
   \put(0,8){
     \setlength{\unitlength}{.75pt}\put(-18,-11){
     \put(18, 120)     {\scriptsize $ U_i $}
     \put(46, 120)     {\scriptsize $ U_j $}
     \put(32, 68)       {\scriptsize $\alpha$} 
     \put(32, 2)         {\scriptsize $ U_k $ }
     }\setlength{\unitlength}{1pt}}
  \end{picture}}
~\quad \frac{\dim U_i \,\dim U_j}{\dim U_k \,\Dim\,\Cc}~.
\ee

If $\Cc$ has more than one simple object, then $R$ does not take the tensor unit of $\Cc$ to the tensor unit of $\CcC$ and so is clearly not a tensor functor. However, we will show that $R$ is still a Frobenius functor.
This will imply that if $A$ is a Frobenius algebra in $\Cc$, then
\be
   R(A) = ( R(A), m_{R(A)}, \eta_{R(A)},\Delta_{R(A)}, \eps_{R(A)} )
\labl{eq:RA-Frobalg}
is a Frobenius algebra in $\CcC$, where the structure morphisms were given
in \eqref{eq:FA-alg-st} and \eqref{eq:FA-coalg-st}. 
In the case $A=\one$ 
it was proved in  \cite[prop.\,4.1]{mug2} (see also 
\cite[lem.\,6.19]{corr} and \cite[thm.\,5.2]{ko-ffa}) that 
\eqref{eq:RA-Frobalg} 
is a commutative simple symmetric normalised-special Frobenius algebra 
in $\CcC$. 
In fact, given a Frobenius algebra $A$ in $\Cc$, 
it is straightforward to verify 
that the structure morphisms in \eqref{eq:RA-Frobalg} are precisely those of
$(A\times \one)\otimes R(\one)$, cf.\ section \ref{sec:algebra}.

\begin{prop}  \label{prop:R-Frob}
$(R, \phi_2^R, \phi_0^R, \psi_2^R, \psi_0^R)$ is a Frobenius functor. 
\end{prop}
\pf
Using the explicit graphical expression of $\phi_2^R, \phi_0^R, \psi_2^R, \psi_0^R$, 
it is easy to see that the commutativity of the diagrams 
(\ref{Frob-ten-cat-1}) and (\ref{Frob-ten-cat-2}) are equivalent to 
the statement that $R(\one)$ with structure morphisms as in \eqref{eq:RA-Frobalg}
is a Frobenius algebra in $\CcC$. The latter statement is true by 
\cite[prop.\,4.1]{mug2}.
\epf

\medskip

From lemma \ref{lem:T-is-tensor} 
and proposition \ref{prop:R-Frob} we see that $T$ and $R$
take Frobenius algebras to Frobenius algebras. The following two
propositions show how the properties of Frobenius algebras are 
transported.

\begin{prop}  \label{prop:symm-T-symm}
Let $A$ be a Frobenius algebra in $\CcC$. Then $T(A)$ is a 
Frobenius algebra in $\Cc$ and
\\
(i)\phantom{i} $A$ is symmetric iff $T(A)$ is symmetric.
\\
(ii) $A$ is (normalised-)special iff $T(A)$ is (normalised-)special.
\end{prop}

\medskip\noindent
\pf
For part (i) write $A$ as a direct sum 
$\oplus_{n=1}^N A_n^{l} \ti A_n^{r}$. Then 
the maps $I_{T(A)}, I_{T(A)}': T(A^{\vee}) \rightarrow T(A)^{\vee}$ defined
in \eqref{equ:II'-def} are given by: 
\be  \label{eq:I-I'-T-c}
I_{T(A)} = I_{T(A)}' = 
\oplus_{n=1}^N c_{(A_n^{l})^{\vee},(A_n^{r})^{\vee}}~.
\ee 
Therefore, by (\ref{eq:Phi-I}), 
$\Phi_{T(A)}=\Phi_{T(A)}'$ is equivalent
to $T(\Phi_A) = T(\Phi_A')$. Since by lemma \ref{lem:RT-injective}, $T$ is
injective on morphisms, this proves part (i). 
Part (ii) can be checked in the same way, for example the condition
$m_{T(A)} \circ \Delta_{T(A)} = \zeta \, \id_{T(A)}$ is easily
checked to be equivalent to 
$T(m_A \circ \Delta_A) = \zeta \, T(\id_A)$.
\epf

\begin{prop}  \label{prop:symm-R-symm}
Let $A$ be a Frobenius algebra $A$ in $\Cc$. Then $R(A)$ is a 
Frobenius algebra in $\CcC$ and
\\
(i)\phantom{i} $A$ is symmetric iff $R(A)$ is symmetric.
\\
(ii) $A$ is (normalised-)special iff $R(A)$ is (normalised-)special.
\end{prop}
\pf
Recall that the structure morphisms of the Frobenius
algebra $R(A)$ are equal to those of 
$(A\times \one)\otimes R(\one)$. 
Using this equality, part (i) and (ii)
follow because $R(\one)$ is symmetric and normalised-special.
For example, 
\be
  m_{R(A)} \circ \Delta_{R(A)} = 
  \big[(m_A \circ \Delta_A) \times \id_\one\big] 
  \otimes \id_{R(\one)}
  = R(m_A \circ \Delta_A) ~,
\ee  
so that $m_{R(A)} \circ \Delta_{R(A)} = \zeta \, \id_{R(A)}$
is equivalent to $R(m_A \circ \Delta_A) = \zeta \, R(\id_A)$,
which by lemma \ref{lem:RT-injective} is equivalent to
$m_A \circ \Delta_A = \zeta \,\id_A$.
\epf

\medskip

The functor $R$ has one additional property not shared by $T$, namely 
$R$ takes absolutely simple algebras to 
absolutely simple algebras.  We will see 
explicitly in section \ref{sec:cardy-reconstr} that this is
not true for $T$.

\begin{lemma}  \label{lem:R-iso-hom}
For a $\Cc$-algebra $A$, the map 
\be   \label{eq:R-AA-RARA}
R: \Hom_{A|A}(A,A) \rightarrow \Hom_{R(A)|R(A)}(R(A), R(A))
\ee
given by $f\mapsto R(f)$ is well-defined and an isomorphism. 
\end{lemma}
\pf
Since $R$ is a lax tensor functor, $R(A)$ is naturally a $R(A)$-bimodule. 
It is easy to see that $R$ in (\ref{eq:R-AA-RARA}) is a well-defined map. 
$R(A)$ is also naturally a $R(\one)$-bimodule, which can be identified with 
the induced $R(\one)$-bimodule structure on 
$(A\times \one)\otimes R(\one)$, where the left $R(\one)$ action on 
$(A\times \one)\otimes R(\one)$
is given by $(\id_{A\times \one} \otimes m_{R(\one)}) \circ (c_{A\times 1, R(\one)}^{-1}\otimes \id_{R(\one)})$. We have the following natural isomorphisms: 
\bea
\Hom_{R(\one)|R(\one)}(R(A), R(A)) \xrightarrow{\cong} 
\Hom_{\Cc_{\pm}^2}(A\times \one, R(A)) 
\xrightarrow{\hat{\chi}^{-1}} \Hom_{\Cc}(A, A).
\eea
which, by (\ref{eq:chi-inv}), 
are given by, for $f \in \Hom_{R(\one)|R(\one)}(R(A), R(A))$, 
\be  \label{eq:f-f'-f''}
f\mapsto f'= f\circ (\id_{A\times \one} \otimes \eta_{R(\one)})
\mapsto f''= \hat{\rho} \circ T(f'),
\ee
and its inverse is given by, for $g\in \Hom_{\Cc}(A,A)$,
\be  \label{eq:g-g'-g''}
g \mapsto g'=R(g) \circ \hat{\delta} \mapsto g'' = 
(\id_{A\times \one}\otimes m_{R(\one)}) \circ (g' \otimes \id_{R(\one)})
\ee
where $g''$ is indeed a $R(\one)$-bimodule map due to the commutativity of $R(\one)$.  It is easy to check that $g''=R(g)$ in (\ref{eq:g-g'-g''}).  
Therefore $R$ gives an isomorphism from $\Hom_{\Cc}(A, A)$ to 
$\Hom_{R(\one)|R(\one)}(R(A), R(A))$. Moreover, 
one verifies that $R(g)$ is an $R(A)$-bimodule map iff $g$ is 
an $A$-bimodule map. In other words, 
$R: g\mapsto R(g)$ gives an isomorphism 
$\Hom_{A|A}(A,A) \xrightarrow{\cong} \Hom_{R(A)|R(A)}(R(A), R(A))$. 
\epf

\begin{cor}\label{cor:A-simp-RA-simp}
Let $A$ be a $\Cc$-algebra. \\
(i)\phantom{i} 
$A$ is absolutely simple iff $R(A)$ is absolutely simple. \\
(ii) Let $A$ be in addition Frobenius. Then $A$ is simple and 
special iff $R(A)$ is simple and special.
\end{cor}

\pf
Part (i) immediately follows from lemma \ref{lem:R-iso-hom}.
The statement of part (ii) without the qualifier `simple' is proved in 
proposition \ref{prop:symm-R-symm}. But, as in
the proof of (iii)$\Rightarrow$(i) in lemma 
\ref{lem:modinv-simp-hapl}, a special
Frobenius algebra in a semi-simple category has a semi-simple
category of bimodules, and
for a semi-simple $\Cb$-linear
category, simple and absolutely
simple are equivalent.
Part (ii) then follows from part (i).  
\epf

\medskip

The following lemma will be needed in section \ref{sec:cardy-def} below to discuss the properties of Cardy algebras.

\begin{lemma}\label{lem:delta-star-rho}
Let $A$ be a Frobenius algebra in $\Cc$. Then
$(\check\delta_A)^* = \hat\rho_A$.
\end{lemma}

\pf
Recall from \eqref{eq:delta-rho-hat-check} that $\check\delta_A$ is a 
morphism $A \rightarrow TR(A)$.
Since $T$ and $R$ are both Frobenius functors, $TR(A)$ is a Frobenius algebra in $\Cc$. Substituting the definitions, after a short calculation one finds
\be
  \eps_{TR(A)} \circ m_{TR(A)} ~=~ \text{Dim}(\Cc)~ \bigoplus_{i \in \Ic}
  ~~
\raisebox{-40pt}{
  \begin{picture}(150,100)
   \put(0,8){\scalebox{.75}{\includegraphics{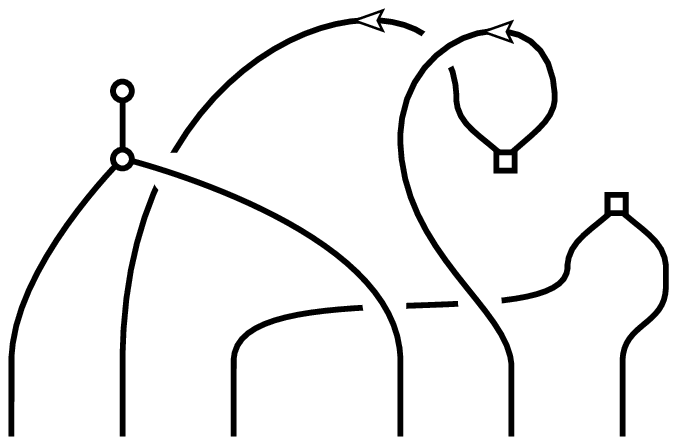}}}
   \put(0,8){
     \setlength{\unitlength}{.75pt}\put(-34,-37){
     \put(31, 28)  {\scriptsize $ A $}
     \put(61, 28)  {\scriptsize $ U_i^\vee $}
     \put(94, 28)  {\scriptsize $ U_i $}
     \put(141, 28)  {\scriptsize $ A $}
     \put(174, 28)  {\scriptsize $ U_{\bar\imath}^\vee $}
     \put(209, 28)  {\scriptsize $ U_{\bar\imath} $}
     }\setlength{\unitlength}{1pt}}
  \end{picture}}
\ee
Substituting this in the definition of $(\check\delta_A)^*$ gives, again after a short calculation, the morphism $\hat\rho_A$. At an intermediate step one uses that the part of the morphism 
$(\check\delta_A)^*$, which is made up of $U_i$ and $U_{\bar\imath}$ ribbons, their duals, and the basis morphisms 
$\lambda_{(i,\bar\imath)0}$ and
$\Upsilon^{(i,\bar\imath)0}$, can be replaced by 
$\frac{1}{\dim U_i} \cdot d_{U_i}$. 
\epf

\newpage

\sect{Cardy algebras}

In this section we start by investigating the properties of Frobenius algebras which satisfy the so-called modular invariance condition. We then give two definitions of a Cardy algebra and prove their equivalence. Finally, in section \ref{sec:cardy-reconstr}, 
we study the properties of these algebras
and state our main results. 

We fix a modular tensor category $\Cc$. Recall that $\CcC$ is an abbreviation 
for $\Cc_+ \boxtimes \Cc_-$.

\subsection{Modular invariance}\label{sec:mod-inv}

In $\CcC$, we define 
the object $K$ and the morphism $\omega: K\rightarrow K$ as 
\be
K = \bigoplus_{i, j\in \Ic} \,\, U_i \times U_j~~,  \qquad
\omega = \sum_{i,j\in \Ic} 
\frac{\dim U_i \, \dim U_j}{\Dim\,\Cc} \, \id_{U_i\times U_j} ~ . 
\ee
They have the property (see e.g.\ \cite[cor.\,3.1.11]{baki})
\be
 \raisebox{-36pt}{\begin{picture}(85,78)
  \put(0,8){\scalebox{0.75}{\includegraphics{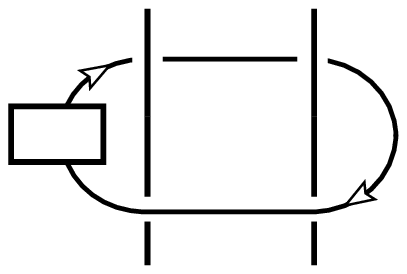}}}
  \put(0,8){
     \setlength{\unitlength}{.75pt}\put(-65,-33){
     \put(74,70)   {\scriptsize $\omega$ }
     \put(78,93)   {\scriptsize $K$ }
     \put(80,25)   {\scriptsize $(U_k\,{\times}\,U_l)^\vee$ }
     \put(80,113)   {\scriptsize $(U_k\,{\times}\,U_l)^\vee$ }
     \put(140,25)   {\scriptsize $U_i\,{\times}\,U_j$ }
     \put(140,113)   {\scriptsize $U_i\,{\times}\,U_j$ }
     }\setlength{\unitlength}{1pt}}
  \end{picture}}
  ~=~ 
 \delta_{i,k} \, \delta_{j,l} \,
 \frac{\Dim\,\Cc}{\dim U_i \,\dim U_j} ~
 \tilde b_{U_i \times U_j} \circ d_{U_i \times U_j}
\labl{eq:Kloop-prop}

\begin{defn}\label{def:mod-inv}
\rule{0pt}{1em}\\[.1em]
(i) Let $A,B$ be objects of $\CcC$. A morphism $f : A \oti B \rightarrow B$ is called 
{\em S-invariant} iff
\be   
  \raisebox{-55pt}{
  \begin{picture}(102,120)
   \put(0,8){\scalebox{.75}{\includegraphics{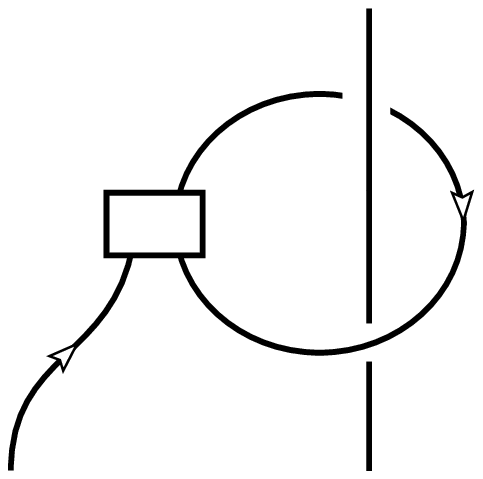}}}
   \put(0,8){
     \setlength{\unitlength}{.75pt}\put(-32,-25){
     \put( 27, 15)  {\scriptsize $ A $}
     \put( 88,133)  {\scriptsize $ B $}
     \put(130, 15)  {\scriptsize $W$ }
     \put(130,165)  {\scriptsize $W$ }
     \put(71,  95)  {\scriptsize $f$} 
     }\setlength{\unitlength}{1pt}}
  \end{picture}}
\quad = \quad 
 \raisebox{-55pt}{
  \begin{picture}(110,120)
   \put(0,8){\scalebox{.75}{\includegraphics{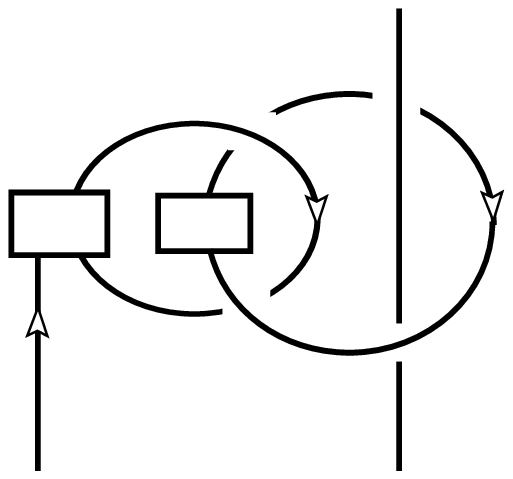}}}
   \put(0,8){
     \setlength{\unitlength}{.75pt}\put(-23,-25){
     \put( 27, 15)  {\scriptsize $ A $}
     \put( 45,124)  {\scriptsize $ B $}
     \put(130, 15)  {\scriptsize $W$ }
     \put(130,165)  {\scriptsize $W$ }
     \put(35,  95)  {\scriptsize $f$} 
     \put(108,139)  {\scriptsize $ K $}
     \put( 75, 96)  {\scriptsize $\omega$}
}\setlength{\unitlength}{1pt}}
  \end{picture}}
\labl{eq:mod-inv}
holds for all for $W\in \CcC$.
\\[.3em]
(ii) A $\CcC$-algebra $(\Acl, m_\text{\rm{cl}},\eta_\text{\rm{cl}})$ 
is called {\em modular invariant} iff 
$\theta_{\Acl} = \id_{\Acl}$ and $m_\text{\rm{cl}}$ is S-invariant.
\end{defn}

\begin{lemma}  \label{lemma:S-inv}
The morphism $f : A \oti B \rightarrow B$ is S-invariant if and only if
\be   
  \raisebox{-55pt}{
  \begin{picture}(102,120)
   \put(0,8){\scalebox{.75}{\includegraphics{pic-modinvf-L.eps}}}
   \put(0,8){
     \setlength{\unitlength}{.75pt}\put(-32,-25){
     \put( 27, 15)  {\scriptsize $ A $}
     \put( 88,133)  {\scriptsize $ B $}
     \put(120, 15)  {\scriptsize $ U_i\ti U_j $}
     \put(120,165)  {\scriptsize $ U_i\ti U_j $}
     \put(71,  95)  {\scriptsize $f$} 
     }\setlength{\unitlength}{1pt}}
  \end{picture}}
\quad = \quad 
\frac{\Dim\,\Cc}{\dim U_i \,\dim U_j} ~ \sum_{\alpha} \hspace{-.3em}
 \raisebox{-55pt}{
  \begin{picture}(60,120)
   \put(0,8){\scalebox{.75}{\includegraphics{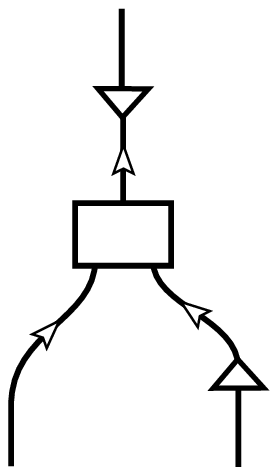}}}
   \put(0,8){
     \setlength{\unitlength}{.75pt}\put(-46,-26){
     \put( 43, 15)  {\scriptsize $ A $}
     \put(102, 76)  {\scriptsize $ B $}
     \put( 85,109)  {\scriptsize $ B $}
     \put( 77, 91)  {\scriptsize $f$} 
     \put( 96, 15)  {\scriptsize $ U_i\ti U_j $}
     \put( 63,165)  {\scriptsize $ U_i\ti U_j $}
     \put( 63,129)  {\scriptsize $\alpha$}
     \put( 95, 51)  {\scriptsize $\alpha$}
}\setlength{\unitlength}{1pt}}
  \end{picture}}
\labl{eq:mod-inv-2}
holds for all $i,j \in \Ic$.
\end{lemma}

\pf
Condition \eqref{eq:mod-inv} holds for
all $W$ iff it holds for all $W = U_i \ti U_j$, $i,j \in \Ic$, so it
is enough to show that the right hand side of
\eqref{eq:mod-inv} with $W = U_i \ti U_j$
is equal to the right hand side of \eqref{eq:mod-inv-2}.
Recall the notation for basis morphisms in
\eqref{eq:b-dual-b}.
Starting from \eqref{eq:mod-inv}, write
\be
  \id_{B^\vee} \otimes \id_{U_i \times U_j}
  = \Big(\sum_{k,l,\alpha} 
\big(b_{B}^{(k\times l; \alpha)}\big)^\vee \circ  
\big(b_{(k\times l;\alpha)}^{B}\big)^\vee
\Big) \otimes \id_{U_i \times U_j} ~,
\ee
and then apply \eqref{eq:Kloop-prop}. The 
graphical representation of the resulting morphism 
can be deformed to give \eqref{eq:mod-inv-2}.
\epf

\begin{rema} {\rm  
As shown in \cite[sect.\,6.1]{cardy}, the modular invariance condition
of a $\CcC$-algebra exactly coincides with the modular invariance 
condition for torus 1-point correlation functions of a genus-0,1 closed
CFT. In particular, the condition $\theta_{\Acl}=\id_{\Acl}$ is equivalent to 
invariance under the modular transformation $T: \tau \mapsto \tau +1$, and 
the condition (\ref{eq:mod-inv-2}) with $f=m_\text{\rm{cl}}$ 
is equivalent to invariance under $S: \tau \mapsto -\frac{1}{\tau}$. 
Combining the modular invariance condition with the genus zero 
properties of a genus-0,1 closed CFT results in a
modular invariant commutative symmetric Frobenius algebra in $\CcC$. 
}
\end{rema}

Let $\Acl$ be a modular invariant $\CcC$-algebra. Evaluating (\ref{eq:mod-inv-2}) for $f=m_\text{\rm{cl}}$, composing it with $\eta_\text{\rm{cl}} \otimes \id_{U_i \times U_j}$ and taking the trace implies the following identity: 
\be  \label{eq:Z-ij}
 Z_{ij} ~=~  \frac{1}{\Dim \Cc}
\raisebox{-30pt}{
  \begin{picture}(90,65)
   \put(0,8){\scalebox{.75}{\includegraphics{pic-S-ij.eps}}}
   \put(0,8){
     \setlength{\unitlength}{.75pt}\put(-18,-19){
     \put(98, 52)     {\scriptsize $ \Acl $}
     \put(30, 52)     {\scriptsize $U_i{\times}U_j$ }
     }\setlength{\unitlength}{1pt}}
  \end{picture}}
  \qquad \text{where} ~~
Z_{ij} = \dim_\Cb \Hom_{\CcC}(U_i{\times}U_j,\Acl) ~.
\ee
Decomposing $\Acl$ into simple objects, this gives
\be  \label{eq:S-Z-S}
  Z_{ij} = \sum_{k,l\in \Ic}  S_{ik} \, Z_{kl} \, S_{lj}^{-1}
  \qquad
  \text{where}
  ~~
  S_{ij} = s_{i,j} / \sqrt{\Dim \Cc} ~,
\ee
which in CFT terms is of course nothing but the invariance of the
torus partition function under the modular S-transformation.

The following theorem gives a simple criterion for
modular invariance.

\begin{thm}
\label{thm:modinv-dim}
Let $\Acl$ be a haploid commutative symmetric 
Frobenius algebra in $\CcC$.
\\
(i)\phantom{i} 
If $\Acl$ is modular invariant, then $\dim \Acl = \Dim\,\Cc$.
\\
(ii) 
If $\dim \Acl = \Dim\,\Cc$, then $\Acl$ is
special and modular invariant.
\end{thm}
\pf
Part (i):
Since $\Acl$ is haploid, for $i=j=0$, equation \eqref{eq:Z-ij} 
reduces to $1 = \dim \Acl / \Dim \Cc$.
\\[.3em]
Part (ii):
By the same reasoning as in the proof of (iii)$\Rightarrow$(i)
in lemma \ref{lem:modinv-simp-hapl} one shows that
$\Acl$ is special. Thus 
$m_\text{\rm{cl}} \circ \Delta_\text{\rm{cl}} = \zeta_\text{\rm{cl}} \id_{\Acl}$ for some $\zeta_\text{\rm{cl}} \neq 0$.

By \cite[thm.\,4.5]{kios}, the category $(\CcC)_{\Acl}^{loc}$ of 
local $\Acl$-modules is again a 
modular tensor category 
and $\Dim\,(\CcC)_{\Acl}^{loc}=(\Dim\,\Cc / \dim \Acl)^2$
(see \cite[prop.\,3.21 \& rem.\,3.23]{corr} for the same statement 
in the notation used here). 
Thus by assumption  
we have $\Dim (\CcC)_{\Acl}^{loc}=1$. It then follows from 
\cite[thm.\,2.3]{eno} that up to isomorphism, $(\CcC)_{\Acl}^{loc}$
has a unique simple object (namely the tensor unit).
In other words, every simple local $\Acl$-module is isomorphic to $\Acl$
(seen as a left-module over itself).

We have the following isomorphisms between morphism spaces
\cite[prop.\,4.7\,\&\,4.11]{fs-cat},
\bea
\Hom_{\Acl}(\Acl\otimes (U_i\times U_j), \Acl)
&\cong& \Hom_{\CcC}(U_i\times U_j, \Acl)~, \nn
\Hom_{\Acl}(\Acl, \Acl\otimes (U_i\times U_j))
&\cong& \Hom_{\CcC}(\Acl, U_i\times U_j)~. \nonumber
\eea
Using these to transport the bases \eqref{eq:b-dual-b} from the right to
the left, we obtain bases
$\{ b_{(ij)}^{\alpha} \}_{\alpha}$ 
of $\Hom_{\Acl}(\Acl\otimes (U_i\times U_j), \Acl)$
and $\{ b_{\beta}^{(ij)} \}_{\beta}$
of $\Hom_{\Acl}(\Acl, \Acl\otimes (U_i\times U_j))$. These can 
be expressed graphically as
\be   \label{eq:mod-inv-4}
b_{(ij)}^{\alpha} = \raisebox{-40pt}{
  \begin{picture}(60, 95)
   \put(0,8){\scalebox{.75}{\includegraphics{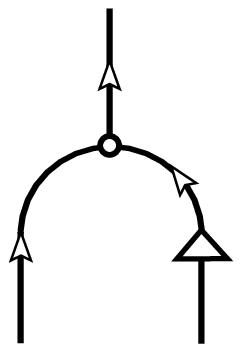}}}
   \put(0,8){
     \setlength{\unitlength}{.75pt}\put(-18,-19){
     \put(18, 8)     {\scriptsize $ \Acl $}
     \put(40, 125)  {\scriptsize $ \Acl $}
     \put(74, 70)   {\scriptsize $ \Acl $}
     \put(60,8)    {\scriptsize $ U_i\times U_j $ }
     \put(87, 46)    {\scriptsize $\alpha$}
     }\setlength{\unitlength}{1pt}}
  \end{picture}}
~~,
\quad\quad\quad
b_{\beta}^{(ij)} = \frac{\Dim \Cc}{\dim U_i \dim U_j}
\frac{1}{\zeta_\text{\rm{cl}}}
 \raisebox{-40pt}{
  \begin{picture}(60,95)
   \put(8,8){\scalebox{.75}{\includegraphics{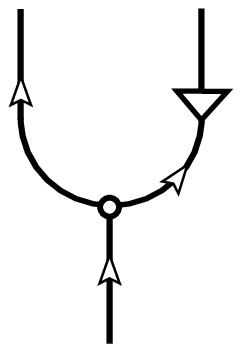}}}
   \put(0,8){
     \setlength{\unitlength}{.75pt}\put(-18,-19){
     \put(50, 8)    {\scriptsize $ \Acl $}
     \put(64, 125)  {\scriptsize $ U_i\times U_j $}
     \put(24, 125)  {\scriptsize $ \Acl $}
     \put(83, 62)   {\scriptsize $ \Acl $} 
     \put(97, 90)  {\scriptsize $\beta$}     
}\setlength{\unitlength}{1pt}}
  \end{picture}}
~~,
\ee
where the nonzero factor in $b_{\beta}^{(ij)}$ is included for
convenience. Notice that $b_{(ij)}^{\alpha} \circ b_{\beta}^{(ij)}$
is a left $\Acl$-module map.
Since $\Acl$ is simple as a left module over itself, we
have
\be
b_{(ij)}^{\alpha} \circ b_{\beta}^{(ij)} = \lambda_{\alpha\beta}\, \id_{\Acl}
\ee
for some $\lambda_{\alpha\beta} \in \mathbb{C}$. 
By computing $\text{tr}(b_{(ij)}^{\alpha} \circ b_{\beta}^{(ij)})$, 
it is easy to verify that $\lambda_{\alpha\beta} = \delta_{\alpha\beta}$. 

We will now prove the following identity:
\be   \label{eq:mod-inv-3}
\frac{1}{\zeta_\text{\rm{cl}}} \quad
  \raisebox{-50pt}{
  \begin{picture}(60, 110)
   \put(0,8){\scalebox{.75}{\includegraphics{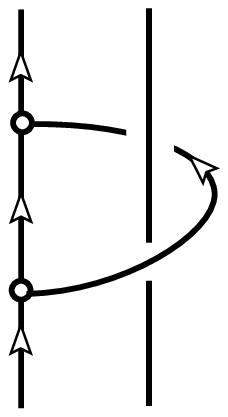}}}
   \put(0,8){
     \setlength{\unitlength}{.75pt}\put(-18,-19){
     \put(15, 8)     {\scriptsize $ \Acl $}
     \put(15, 145)  {\scriptsize $ \Acl $}
     \put(50,8)    {\scriptsize $ U_i\times U_j $ }
     \put(50,145){\scriptsize $ U_i\times U_j $ }
     \put(77, 95)  {\scriptsize $\Acl$} 
     }\setlength{\unitlength}{1pt}}
  \end{picture}}
\quad = \quad \sum_{\alpha} \,\,
\frac{\Dim \Cc}{\dim U_i \dim U_j}\, \frac{1}{\zeta_\text{\rm{cl}}}
 \raisebox{-50pt}{
  \begin{picture}(60,110)
   \put(8,8){\scalebox{.75}{\includegraphics{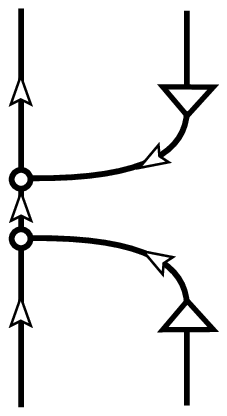}}}
   \put(0,8){
     \setlength{\unitlength}{.75pt}\put(-18,-19){
     \put(24, 8)    {\scriptsize $ \Acl $}
     \put(64, 145)  {\scriptsize $ U_i\times U_j $}
     \put(64, 8)    {\scriptsize $ U_i\times U_j $}
     \put(24, 145)  {\scriptsize $ \Acl $}
     \put(76, 85)   {\scriptsize $ \Acl $}
     \put(78, 65)   {\scriptsize $ \Acl $} 
     \put(91, 108)  {\scriptsize $\alpha$}
     \put(91, 45)    {\scriptsize $\alpha$}
}\setlength{\unitlength}{1pt}}
  \end{picture}}
  \quad.
\ee
One checks that the left hand side of this equation is an idempotent,
which we denote by 
$\zeta_\text{\rm{cl}}^{-1} P_{\Acl}^l(U_i\times U_j)$, 
cf.\ \cite[sect.\,3.1]{corr}.
By \cite[prop.\,4.1]{corr} the image 
$\text{Im}(\zeta_\text{\rm{cl}}^{-1} P_{\Acl}^l(U_i\times U_j))$ 
is a local $\Acl$-module, and
hence isomorphic to $\Acl^{\oplus N}$ for some 
$N\in \Zb_{\ge 0}$.

All left-module morphisms from $\Acl$ to $\Acl\otimes (U_i\times U_j)$
are linear combinations of the $b_{\beta}^{(ij)}$. Furthermore
one verifies that
$\zeta_\text{\rm{cl}}^{-1} P_{\Acl}^l(U_i\times U_j) \circ b_{\beta}^{(ij)} 
= b_{\beta}^{(ij)}$.
Therefore, the $b_{\beta}^{(ij)}$ describe precisely the image of
the idempotent, i.e.\ $\zeta_\text{\rm{cl}}^{-1} P_{\Acl}^l(U_i\times U_j)
= \sum_\alpha b^{(ij)}_\alpha \circ b_{(ij)}^\alpha$, which is nothing but
\eqref{eq:mod-inv-3}. Composing \eqref{eq:mod-inv-3} with 
$\zeta_\text{\rm{cl}} \cdot \eps_\text{\rm{cl}} 
\otimes \id_{U_i \times U_j}$ from 
the left (i.e.\ from the top) produces \eqref{eq:mod-inv-2}. 

In addition, since $\Acl$ is commutative symmetric Frobenius 
it satisfies $\theta_{\Acl} = \id_{\Acl}$ \cite[Prop.\,2.25]{corr}.
Altogether, this shows that $\Acl$ is modular invariant.
\epf

\begin{rema}  {\rm 
As we were writing this paper, we heard that the results
in theorem \ref{thm:modinv-dim} were obtained 
independently by Kitaev and M\"{u}ger \cite{kimu-priv}. 
}
\end{rema}

\begin{rema} \label{rem:Dim-Z00} {\rm 
Setting $i=j=0$ in \eqref{eq:Z-ij} gives the identity
$\dim \Acl = Z_{00}\, \Dim\,\Cc$ \cite[prop.\,2.3]{morita}. 
Combining this with theorem \ref{thm:modinv-dim}\,(ii) one may wonder 
if a general modular invariant commutative symmetric Frobenius algebra
$\Acl$ in $\CcC$ is isomorphic to a direct sum of simple such algebras.
However, this is not so. For example, one can take the
commutative symmetric Frobenius algebra $\Acl=\mathbb{C}[x]/x^2$
in the category of vector spaces equipped with  
the non-degenerate trace $\eps(ax+b)=a$. In this case
the modular invariance condition holds automatically, but 
$\Acl$ is clearly not a direct sum of two algebras. 
For a general modular tensor category $\Cc$, 
the algebra $\mathbb{C}[x]/x^2\,\boxtimes R(\one)$,
understood as an algebra in $\CcC$ via the braided monoidal
isomorphism 
${\mathcal V}\mbox{\sl ect}_f(\Cb) \boxtimes \CcC \rightarrow \CcC$,
provides another counter-example.}
\end{rema}

\subsection{Two definitions}\label{sec:cardy-def}

Define a morphism $P^l_A : A\rightarrow A$ for a Frobenius algebra $A$
in $\Cc$ or $\CcC$ 
as follows \cite[sect.\,2.4]{corr},
\be
  P_A^l ~=~
  \raisebox{-40pt}{
  \begin{picture}(54,80)
   \put(0,8){\scalebox{.75}{\includegraphics{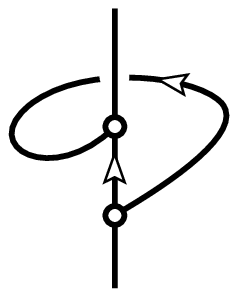}}}
   \put(0,8){
     \setlength{\unitlength}{.75pt}\put(-18,-19){
     \put(45, 10)  {\scriptsize $ A $}
     \put(45,105)  {\scriptsize $ A $}
     \put(78, 79)  {\scriptsize $ A $}
     \put(55, 65)  {\scriptsize $ m $}
     \put(36, 36)  {\scriptsize $ \Delta $}
     }\setlength{\unitlength}{1pt}}
  \end{picture}}
\labl{eq:Pl(A)-def}
If $A$ is also commutative and obeys 
$m_A \circ \Delta_A = \zeta_A\, \id_A$, we have $P_A^l = \zeta_A \, \id_A$. In particular, this holds if $A$ is commutative and special.
Using the fact that the Frobenius algebra
$R(\one)$ is commutative and normalised-special, one can check that 
$P_{R(A)}^l: R(A) \rightarrow R(A)$ takes the following form,
\be   \label{eq:P-RA-l}
  P_{R(A)}^l ~=~  \bigoplus_{i\in \Ic} 
  \raisebox{-40pt}{
  \begin{picture}(75,80)
   \put(0,8){\scalebox{.75}{\includegraphics{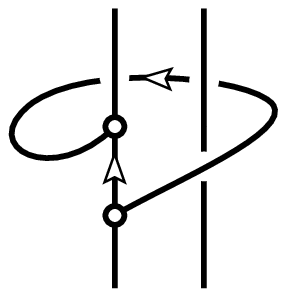}}}
   \put(0,8){
     \setlength{\unitlength}{.75pt}\put(-18,-19){
     \put(45, 10)  {\scriptsize $ A $}
     \put(45,105)  {\scriptsize $ A $}
     \put(72, 10)  {\scriptsize $ U_i^{\vee} $}
     \put(72,105)  {\scriptsize $ U_i^{\vee} $}
     \put(93, 79)  {\scriptsize $ A $}
     \put(55, 65)  {\scriptsize $ m $}
     \put(36, 36)  {\scriptsize $ \Delta $}
     }\setlength{\unitlength}{1pt}}
  \end{picture}}
\times  
 \raisebox{-40pt}{
  \begin{picture}(25,80)
   \put(8,8){\scalebox{.75}{\includegraphics{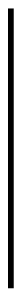}}}
   \put(0,8){
     \setlength{\unitlength}{.75pt}\put(-18,-19){
     \put(30, 10)  {\scriptsize $ U_i $}
     \put(30,105)  {\scriptsize $ U_i $}     
}\setlength{\unitlength}{1pt}}
  \end{picture}}
  ~~.
\ee
With these ingredients, we can now give the first definition of a
Cardy $\Cc|\CcC$-algebra, which was introduced in \cite[def.\ 5.14]{cardy}, 
cf.\ remark \ref{rem:cardy-rescale} below.

\begin{defn}[Cardy $\Cc|\CcC$-algebra I]  \label{def-I}
A Cardy $\Cc|\CcC$-algebra is a triple
$(\Aop|\Acl, \ico)$, where 
$(\Acl, m_\text{\rm{cl}}, \eta_\text{\rm{cl}}, \Delta_\text{\rm{cl}}, \eps_\text{\rm{cl}})$ is a 
modular invariant commutative symmetric Frobenius $\CcC$-algebra,
$(\Aop, m_\text{\rm{op}}, \eta_\text{\rm{op}}, \Delta_\text{\rm{op}}, \eps_\text{\rm{op}})$ is
a symmetric Frobenius $\mathcal{C}$-algebra, 
and $\ico: \Acl \rightarrow R(\Aop)$ an algebra homomorphism, 
such that the following conditions are satisfied: 
\\[.3em]
(i) Centre condition:  
\be  \label{eq:left-comm}
 \raisebox{-50pt}{
  \begin{picture}(80,110)
   \put(0,8){\scalebox{.75}{\includegraphics{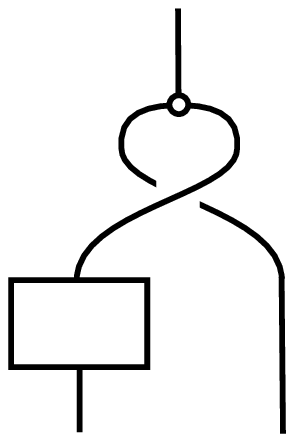}}}
   \put(0,8){
     \setlength{\unitlength}{.75pt}\put(-18,-19){
     \put(30, 10)  {\scriptsize $ \Acl $}
     \put(80, 10)  {\scriptsize $ R(\Aop) $}     
     \put(60, 150) {\scriptsize $R(\Aop)$}
     \put(2, 77)  {\scriptsize $ R(\Aop) $}
     \put(25, 50) {\scriptsize  $\ico$ }
     \put(20, 125) {\scriptsize $m_{R(\Aop)}$ }
     }\setlength{\unitlength}{1pt}}
  \end{picture}}
= \quad\quad
\raisebox{-50pt}{
  \begin{picture}(80,110)
   \put(0,8){\scalebox{.75}{\includegraphics{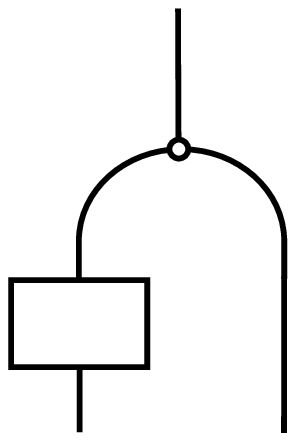}}}
   \put(0,8){
     \setlength{\unitlength}{.75pt}\put(-18,-19){
     \put(30, 10)  {\scriptsize $ \Acl $}
     \put(25, 50) {\scriptsize  $\ico$ }
     \put(80,10)  {\scriptsize $ R(\Aop) $}
     \put(-2, 73)  {\scriptsize $ R(\Aop) $}
      \put(60, 150) {\scriptsize $R(\Aop)$}
     \put(20, 110)  {\scriptsize $m_{R(\Aop)}$}
     }\setlength{\unitlength}{1pt}}
  \end{picture}}
~~.
\ee
(ii) Cardy condition: 
\be   \label{eq:cardy-CC}
  \ico \circ \ico^* ~=~
  \raisebox{-40pt}{
  \begin{picture}(54,80)
   \put(0,8){\scalebox{.75}{\includegraphics{pic-PC_l.eps}}}
   \put(0,8){
     \setlength{\unitlength}{.75pt}\put(-18,-19){
     \put(45, 10)  {\scriptsize $ R(\Aop) $}
     \put(45,105)  {\scriptsize $ R(\Aop) $}
     \put(80, 80)  {\scriptsize $ R(\Aop) $}
     }\setlength{\unitlength}{1pt}}
  \end{picture}}
~~.
\ee
\end{defn}

\begin{rema}  {\rm
\rule{0pt}{1em}\\[.1em]
(i) The name ``Cardy $\Cc|\CcC$-algebra'' in definition \ref{def-I} was chosen because many of the important ingredients were first                  
studied by Cardy: the modular invariance of the closed theory \cite{Cardy:1986ie}, the consistency of the annulus amplitude 
\cite{Cardy:1989ir}, and the bulk-boundary OPE \cite{Cardy:1991tv}.  On the other hand, the boundary-boundary OPE and the OPE analogue of the centre condition were first considered in \cite{Lew}.  
\\[.3em]
(ii) One can easily see that in the special case that
$\Cc$ is the category 
${\mathcal V}\mbox{\sl ect}_f(\Cb)$
of finite-dimensional $\Cb$-vector spaces,
a Cardy $\Cc|\CcC$-algebra gives exactly the
algebraic formulation of two-dimensional open-closed
topological field theory over $\Cb$ (cf.\ remark 6.14 in \cite{cardy}), see
\cite[sect.\,4.8]{Laz}, 
\cite[thm.\,1.1]{Mo1},
\cite[thm.\,4.5]{AN}, 
\cite[cor.\,4.3]{LP},
\cite[sect.\,2.2]{MSeg}.
When passing to a general modular tensor category $\Cc$ there are two
important differences to the two-dimensional topological field theory.
Firstly, the algebras $\Acl$ and $\Aop$ now live in {\em different} 
categories, which in particular affects the formulation of the centre
condition and the Cardy condition. Secondly, the modular invariance 
condition has to be imposed on $\Acl$. 
In the case $\Cc = {\mathcal V}\mbox{\sl ect}_f(\Cb)$,
modular invariance holds automatically. 
}
\end{rema}

\begin{defn}
A homomorphism of Cardy $\Cc|\CcC$-algebras 
$(\Aop^{(1)}|\Acl^{(1)}, \ico^{(1)}) \rightarrow 
(\Aop^{(2)}|\Acl^{(2)}, \ico^{(2)})$ is a pair $(f_\text{\rm{op}}, f_\text{\rm{cl}})$
of Frobenius algebra homomorphisms
$f_\text{\rm{op}}: \Aop^{(1)}\rightarrow \Aop^{(2)}$ and $f_\text{\rm{cl}}: \Acl^{(1)}\rightarrow \Acl^{(2)}$
such that the diagram
\be  \label{eq:cardy-hom-def}
\xymatrix{
\Acl^{(1)}  \ar[r]^{f_\text{\rm{cl}}}\ar[d]_{\ico^{(1)}} & \Acl^{(2)}
\ar[d]^{\ico^{(2)}}  \\   
R(\Aop^{(1)}) \ar[r]^{R(f_\text{\rm{op}})}  & R(\Aop^{(2)})
}  
\ee
commutes. 
\end{defn}

\begin{rema}{\rm
Since a homomorphism of Frobenius algebras is invertible 
(cf.\ lemma \ref{lem:f-star-prop}\,(iv)), a homomorphism of 
Cardy algebras is always an isomorphism.
}
\end{rema}

For a homomorphism $(f_\text{\rm{op}}, f_\text{\rm{cl}})$ 
of Cardy $\mathcal{C}|\CcC$-algebras,
using the commutativity of (\ref{eq:cardy-hom-def}) and the fact that    
$f_\text{\rm{cl}}$ and $f_\text{\rm{op}}$ are both algebra and coalgebra homomorphisms, 
it is easy to show that \eqref{eq:cardy-hom-def} commutes iff  
\be
\xymatrix{
R(\Aop^{(1)}) \ar[r]^{R(f_\text{\rm{op}})} \ar[d]_{ (\ico^{(1)})^* } & 
R(\Aop^{(2)}) \ar[d]^{ (\ico^{(2)})^* } \\
\Acl^{(1)}  \ar[r]^{f_\text{\rm{cl}}} & \Acl^{(2)}
}  
\ee
commutes. 

Let  $(\Aop|\Acl, \ico)$ be a Cardy $\mathcal{C}|\CcC$-algebra.  
Define the morphism 
\be
\itco = \hat{\chi}^{-1} (\ico) ~:~ 
T(\Acl)  \longrightarrow \Aop~. 
\labl{eq:iota-tilde-def}
Decompose $\Acl$ as $\Acl=\oplus_{n=1}^N C_n^l \times C_n^r$ such
that $C_1^l\times C_1^r = \eta_\text{\rm{cl}}(\one \times \one)$. 
We use $\ico^{(n)}$ to denote the restriction of $\ico$ 
to $C_n^l\times C_n^r$ and $\itco^{(n)}$ to denote 
the restriction of $\itco$ to $C_n^l\otimes C_n^r$. 
We introduce the following graphical notation: 
\be   \label{iota-*-def}
\itco^{\, (n)} ~=~
  \raisebox{-30pt}{
  \begin{picture}(54,73)
   \put(0,8){\scalebox{.75}{\includegraphics{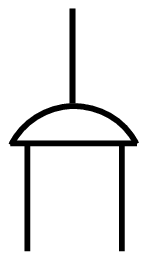}}}
   \put(0,8){
     \setlength{\unitlength}{.75pt}\put(-18,-19){
     \put(19, 8)  {\scriptsize $ C_n^l $}
     \put(51, 8)  {\scriptsize $ C_n^r $}
     \put(35, 98)  {\scriptsize $ \Aop $}
     }\setlength{\unitlength}{1pt}}
  \end{picture}}
~~.
\ee
By (\ref{eq:hat-chi}), $\ico$ can be expressed in terms of 
$\itco$ as follows:
\be   \label{eq:iota-clop-graph}
\setlength{\unitlength}{1pt}
\ico =  \bigoplus_{n=1}^N \bigoplus_{i\in \Ic} \sum_{\alpha}\,\, 
\raisebox{-45pt}{
  \begin{picture}(60,90)
   \put(0,8){\scalebox{.75}{\includegraphics{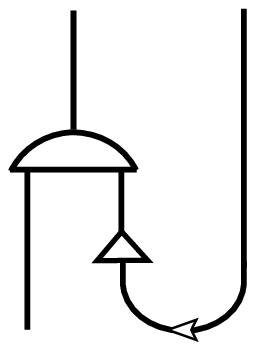}}}
   \put(0,8){
     \setlength{\unitlength}{.75pt}\put(-18,-11){
     \put(23, 2)    {\scriptsize $ C_n^l $}
     \put(62, 46)  {\scriptsize $ C_n^r $}
     \put(37, 110){\scriptsize $\Aop$}
     \put(38, 35)  {\scriptsize $\alpha$} 
     \put(87,110) {\scriptsize $U_i^{\vee}$}
     }\setlength{\unitlength}{1pt}}
  \end{picture}}
\,\, \times \,\, 
\raisebox{-45pt}{
  \begin{picture}(35,90)
  \put(0,8){\scalebox{.75}{\includegraphics{pic-hatchi-R-2.eps}}}
   \put(0,8){
     \setlength{\unitlength}{.75pt}\put(-18,-11){
     \put(37,57)     {\scriptsize $ \alpha $}
     \put(23,2)       {\scriptsize $ C_n^r $}
     \put(20,110)   {\scriptsize $ U_i $}
  }\setlength{\unitlength}{1pt}}
  \end{picture}} 
\hspace{0.7cm}
  ~~.
\ee

\begin{lemma}  \label{lemma:left-comm}
The centre condition \eqref{eq:left-comm} is equivalent to 
the following condition in $\Cc$,
\be   \label{eq:comm-C}
\setlength{\unitlength}{1pt}
\raisebox{-45pt}{
  \begin{picture}(60,90)
   \put(0,8){\scalebox{.75}{\includegraphics{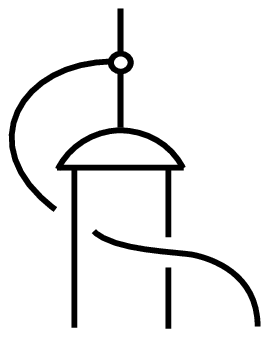}}}
   \put(0,8){
     \setlength{\unitlength}{.75pt}\put(-18,-11){
     \put(33, 2)     {\scriptsize $ C_n^l $}
     \put(63, 2)     {\scriptsize $ C_n^r $}
     \put(45, 110) {\scriptsize $\Aop$}
     \put(85, 2)     {\scriptsize $ \Aop$} 
     \put(58, 78)   {\scriptsize $ \Aop$}
     }\setlength{\unitlength}{1pt}}
  \end{picture}}
~=~ 
\raisebox{-45pt}{
  \begin{picture}(35,90)
  \put(0,8){\scalebox{.75}{\includegraphics{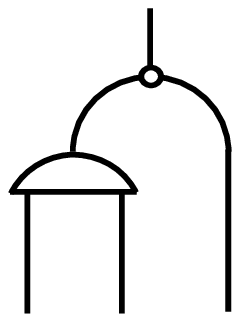}}}
   \put(0,8){
     \setlength{\unitlength}{.75pt}\put(-18,-11){
     \put(23,2)       {\scriptsize $ C_n^l $}
     \put(50,2)       {\scriptsize $ C_n^r $}
     \put(58,110)   {\scriptsize $ \Aop $}
     \put(80,2)       {\scriptsize $ \Aop $}
     \put(23,75)     {\scriptsize $ \Aop $}
  }\setlength{\unitlength}{1pt}}
  \end{picture}} 
\qquad \qquad
\text{for}~n=1, \dots, N~~. 
\ee
\end{lemma}
\pf
First, insert (\ref{eq:iota-clop-graph}) and the definition \eqref{eq:FA-alg-st} of $m_{R(\Aop)}$ into (\ref{eq:left-comm}).  Then apply the commutativity of $R(\one)$ to the left hand side of (\ref{eq:left-comm}). The equivalence between (\ref{eq:left-comm}) and (\ref{eq:comm-C}) follows immediately. 
\epf

\begin{rema} {\rm 
The centre condition (\ref{eq:comm-C})
is very natural from the open-closed conformal field theory point of view. Correlators on 
the upper half plane are expressed in terms of conformal blocks on the full complex 
plane. The objects $C_n^l$ and $C_n^r$ are associated to the field insertion at a
point $z$ in the upper half plane and at the complex conjugate point $\bar z$ in the 
lower half plane, respectively. The object $\Aop$ corresponds to a field inserted at a point $r$ on the 
real axis. The centre condition (\ref{eq:comm-C}) simply says that the correlation 
functions in the disjoint domains $|z|>r>0$ and $r>|z|>0$ are analytic continuations of 
each other, see \cite[prop. 1.18]{ocfa}. 
}
\end{rema}

Recall that we define $\itco^*: \Aop \rightarrow T(\Acl)$ 
as in (\ref{f-star-def}).  We introduce the graphical notation
\be   \label{eq:tilde-iota}
\setlength{\unitlength}{1pt}
\itco^* = \bigoplus_{n=1}^N \,\,
\raisebox{-35pt}{
  \begin{picture}(40,75)
   \put(0,8){\scalebox{.75}{\includegraphics{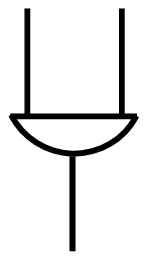}}}
   \put(0,8){
     \setlength{\unitlength}{.75pt}\put(-18,-11){
     \put(20, 90)     {\scriptsize $ C_n^l $}
     \put(50, 90)     {\scriptsize $ C_n^r $}
     \put(35, 2)     {\scriptsize $\Aop$}
     }\setlength{\unitlength}{1pt}}
  \end{picture}}
~=~~ \bigoplus_{n=1}^N \,\,
\raisebox{-35pt}{
  \begin{picture}(85,70)
  \put(0,8){\scalebox{.75}{\includegraphics{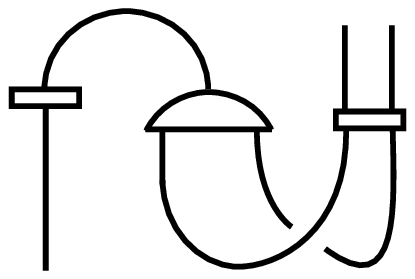}}}
   \put(0,8){
     \setlength{\unitlength}{.75pt}\put(-18,-11){
     \put(130,90)       {\scriptsize $ C_n^l $}
     \put(150,90)       {\scriptsize $ C_n^r $}
     \put(47,2)           {\scriptsize $ \Aop $}
     \put(13,60)         {\scriptsize $ \Phi_{\Aop} $}
     \put(160, 52)      {\scriptsize $ T(\Phi_{\Acl}^{-1}) $ }
  }\setlength{\unitlength}{1pt}}
  \end{picture}} 
\hspace{0.7cm}
  \qquad\quad,
\ee
where the second equality follows from  for 
(\ref{eq:Frob-sym-cond}), (\ref{eq:Phi-I}) 
and (\ref{eq:I-I'-T-c}).

\begin{lemma}  \label{lemma:cardy-cond}
The Cardy condition (\ref{eq:cardy-CC}) 
is equivalent to the following identity in $\Cc$: 
\be   \label{eq:cardy-C}
\bigoplus_{n=1}^N \frac{\Dim \Cc }{\dim U_i}
      \, \sum_\alpha  
\quad
\raisebox{-70pt}{
  \begin{picture}(70,130)
   \put(0,8){\scalebox{.75}{\includegraphics{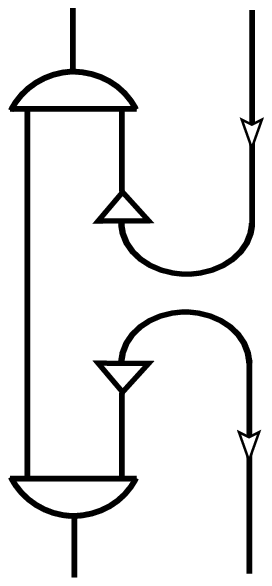}}}
   \put(0,8){
     \setlength{\unitlength}{.75pt}\put(-18,-11){
     \put(9, 95)     {\scriptsize $ C_n^l $}
     \put(60, 54)     {\scriptsize $ C_n^r $}
     \put(60, 134)   {\scriptsize $ C_n^r $}
     \put(40, 118)   {\scriptsize $ \alpha $}
     \put(40, 70)   {\scriptsize $ \alpha $}
     \put(35, 2)       {\scriptsize $\Aop$}
     \put(35, 185)   {\scriptsize $\Aop$}
     \put(85, 185)   {\scriptsize $ U_i^{\vee} $}
     \put(85, 2)       {\scriptsize $ U_i^{\vee} $}
     }\setlength{\unitlength}{1pt}}
  \end{picture}}
=~
\raisebox{-55pt}{
  \begin{picture}(85,130)
  \put(0,8){\scalebox{.75}{\includegraphics{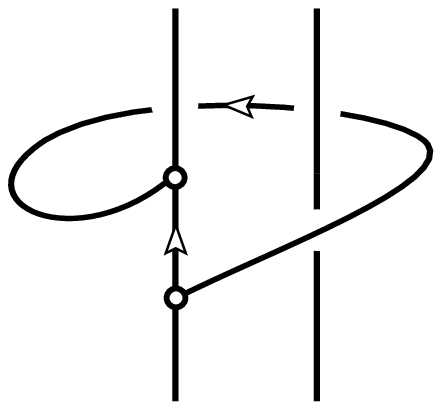}}}
   \put(0,8){
     \setlength{\unitlength}{.75pt}\put(-18,-11){ 
     \put(140,90)       {\scriptsize $ \Aop $}
     \put(60,2)           {\scriptsize $ \Aop $}
     \put(60,135)       {\scriptsize $ \Aop $}
     \put(100, 2)        {\scriptsize $ U_i^{\vee} $}
     \put(100, 135)    {\scriptsize $ U_i^{\vee} $}
     \put(40, 42)      {\scriptsize $\Delta_\text{\rm{op}}$ }
     \put(72, 77)     {\scriptsize $m_\text{\rm{op}}$ }
  }\setlength{\unitlength}{1pt}}
  \end{picture}} 
\qquad \quad \text{for~all}~ i\in \mathcal{I}~.
\ee
\end{lemma}
\pf
By \eqref{eq:check-chi}, \eqref{eq:P-RA-l} 
and \eqref{eq:iota-clop-graph}, it is easy to 
see that \eqref{eq:cardy-C} is equivalent to the following identity:
\be
\ico \circ \check{\chi}(\itco^*) = P_{R(\Aop)}^l~.
\ee
Therefore, it is enough to show that 
\be  \label{eq:iota-iota}
\check{\chi}(\itco^*) = \ico^*~. 
\ee
We have
\be
\check\chi^{-1}( \ico^* )
\overset{(1)}=
T( \ico^* ) \circ \check\delta
\overset{(2)}=
T( \ico )^* \circ \check\delta^{**}
\overset{(3)}=
\big(\hat\rho \circ T( \ico )\big)^*
\overset{(4)}=
\big( \itco \big)^* ~.
\ee
In step (1) we use the expression \eqref{eq:chi-inv} for $\check\chi^{-1}$, step (2) follows from lemma \ref{lem:f-star-prop}\,(v) and lemma \ref{lem:Fstar-starF}. Step (3) is lemma \ref{lem:f-star-prop}\,(i) and lemma \ref{lem:delta-star-rho}, and finally step (4) amounts to substituting \eqref{eq:chi-inv} and \eqref{eq:iota-tilde-def}. Acting with $\check\chi$ on both sides of the above equality produces \eqref{eq:iota-iota}.
\epf

\medskip

Combining lemmas \ref{lemma:left-comm} and \ref{lemma:cardy-cond}, and proposition \ref{prop:f-til-f}, we obtain the following equivalent definition of Cardy $\Cc|\CcC$-algebra (recall the graphical
notation \eqref{iota-*-def} for $\itco$ and \eqref{eq:tilde-iota} for $\itco^*$). 

\begin{defn}[Cardy $\Cc|\CcC$-algebra II]   \label{def-II} 
A Cardy $\Cc$-algebra is a triple $(\Aop|\Acl, \itco)$, where
$\Acl$ is a commutative symmetric
Frobenius $\CcC$-algebra satisfying property \eqref{eq:mod-inv-2}
with $f=m_\text{\rm cl}$, $\Aop$ is  
a symmetric Frobenius $\mathcal{C}$-algebra, and
$\itco: T(\Acl) \rightarrow \Aop$ is an algebra homomorphism 
satisfying the conditions 
\eqref{eq:comm-C} and \eqref{eq:cardy-C}. 
\end{defn}

\begin{rema}\label{rem:cardy-rescale} {\rm 
Up to a choice of normalisation, definition \ref{def-II} 
is the same as the original one in \cite[def.\,6.13]{cardy}. 
The difference between the two definitions is 
the factor $\Dim\,\Cc / \dim U_i$ on the left hand side
of \eqref{eq:cardy-C}, which in \cite[def.\,6.13]{cardy} is
given by $\sqrt{\Dim\,\Cc} / \dim U_i$. The two definitions
are related by rescaling the coproduct $\Delta_\text{\rm{cl}}$ 
and counit $\eps_\text{\rm{cl}}$ of $\Acl$ by 
$1/\sqrt{\Dim\,\Cc}$ and $\sqrt{\Dim\,\Cc}$, 
respectively. 
We chose the convention in \eqref{eq:cardy-C} to remove all
dimension factors from the expression \eqref{eq:cardy-CC}
for the Cardy condition.
}
\end{rema}

\subsection{Uniqueness and existence theorems}
\label{sec:cardy-reconstr}

In this subsection we investigate the structure of Cardy
algebras. We start with the following proposition,
which, when combined with the results of part II, 
provides an alternative proof of 
\cite[prop.\,4.22]{unique}.

\begin{prop}  \label{prop:Aop-special}
Let $(\Aop|\Acl, \ico)$ be a Cardy $\Cc|\CcC$-algebra.
If $\Acl$ is simple and $\dim \Aop \neq 0$,  
then $\Aop$ is simple and special. 
\end{prop}
\pf
By remark \ref{rem:Dim-Z00}, we have 
$\dim \Acl = Z_{00} \, \Dim\,\Cc \neq 0$,
and by lemma \ref{lem:modinv-simp-hapl}, $\Acl$ 
is therefore haploid. 
Restricting the Cardy condition \eqref{eq:cardy-C} 
to the case $U_i=\one$ and composing
both sides with $\eps_\text{\rm{op}}$ from the left, we see that $\eps_\text{\rm{op}}$ 
kills all terms associated to 
$U_j \times \one\in \Acl$ in the sum except for
a single $\one \times \one$ term. 
Thus we obtain the following identity,
\be \label{eq:cardy-a-e-1}
\beta \, \eps_\text{\rm{op}}  = 
   \tilde d_{\Aop} \circ (m_\text{\rm{op}} \otimes \id_{\Aop^\vee})
    \circ (\id_{\Aop} \otimes b_{\Aop}) ~,
\ee
where $\beta \in \mathbb{C}$.
Composing with $\eta_\text{\rm{op}}$ from the right in turn implies that
$\beta \eps_\text{\rm{op}} \circ \eta_\text{\rm{op}} = \dim \Aop$, which
is nonzero by assumption. Thus also $\beta \neq 0$ and
$\eps_\text{\rm{op}}$ is a nonzero multiple of the morphism on the
right hand side of \eqref{eq:cardy-a-e-1}. By 
\cite[lem.\,3.11]{tft1}, $\Aop$ is special. 

Since $\Aop$ is a special Frobenius algebra, $\Aop$ is
semi-simple as an $\Aop$-bimodule 
(apply \cite[prop.\,5.24]{fs-cat} to $\Aop$ tensored with 
its opposite algebra). Suppose $\Aop$ is not simple, so 
that we can write
$\Aop= \Aop^{(1)}\oplus \Aop^{(2)}$ for nonzero
$\Aop$-bimodules $\Aop^{(1)}$ and $\Aop^{(2)}$. 
We denote the canonical 
embeddings and projections associated to this decomposition 
as $\iota_{1,2}$ and $\pi_{1,2}$.  We have the identities
\be  
m_\text{\rm{op}} \circ (\iota_1 \otimes \iota_2) = 0
\quad , \quad
\eps_\text{\rm{op}} \circ \eta_\text{\rm{op}}= 
\sum_{i=1}^2 \eps_\text{\rm{op}}\circ \iota_i \circ \pi_i 
\circ \eta_\text{\rm{op}}~.
\labl{eq:Aop-simple-aux3}
The first identity follows since $\pi_1 \circ m_\text{\rm{op}} \circ (\iota_1 \otimes \iota_2) = 0$ (as $m_\text{\rm{op}}$ gives the left action of $\Aop$ on $\Aop$ and hence it preserves $\Aop^{(2)}$), and similarly $\pi_2 \circ m_\text{\rm{op}} \circ (\iota_1 \otimes \iota_2) = 0$. The second identity is just the completeness of $\iota_{1,2}, \pi_{1,2}$. 

Since $\eps_\text{\rm{op}} \circ \eta_\text{\rm{op}} \neq 0$,
without losing generality we can assume 
$\eps_\text{\rm{op}}\circ \iota_1 \circ \pi_1 \circ \eta_\text{\rm{op}} \neq 0$.
Using that $\pi_2$ is a bimodule map we compute
\be\begin{array}{l}\displaystyle  
\pi_2 \circ \big[\mbox{LHS of \eqref{eq:cardy-C}}\big]_{U_i=\one} 
\circ \iota_1 \circ \pi_1 \circ \eta_\text{\rm{op}} 
= \pi_2 \circ \big[\mbox{RHS of (\ref{eq:cardy-C})}\big]_{U_i=\one} 
\circ \iota_1 \circ \pi_1
\circ \eta_\text{\rm{op}} 
\enl
=     P_{\!A_{\rm{op}} ^ {(2)} }^l      
\circ \pi_2 \circ \iota_1 \circ \pi_1
\circ \eta_\text{\rm{op}} = 0~. 
\eear
\labl{prop-equ-2}
On the other hand, using that $\Acl$ is haploid, 
that $\itco$ is an algebra map, 
and that $\itco^{\,*}$ is a coalgebra map, one can check that 
the left hand side of (\ref{prop-equ-2}) is equal to
$\lambda (\eps_\text{\rm{op}}\circ \iota_1 \circ \pi_1 \circ \eta_\text{\rm{op}}) \,
\pi_2 \circ \eta_\text{\rm{op}} $ for some $\lambda\neq 0$. This implies that
$\pi_2\circ \eta_\text{\rm{op}}=0$. 
Thus 
$\eta_\text{\rm{op}} = \sum_i \iota_i \circ \pi_i \circ \eta_\text{\rm{op}}
= \iota_1 \circ \pi_1 \circ \eta_\text{\rm{op}}$.
Hence, we have
\be
0\neq \pi_2 \circ \iota_2 =
\pi_2 \circ m_\text{\rm{op}} \circ (\eta_\text{\rm{op}} \otimes \iota_2) 
=  \pi_2 \circ m_\text{\rm{op}} \circ 
((\iota_1 \circ \pi_1 \circ \eta_\text{\rm{op}}) \otimes \iota_2)  ~.
\ee
However, the right hand side is zero by 
\eqref{eq:Aop-simple-aux3}. This
is a contradiction and hence $\Aop$ must be simple. 
\epf

\medskip

To formulate the next theorem we need the notion of the full centre of an algebra \cite[def.\,4.9]{unique}. Recall that an algebra $A$ in a braided tensor category has  
a left centre and a right centre \cite{vz,ost}, 
both of which are sub-algebras of $A$. Of these two, 
we will only need the left centre. 
The following definition is \cite[def.\,2.31]{corr}, 
which in our setting is equivalent to that of \cite{vz,ost}. 

\begin{defn} 
Let $A$ be a symmetric special Frobenius algebra such that
$m_A \circ \Delta_A = \zeta_A\,\id_A$. 
\\
(i)\phantom{i}~The {\em left centre} $C_l(A)$ of 
$A$ is the image of the idempotent $\zeta_A^{-1} \, P^l_A$.
\\
(ii)~The {\em full centre} $Z(A)$ is $C_l(R(A))$.
\end{defn}

That $\zeta_A^{-1} \, P^l_A$ is an idempotent follows from
\cite[lem.\,5.2]{tft1} when keeping track of the factors $\zeta_A$ 
(\cite{tft1} assumes normalised-special, i.e.\ $\zeta_A=1$).
Note that
$C_l(A)$ is again an object of $\Cc$, while $Z(A)$ is an object 
of $\CcC$. Let $e_A^l : C_l(A) \rightarrow A$ be the embedding of
$C_l(A)$ into $A$. The left centre is in fact the maximal subobject
of $A$ such that
\be
  m_A \circ c_{A,A} \circ (e_A^l \otimes \id_A) = m_A ~,
\labl{eq:left-centre-prop}
see \cite[lem.\,2.32]{corr}. 
This observation explains the name left centre and
also makes the connection to \cite[def.\,15]{ost}.

The full centre is by definition the image of the
idempotent $\zeta_A^{-1} P_{R(A)}^l: R(A) \rightarrow R(A)$. 
Since $\CcC$ is abelian, the idempotent splits and we obtain
the embedding and restriction morphisms
\be
  e : Z(A) \hookrightarrow R(A)
  \quad \text{and} \quad
  r : R(A) \twoheadrightarrow Z(A)
\labl{eq:eZ-rZ-def}
which obey $r \circ e = \id_{Z(A)}$ and
$e \circ r = \zeta_A^{-1} P_{R(A)}^l$.
It follows from proposition \ref{prop:symm-R-symm} and
\cite[prop.\,2.37]{corr} that $Z(A)$ is           
a commutative symmetric Frobenius algebra in $\CcC$ with
structure morphisms\footnote{
  The normalisation of product and unit is the standard one.
  The factors in the coproduct and counit have to be included
  in order for $(A|Z(A), e)$ to be a Cardy algebra,
  see theorem \ref{thm:reconst} below. The normalisation of
  the counit enters the Cardy condition \eqref{eq:cardy-CC} through
  the definition of $(\,\cdot\,)^*$.}
\be \begin{array}{ll}\displaystyle
m_{Z(A)} = r \circ m_{R(A)} \circ 
(e \otimes e) ~~, \quad \etb 
\eta_{Z(A)} = r \circ \eta_{R(A)},  \enl
\Delta_{Z(A)} = \zeta_A \cdot 
(r \otimes r) \circ 
\Delta_{R(A)} \circ e~~, \quad 
\etb \eps_{Z(A)} = \zeta_A^{-1} \cdot 
\eps_{R(A)} \circ e ~.
\eear
\labl{eq:ZA-alg}
Moreover,
if $A$ is simple then $Z(A)$ is simple, and if
$A$ is simple and $\dim A\neq0$, then $Z(A)$ is simple and special. 
The normalisation of the counit is such that
\be  \label{eq:ZA-counit-norm}
\eps_{Z(A)} \circ \eta_{Z(A)} = \zeta_A^{-2}\,\dim A \, \Dim \Cc~.
\ee

\begin{thm}  \label{thm:reconst}
Let $A$ be a special symmetric Frobenius 
$\Cc$-algebra. Then $(A|Z(A), e)$ is a Cardy $\Cc|\CcC$-algebra. 
\end{thm}

The proof of this theorem makes use the following two lemmas.

\begin{lemma}  \label{lem:ZA-e-star}
$e: Z(A) \hookrightarrow R(A)$ is an algebra map, and
$e^* = \zeta_A \cdot r$.
\end{lemma}

\pf
It follows from \cite[lem.\,2.29]{corr} (or by direct calculation, 
using in particular
$m_{R(A)} \circ \Delta_{R(A)} = \zeta_A \id_{R(A)}$) that
\be
  m_{R(A)} \circ (e  \otimes e )
  = \zeta_A^{-1} \cdot
  P^l_{R(A)} \circ m_{R(A)} \circ (e  \otimes e ) ~.
\labl{eq:leave-centre-proj}
Substituting $e  \circ r  = \zeta_A^{-1} P_{R(A)}^l$
shows that $e $ is compatible with multiplication.
For the unit one finds
\be
  e  \circ \eta_{Z(A)}
  = e  \circ r  \circ \eta_{R(A)}
  = \zeta_A^{-1} \, P_{R(A)}^l \circ \eta_{R(A)}
  = \eta_{R(A)}
\ee
Thus $e $ is an algebra map. For the second
statement one computes
\be
 e^* 
 \overset{(1)}{=} 
 \zeta_A \,\,
 \raisebox{-60pt}{
  \begin{picture}(100,130)
   \put(0,8){\scalebox{.75}{\includegraphics{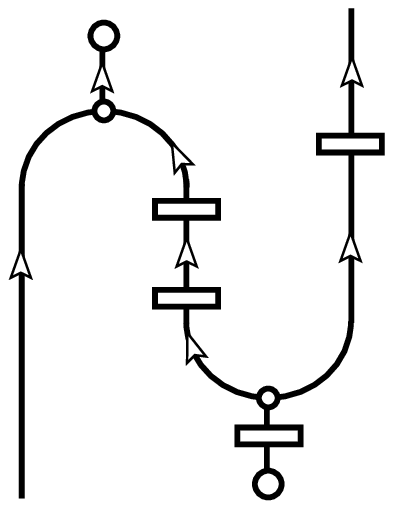}}}
   \put(0,8){
     \setlength{\unitlength}{.75pt}\put(-18,-19){
     \put(15, 10)       {\scriptsize $ R(A) $}
     \put(130,120)    {\scriptsize $ r $}
     \put(123,90)      {\scriptsize $ R(A) $}
     \put(74, 112)     {\scriptsize $ R(A) $}
     \put(74,  88)      {\scriptsize $ Z(A) $} 
     \put(48, 47)       {\scriptsize $ R(A) $}
     \put(104,168)    {\scriptsize $ Z(A) $}
     \put(51, 77)       {\scriptsize $ r $}
     \put(51, 102)     {\scriptsize $ e $}
     \put(76,37)        {\scriptsize $ e $}
     }\setlength{\unitlength}{1pt}}
  \end{picture}}
 \overset{(2)}{=} 
\,\,
  \raisebox{-60pt}{
  \begin{picture}(100,130)
   \put(0,8){\scalebox{.75}{\includegraphics{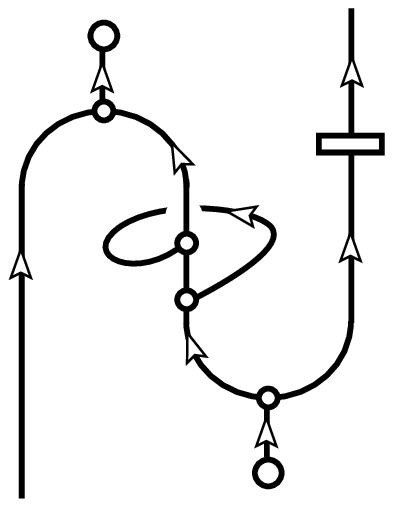}}}
   \put(0,8){
     \setlength{\unitlength}{.75pt}\put(-18,-19){
     \put(15, 10)       {\scriptsize $ R(A) $}
     \put(130,120)    {\scriptsize $ r $}
     \put(123,90)      {\scriptsize $ R(A) $}
     \put(74, 112)     {\scriptsize $ R(A) $}
     \put(58, 40)       {\scriptsize $ R(A) $}
     \put(104,168)    {\scriptsize $ Z(A) $ }
     }\setlength{\unitlength}{1pt}}
  \end{picture}}
 \overset{(3)}{=} 
\, \,
\raisebox{-38pt}{
  \begin{picture}(45,80)
   \put(0,8){\scalebox{.75}{\includegraphics{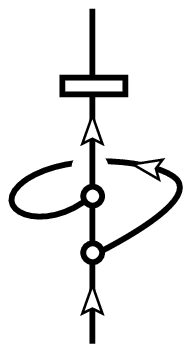}}}
   \put(0,8){
     \setlength{\unitlength}{.75pt}\put(-18,-19){
     \put(28, 10)       {\scriptsize $ R(A) $}
     \put(24,93)    {\scriptsize $ r $}
     \put(55, 80)       {\scriptsize $ R(A) $}
     \put(28,123)    {\scriptsize $ Z(A) $ }
     }\setlength{\unitlength}{1pt}}
  \end{picture}}
 \overset{(4)}{=} 
\,\, \zeta_A \cdot r 
~~,
\ee
where in (1) the definitions \eqref{f-star-def} and 
\eqref{eq:ZA-alg} have been substituted, step (2) is
$e  \circ r  = \zeta_A^{-1} P_{R(A)}^l$,
step (3) uses that $R(A)$ is symmetric Frobenius, and
step (4) is again
$e  \circ r = \zeta_A^{-1} P_{R(A)}^l$.
\epf

\begin{lemma} \label{lemma:mod-inv-RA}
Let $A$ be a symmetric Frobenius algebra in $\Cc$. 
The morphism
\be
  P^l_{R(A)} \circ m_{R(A)} \circ
  \big( P^l_{R(A)} \otimes P^l_{R(A)} \big)  
  ~ : ~ R(A) \otimes R(A) \longrightarrow R(A)
\labl{eq:RA-S-inv}
is S-invariant.
\end{lemma}

The proof of this lemma is a slightly lengthy
explicit calculation and has been deferred to appendix 
\ref{app:RA-Sinv}.

\medskip

\noindent{\it Proof of theorem \ref{thm:reconst}.}\hspace{2ex}
That $e$ is an algebra map was proved in lemma
\ref{lem:ZA-e-star}. The centre condition 
\eqref{eq:left-comm} holds by
property \eqref{eq:left-centre-prop} of the left centre.
The Cardy condition \eqref{eq:cardy-CC} also is an
immediate consequence of lemma \ref{lem:ZA-e-star},
\be
  \ico \circ \ico^* = e \circ (\zeta_A \, r)
  = P_{R(A)}^l ~.
\ee  
The full centre $Z(A)$ is a commutative symmetric Frobenius algebra.
It remains to prove modular invariance. That 
$\theta_{Z(A)} = \id_{Z(A)}$ is implied by commutativity and
symmetry of $Z(A)$ \cite[prop.\,2.25]{corr}. 
The S-invariance condition \eqref{eq:mod-inv} follows from 
lemma \ref{lemma:mod-inv-RA}: In 
\eqref{eq:RA-S-inv} substitute 
$P_{R(A)}^l = \zeta_A \, e \circ r$
and then put the resulting morphism into 
\eqref{eq:mod-inv}. Compose the resulting equation
with $e \otimes \id_W$ from the right 
(i.e.\ from the bottom) and substitute the definition
\eqref{eq:ZA-alg} of $m_{Z(A)}$. This results in the statement
that $m_{Z(A)}$ is S-invariant.
\epf

\medskip

The following theorem is analogous to 
\cite[prop.\,2.9]{Longo:2004fu} and
\cite[thm.\,4.26]{unique}, which, roughly speaking, answer the question
under which circumstances the restriction of a
two-dimensional conformal field theory to 
the boundary already determines the entire conformal field theory. The
first work is set in Minkowski space and uses operator
algebras and subfactors, while the second work is set in 
Euclidean space and uses modular tensor categories.

\begin{thm}  \label{thm:unique}
Let $(A|\Acl,\ico)$ be a Cardy $\Cc|\CcC$-algebra such that
$\dim A\neq 0$ and $\Acl$ is simple. Then $A$ is special and
$(A|\Acl, \ico) \cong (A|Z(A), e)$ as Cardy algebras.
\end{thm}

\pf
By proposition \ref{prop:Aop-special}, $A$ is simple and special.
Since $\Acl$ is simple, the algebra map 
$\ico : \Acl \rightarrow R(A)$ is either zero or a 
monomorphism. But $\ico \circ \eta_\text{\rm{cl}} = \eta_{R(A)}$,
and so $\ico \neq 0$. Thus $\ico$ is monic. 
By lemma \ref{lem:f-star-prop}\,(ii), $\zeta_A^{-1}\, \ico^*$ is epi.
The Cardy condition \eqref{eq:cardy-CC} implies
\be
  \ico \circ \zeta_A^{-1}\, \ico^* = \zeta_A^{-1} P^l_{R(A)} = 
  e \circ r ~.
\ee
Composing this with $e \circ r$ from the left yields 
$e \circ r \circ \ico \circ \zeta_A^{-1}\, \ico^* = e \circ r = 
\ico \circ \zeta_A^{-1}\, \ico^*$. Since $\zeta_A^{-1}\, \ico^*$ is epi, we have
\be
  e \circ r \circ \ico = \ico ~.
\labl{eq:unique-aux1}
Actually, (\ref{eq:unique-aux1}) also follows from (\ref{eq:left-comm}) and 
specialness of $R(A)$. We will prove that
\be
  (f_\text{\rm{op}}, f_\text{\rm{cl}})
  : (A|\Acl,\ico) \longrightarrow (A|Z(A), e)
\quad \text{where} ~
  f_\text{\rm{op}} = \id_A 
  ~,~
  f_\text{\rm{cl}} = r \circ \ico
  ~,
\ee
is an isomorphism of Cardy algebras. 
\\[.3em]
{\it $f_\text{\rm{cl}}$ is an algebra map:}
Compatibility with the units follows since $\ico$ is an algebra map,
\be
f_\text{\rm{cl}} \circ \eta_\text{\rm{cl}}
=
r \circ \ico \circ \eta_\text{\rm{cl}}
=
r \circ \eta_{R(A)}
=
\eta_{Z(A)} ~.
\labl{eq:unique-fcl-unit}
Compatibility with the multiplication also follows since $\ico$
is an algebra map,
\be\begin{array}{l}\displaystyle
m_{Z(A)} \circ (f_\text{\rm{cl}} \otimes f_\text{\rm{cl}})
= 
r \circ m_{R(A)} \circ 
(e \oti e) \circ (r \oti r) \circ (\ico \oti \ico)
\enl
= 
r \circ m_{R(A)} \circ 
(\ico \oti \ico)
= 
r \circ  \ico \circ m_\text{\rm{cl}}
= 
f_\text{\rm{cl}} \circ m_\text{\rm{cl}} ~,
\eear\labl{eq:unique-fcl-mult}
where in the second step we used \eqref{eq:unique-aux1}.
\\[.3em]
{\it $f_\text{\rm{cl}}$ is an isomorphism:}
As above, since $f_\text{\rm{cl}}$ is an algebra map and since
$\Acl$ is simple, $f_\text{\rm{cl}}$ has to be monic.
By lemma \ref{lem:ZA-e-star}, 
$r^* = \zeta_A^{-1} \, e$. Thus
$f_\text{\rm{cl}}^* = \ico^* \circ r^* = \zeta_A^{-1}\, \ico^* \circ e$ and
\be
 f_\text{\rm{cl}} \circ f_\text{\rm{cl}}^* =
 r \circ \ico \circ \zeta_A^{-1}\, \ico^* \circ e
 = r \circ e \circ r \circ e
 = \id_{Z(A)} ~,
\labl{eq:unique-aux2}
and so $f_\text{\rm{cl}}$ is also epi, and hence iso.
\\[.3em]
{\it $f_\text{\rm{cl}}$ is a coalgebra map:}
Since $f_\text{\rm{cl}}$ is an algebra map, so is $f_\text{\rm{cl}}^{-1}$. By
\eqref{eq:unique-aux2}, $f_\text{\rm{cl}}^{-1} = f_\text{\rm{cl}}^*$ 
and by lemma \ref{lem:f-star-prop}\,(iii) 
this implies that $f_\text{\rm{cl}}$ is a also coalgebra map.
\\[.3em]
{\it The diagram \eqref{eq:cardy-hom-def} commutes:}
Commutativity of \eqref{eq:cardy-hom-def} is equivalent to
$e \circ f_\text{\rm{cl}} = \ico$, which holds
by \eqref{eq:unique-aux1}.
\epf

\medskip

Let $A$ be a special symmetric Frobenius algebra.
So far we have seen that $(A,Z(A),e)$ is a Cardy algebra,
and that all Cardy algebras with $\Aop=A$ and simple $\Acl$
are of this form. It is now natural to ask if every
simple $\Acl$ does occur as part of a Cardy algebra. 
The following theorem provides an affirmative answer.
Recall that for an $A$-left module $M$, the object
$M^\vee \otimes_A M$ is an algebra 
(see e.g.\ \cite[lem.\,4.2]{morita}).

\begin{thm}\label{thm:open-exists}
If $\Acl$ is a simple modular invariant 
commutative symmetric Frobenius
$\CcC$-algebra, then there exist a simple 
special symmetric Frobenius $\Cc$-algebra $A$ and 
a morphism $\ico: \Acl \rightarrow R(A)$ such that  
\\[.3em]
(i)\phantom{ii}~$\Acl \cong Z(A)$ as Frobenius algebras;
\\[.3em]
(ii)\phantom{i}~$(A|\Acl, \ico)$ is a 
Cardy $\Cc|\CcC$-algebra;
\\[.3em]
(iii)~$T(\Acl) \cong \oplus_{\kappa\in \mathcal{J}}\, 
M_{\kappa}^{\vee} \otimes_{A} M_{\kappa}$ 
as algebras, 
where $\{ M_{\kappa}\}_{\kappa \in \mathcal{J}}$ 
is a set of representatives of the isomorphism classes of
simple $A$-left modules.  
\end{thm}

\pf
By remark \ref{rem:Dim-Z00}, we have 
$\dim \Acl = Z_{00} \, \Dim\,\Cc \neq 0$,
and by lemma \ref{lem:modinv-simp-hapl}, $\Acl$ 
is haploid.  It then follows from 
theorem \ref{thm:modinv-dim} that
$\Acl$ is special. 
By proposition \ref{prop:symm-T-symm}, $T(\Acl)$ is
a special symmetric Frobenius algebra in $\Cc$. 
Thus $T(\Acl) = \oplus_i A_i$, where the $A_i$ are
simple symmetric Frobenius algebras. We will show that at
least one of the $A_i$ is special. Since $T(\Acl)$ is special,
we have 
$m_{T(\Acl)} \circ \Delta_{T(\Acl)} = \zeta \, \id_{T(\Acl)}$
for some $\zeta \in \Cb^\times$. Restricting this to the summand
$A_i$ shows $m_i \circ \Delta_i = \zeta \, \id_{A_i}$. 
Furthermore, $\eps_{T(\Acl)} \circ \eta_{T(\Acl)} = \xi \, \id_\one$
for some $\xi \in \Cb^\times$. But 
$\eps_{T(\Acl)} \circ \eta_{T(\Acl)} = \sum_i \eps_i \circ \eta_i$,
and so at least one of the $\eps_i \circ \eta_i$ has to be nonzero.
Therefore, at least one of the $A_i$ is special; let $A \equiv A_i$
be this summand. We
denote the embedding $A \hookrightarrow T(\Acl)$ by $e_0$
and the restriction $T(\Acl) \twoheadrightarrow A$ by $r_0$. 
Notice that $r_0$ is an algebra homomorphism. Define
\be
  \ico= \hat{\chi}(r_0) ~:~ \Acl\longrightarrow R(A)~.
\ee
By proposition 
\ref{prop:f-til-f}, $\ico$ is an algebra homomorphism.
Next we verify the centre condition \eqref{eq:left-comm}, 
or rather its
equivalent form \eqref{eq:comm-C}. 
By substituting the definitions, one can convince oneself that
the commutativity 
$m_\text{\rm{cl}} \circ c_{\Acl,\Acl} = m_\text{\rm{cl}}$
of $\Acl$ in $\CcC$ implies the condition
$m_{T(\Acl)} \circ \Gamma = m_{T(\Acl)}$ in $\Cc$, see \cite[prop.\,3.6]{ocfa}. Here,
$\Gamma : T(\Acl) \otimes T(\Acl) \rightarrow T(\Acl) \otimes T(\Acl)$ is given by
\be
  \Gamma = 
  \bigoplus_{m,n}
  \big(\id_{C^l_n} \otimes c_{C^l_m,C^r_n} \otimes \id_{C^r_m}\big)
  \circ
  \big(c_{C^l_m,C^l_n} \otimes c_{C^r_n,C^r_m}^{-1}\big)
  \circ
  \big(\id_{C^l_m} \otimes c_{C^l_n,C^r_m}^{-1} \otimes \id_{C^r_n}\big) ~,
\ee
and we decomposed $\Acl$ as $\Acl = \oplus_n C^l_n \ti C^r_n$.
As a consequence we obtain the identity
\be
  r_0 \circ m_{T(\Acl)} \circ \Gamma \circ 
  (\id_{\Acl} \otimes e_0)
  =
  r_0 \circ m_{T(\Acl)} \circ 
  (\id_{\Acl} \otimes e_0) ~.
\ee
Using that $r_0$ is an algebra map, and that 
by definition $\itco = r_0$, we obtain \eqref{eq:comm-C}.

In order to show that $(A|\Acl,\ico)$ is a Cardy algebra,
it remains to show that the Cardy condition \eqref{eq:cardy-CC}
is satisfied. We will demonstrate this via a detour by
first proving that $\Acl \cong Z(A)$ as Frobenius algebras.

Recall the notations $e$ and $r$ given in \eqref{eq:eZ-rZ-def}.
Using the centre condition \eqref{eq:comm-C} one can check that
$P^l_{R(A)} \circ \ico = m_{R(A)} \circ \Delta_{R(A)}
\circ \ico$. By specialness of $A$ we have
$m_{A} \circ \Delta_{A} = \zeta_{A} \, \id_{A}$
and so together with 
$e \circ r = \zeta_{A}^{-1} \, P^l_{R(A)}$
we get,
\be
  e \circ r \circ \ico =  \ico ~.
\labl{eq:exist-aux1}
Next, consider the morphism 
\be
  f_\text{\rm{cl}} = r \circ \ico
  ~:~ \Acl \longrightarrow Z(A) ~.
\ee  
By the same derivation as
in \eqref{eq:unique-fcl-unit} and \eqref{eq:unique-fcl-mult}
one sees that $f_\text{\rm{cl}}$ is an algebra map. 
In particular, $f_\text{\rm{cl}} \circ \eta_\text{\rm{cl}} = 
\eta_{Z(A)} \neq 0$ and so $f_\text{\rm{cl}} \neq 0$.
Since $\Acl$ is simple, $f_\text{\rm{cl}}$ has to be a 
monomorphism. 

By the same argument as used in the proof of theorem
\ref{thm:modinv-dim}\,(ii), up to isomorphism 
$\Acl$ is the unique simple local $\Acl$-(left-)module.
The algebra monomorphism $f_\text{\rm{cl}}$ turns $Z(A)$ into
an $\Acl$-module. 
Since $Z(A)$ is commutative, it is local as an $\Acl$-module,
and so $Z(A) \cong \Acl^{\oplus N}$ for some $N \ge 1$.
By construction, $A$ is a simple special symmetric 
Frobenius algebra. 
Proposition \ref{prop:symm-R-symm} and corollary
\ref{cor:A-simp-RA-simp} show that $R(A)$
inherits all these properties, and thus $Z(A)$ is
simple (see the comment below equation \eqref{eq:ZA-alg}).
By theorem \ref{thm:reconst}, $Z(A)$ is modular
invariant, and then by theorem \ref{thm:modinv-dim}\,(i), 
$\dim Z(A) = \Dim\,\Cc$. This implies that
$N=1$ in $Z(A) \cong \Acl^{\oplus N}$, and so
$f_\text{\rm{cl}}$ is in fact an isomorphism.

Since $\Acl$ and $Z(A)$ are both haploid,
we have $\eps_{Z(A)} \circ f_\text{cl} = \xi \, \eps_\text{cl}$
for some $\xi \in \Cb^\times$.
The counit uniquely determines the Frobenius structure on $\Acl$
and $Z(A)$ (see e.g.\ \cite[lemma 3.7]{tft1}), so that
$f_\text{cl}$ is a coalgebra isomorphism iff $\xi=1$.
To compute $\xi$ we compose the above identity with
$\eta_\text{cl}$ from the right. Defining $\zeta_\text{cl}$ via
$\eps_\text{cl} \circ \eta_\text{cl} = \zeta_\text{cl}^{-1}
\Dim\,\Cc\cdot\id_\one$ and using \eqref{eq:ZA-counit-norm}
gives $\xi = \dim A \,\zeta_\text{cl}/\zeta_A^2$.
By rescaling the comultiplication and the
counit of $A$, and consequently changing $\zeta_A$, we
can always achieve $\xi=1$.
This proves part (i) of the theorem.

Equation \eqref{eq:exist-aux1} implies that 
$\ico = e \circ f_\text{\rm{cl}}$. Since $f_\text{\rm{cl}}$
is an isomorphism of Frobenius algebras, by lemma 
\ref{lem:f-star-prop} and \ref{lem:ZA-e-star} we have
\be
  \ico \circ \ico^*
  = e \circ f_\text{\rm{cl}} \circ f_\text{\rm{cl}}^* \circ e^*
  = \zeta_{A} \,
  e \circ r
  = P^l_{R(A)} ~.
\ee
Thus $(A|\Acl,\ico)$ is a Cardy algebra. This proves
part (ii) of the theorem.

Part (iii) can be seen as follows. 
By \cite[prop.\,4.3]{morita}, 
$TZ(A) \cong \oplus_{\kappa\in \mathcal{J}}\,
M_{\kappa}^{\vee} \otimes_{A} M_{\kappa}$ as algebras. 
Together with the observation that
$T(f_\text{\rm{cl}}) : T(\Acl) \rightarrow TZ(A)$ is
an isomorphism of algebras, this proves part (iii).
\epf

\begin{rema} {\rm
Part (i) of theorem \ref{thm:open-exists} was announced by
M\"uger \cite{mug-conf}. We provide an independent proof 
in the setting of Cardy algebras
}
\end{rema}

The above theorem, together with lemma \ref{lem:modinv-simp-hapl} 
and theorem \ref{thm:modinv-dim},
shows that a simple commutative symmetric Frobenius 
$\CcC$-algebra $\Acl$ with $\dim\Acl=\Dim\,\Cc$ 
is always part of a Cardy algebra $(\Aop|\Acl,\ico)$ for
some simple special symmetric
Frobenius algebra $\Aop$ in $\Cc$. However, the above proof
also illustrates that $\Aop$ is not unique. 
This raises the question how two
Cardy algebras with a given
$\Acl$ can differ. 
This question is answered by \cite[thm.\,1.1]{morita}, which
in the present framework can be restated as follows.

\begin{thm} \label{thm:morita}
If $(\Aop^{(i)}|\Acl^{(i)}, \ico^{(i)}), i=1,2$ 
are two Cardy $\Cc|\CcC$-algebras 
such that $\Acl^{(i)}$ is simple and 
$\dim \Aop^{(i)}\neq 0$ for $i=1,2$, 
then $\Acl^{(1)} \cong \Acl^{(2)}$ as algebras if and only if 
$\Aop^{(1)}$ and $\Aop^{(2)}$ are Morita equivalent. 
\end{thm}
\pf
Theorem 1.1 in \cite{morita} is stated for $\Aop^{(i)}$ being 
non-degenerate algebras and $\Acl^{(i)} = Z(\Aop^{(i)})$ for $i=1,2$. 
By proposition \ref{prop:Aop-special}, $\Aop^{(i)}$ are simple and special for
$i=1,2$. Then by \cite[lem.\,2.1]{morita}, $\Aop^{(i)}$ are non-degenerate 
algebras. By Theorem \ref{thm:unique}, we have 
$\Acl^{(i)} \cong Z(\Aop^{(i)})$ as Frobenius algebras.
Finally, by  \cite[thm.\,1.1]{morita}, 
$Z(\Aop^{(1)}) \cong Z(\Aop^{(2)})$ as algebras iff
$\Aop^{(1)}$ and $\Aop^{(2)}$ are Morita equivalent. 
\epf

\medskip

Let $C_\text{\rm max}(\CcC)$ be the set of equivalence
classes $[B]$
of of simple modular invariant commutative symmetric Frobenius algebras $B$ in $\CcC$. Two such algebras $B$ and $B'$ 
are equivalent if $B$ and $B'$
are isomorphic as {\em algebras} (but not necessarily
as Frobenius algebras). 
Let $M_\text{\rm simp}(\Cc)$ be the set of 
Morita classes of simple special symmetric Frobenius algebras
in $\Cc$. 
Define the map 
$z: M_\text{\rm simp}(\Cc) \rightarrow C_\text{\rm max}(\CcC)$ by 
$z: \{A\} \rightarrow [Z(A)]$, where $\{A\}$ denotes
the Morita class of $A$.
From theorem \ref{thm:open-exists}\,(i) 
and \cite[thm.\,1.1]{morita} we learn:

\begin{cor}
The map $z: M_\text{\rm simp}(\Cc) \rightarrow C_\text{\rm max}(\CcC)$ 
is a bijection. 
\end{cor}

\appendix 

\sect{Appendix}

\subsection{Proof of lemma \ref{lem:lax-colax-adj} }
\label{app:lax-colax}

We will show that if
$(F, \psi_2^F, \psi_0^F)$ is a colax tensor functor from $\Cc_1$ to $\Cc_2$,
then $(G, \phi_2^G, \phi_0^G)$ is a lax tensor functor from $\Cc_2$ to $\Cc_1$. 
Applying this result to the opposed categories then gives
the converse statement.

We need to show that $\phi_0^G$ and $\phi_2^G$ make the
diagrams (\ref{R-ten-fun-1}) and (\ref{R-ten-fun-2})
commute. 
We first prove the commutativity of (\ref{R-ten-fun-1}). 
Consider the following diagram:
\be  \raisebox{4em}{
\xymatrix{
G(A)\otimes (G(B)\otimes G(C)) \ar[d]^{G\psi_2^F \circ \delta}
\ar[rrr]^{\id_{G(A)}\otimes \phi_2^G}  &&& G(A)\otimes G(B\otimes C)  
\ar[d]^{G\psi_2^F \circ \delta} \\  
G(FG(A)\otimes F(G(B)\otimes G(C))) \ar[rrr]^{\hspace{0.3cm} 
G(F(\id_{G(A)}) \otimes F(\phi_2^G))} 
\ar[d]^{G(\id_{FG(A)}\otimes \psi_2^F)} &&&
G(FG(A)\otimes FG(B\otimes C)) 
\ar[d]^{G(\rho_A \otimes \rho_{B\otimes C})}  \\
G(FG(A)\otimes (FG(B)\otimes FG(C))) \ar[rrr]^{\hspace{0.3cm} 
G(\rho_A \otimes 
(\rho_B \otimes \rho_C))} &&& G(A\otimes (B\otimes C))
}} 
\labl{R-ten-fun-1-pf-1}
The top subdiagram is commutative because of the naturality of 
$G\psi_2^F \circ \delta$. The commutativity of bottom subdiagram
follows from the following identities: 
\bea
(\rho_B \otimes \rho_C) \circ \psi_2^F &=&
(\rho_B \otimes \rho_C) \circ \psi_2^F \circ \rho F \circ F\delta  \nn
&=& \rho_{B\otimes C} \circ FG(\rho_B\otimes \rho_C) 
\circ FG(\psi_2^F) \circ F\delta  \nn
&=& \rho_{B\otimes C} \circ F(\phi_2^G) 
\eea
as a map $F(G(B)\otimes G(C)) \rightarrow B\otimes C$. 
The commutativity of (\ref{R-ten-fun-1-pf-1}) implies that
the composition of maps in the left column in (\ref{R-ten-fun-1}) 
can be replaced by
\be  \label{R-ten-fun-1-pf-L} 
G(\rho_A \otimes 
(\rho_B \otimes \rho_C)) \circ G\big(
  (\id_{FG(A)}\otimes \psi_2^F)
  \circ \psi_2^F \big) \circ \delta~.
\ee
Similarly, we can show that the composition of maps in the 
right column in (\ref{R-ten-fun-1}) can be replaced by 
\be  \label{R-ten-fun-1-pf-R}
G( (\rho_A \otimes 
\rho_B) \otimes \rho_C ) \circ G\big(
  (\psi_2^F \otimes \id_{FG(C)}) \circ \psi_2^F 
  \big) \circ \delta~.
\ee
Using the commutativity of (\ref{L-ten-fun-1}), it is easy to see
that (\ref{R-ten-fun-1}) with the left and right columns 
of (\ref{R-ten-fun-1}) replaced by
(\ref{R-ten-fun-1-pf-L}) and (\ref{R-ten-fun-1-pf-R}) respectively
is commutative. Hence (\ref{R-ten-fun-1}) is commutative. 

Now we prove the commutativity of the first diagram in (\ref{R-ten-fun-2}). 
\bea
&& \rule{35em}{0pt} \nonumber \\[-1.1em]
&&\phi_2^G \circ (\phi_0^G \otimes \id_{G(A)}) \nn
&&\hspace{2em}\overset{(1)}{=}   
G(\rho_{\one_2}\otimes \rho_A) \circ G\psi_2^F \circ
\delta \circ \big[ (G\psi_0^F \circ \delta_{\one_1}) \otimes \id_{G(A)}\big] 
\nn 
&&\hspace{2em}\overset{(2)}{=} 
G(\rho_{\one_2}\otimes \rho_A) \circ G\psi_2^F 
\circ
GF(G\psi_0 \otimes \id_{G(A)}) \circ GF(\delta_{\one_1} \otimes \id_{G(A)})
\circ \delta  \nn
&&\hspace{2em}\overset{(3)}{=} G(\rho_{\one_2}\otimes \rho_A) \circ 
G(FG(\psi_0^F)\otimes \id_{FG(A)}) \circ G\psi_2^F \circ 
GF(\delta_{\one_1} \otimes \id_{G(A)}) \circ \delta  \nn
&&\hspace{2em}\overset{(4)}{=} G(\id_{\one_2} \otimes \rho_A) \circ 
G\big( [\psi_0^F \circ \rho_{F(\one_1)} \circ (F\delta)_{\one_1}] 
\otimes \id_{FG(A)}\big) \circ G\psi_2^F \circ \delta \nn
&&\hspace{2em}\overset{(5)}{=} G(\id_{\one_2} \otimes \rho_A) \circ 
G(\psi_0^F \otimes \id_{FG(A)}) \circ G\psi_2^F \circ \delta  
\nn
%
&&\hspace{2em}\overset{(6)}{=} G(\id_{\one_2} \otimes \rho_A) 
\circ G(l_{FG(A)}^{-1}) \circ GF(l_{G(A)}) \circ \delta \nn
&&\hspace{2em}\overset{(7)}{=} G(l_A^{-1}) \circ G \rho_A \circ \delta G 
\circ l_{G(A)}  \nn
&&\hspace{2em}\overset{(8)}{=} G(l_A^{-1}) \circ l_{G(A)}.
\eea
where in step (1) we substituted the definition of $\phi_0^G, \phi_2^G$ given in 
(\ref{eq:phi2F-def-psi}); in step (2) we used the naturality of $\delta$; in
step (3) we used the naturality of 
$G\psi_2^F$; 
in step (4) we switched the position 
between $G\rho_{\one_2}$ and $GFG(\psi_0^F)$ and the position between
$G\psi_2$ and $GF(\delta_{\one_2}\otimes \id_{G(A)})$ using the naturality of
$\rho$ and $F\psi_2^G$ respectively; in step (5) we applied the second identity in 
(\ref{eq:adj-unit-counit}); in step (6) we used (\ref{L-ten-fun-2}); in step (7) we used
the naturality of $l^{-1}$
and $\delta$; in step (8) we used the 
first identity in (\ref{eq:adj-unit-counit}).

The proof of the commutativity of the second diagram
in (\ref{R-ten-fun-2}) is similar. 
Thus we have shown that $G$ is a lax tensor functor. 
\epf

\subsection{Proof of lemma \ref{lemma:mod-inv-RA}}
\label{app:RA-Sinv}

To prepare the proof, recall that for a given
object $B\in \Cc$,
the modular group $PSL(2,\Zb)$ acts on the space
$\oplus_i \Hom_\Cc(B \oti U_i,U_i)$,
see  
e.g.\ \cite[sect\,3.1]{baki} and \cite[eqn.\,(4.55)]{cardy}.
We will only need the action of $S$ and $S^{-1}$. Let
$f \in \oplus_i \Hom_\Cc(B \oti U_i,U_i)$.
Then 
\be   \label{eq:S-map}
S~:~ \quad
  \bigoplus_{i \in \Ic}
  \raisebox{-30pt}{
  \begin{picture}(40,90)
   \put(0,8){\scalebox{.75}{\includegraphics{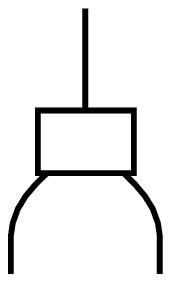}}}
   \put(0,8){
     \setlength{\unitlength}{.75pt}\put(-18,-19){
     \put(12, 8)     {\scriptsize $ B $}
     \put(60,8)      {\scriptsize $U_i$ }
     \put(38,105)  {\scriptsize $U_i$ }
     \put(36,56)    {\scriptsize $f$}
     }\setlength{\unitlength}{1pt}}
  \end{picture}}
\quad \longmapsto \quad 
\bigoplus_{j\in \Ic} 
\,\,\frac{\dim U_j}{\sqrt{\Dim \Cc }}\,\,
\sum_{i\in \Ic} 
 \raisebox{-35pt}{
  \begin{picture}(90,90)
   \put(8,8){\scalebox{.75}{\includegraphics{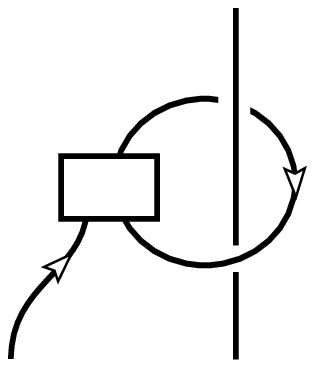}}}
   \put(0,8){
     \setlength{\unitlength}{.75pt}\put(-18,-19){
     \put(24, 8)  {\scriptsize $ B $}
     \put(90, 130)  {\scriptsize $ U_j $}
     \put(90, 8)  {\scriptsize $ U_j $}
     \put(50, 90)    {\scriptsize $ U_i $}
     \put(53, 67)     {\scriptsize $f$}
}\setlength{\unitlength}{1pt}}
  \end{picture}}
  ~~,
\ee
\be   \label{eq:S--1map}
\hspace{0.5cm}
S^{-1}~:~ \quad
  \bigoplus_{i \in \Ic}
  \raisebox{-30pt}{
  \begin{picture}(40,90)
   \put(0,8){\scalebox{.75}{\includegraphics{pic-Smapf-L.eps}}}
   \put(0,8){
     \setlength{\unitlength}{.75pt}\put(-18,-19){
     \put(12, 8)     {\scriptsize $ B $}
     \put(60,8)      {\scriptsize $U_i$ }
     \put(38,105)  {\scriptsize $U_i$ }
     \put(36,56)    {\scriptsize $f$}
     }\setlength{\unitlength}{1pt}}
  \end{picture}}
\quad \longmapsto \quad 
\bigoplus_{j\in \Ic} 
\,\,\frac{\dim U_j}{\sqrt{\Dim \Cc }}\,\,
\sum_{i\in \Ic}
 \raisebox{-35pt}{
  \begin{picture}(90,90)
   \put(8,8){\scalebox{.75}{\includegraphics{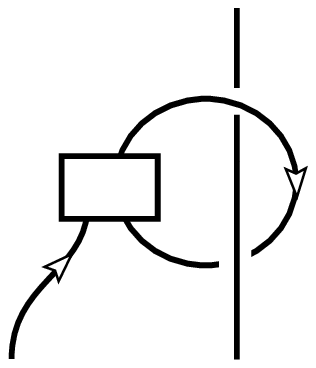}}}
   \put(0,8){
     \setlength{\unitlength}{.75pt}\put(-18,-19){
     \put(24, 8)  {\scriptsize $ B $}
     \put(90, 130)  {\scriptsize $ U_j $}
     \put(90, 8)  {\scriptsize $ U_j $}
     \put(50, 90)    {\scriptsize $ U_i $}
     \put(53, 67)     {\scriptsize $f$}
}\setlength{\unitlength}{1pt}}
  \end{picture}}
  ~~.
\ee

By lemma \ref{lemma:S-inv}, to establish that \eqref{eq:RA-S-inv}
is S-invariant, it is enough to prove the identity \eqref{eq:mod-inv-2}
when $f$ is given by \eqref{eq:RA-S-inv}.
Using \eqref{eq:S-map} and \eqref{eq:S--1map}, we can see that 
equation \eqref{eq:mod-inv-2} simply says 
that $\oplus_{i,j}[ \mbox{RHS of \eqref{eq:mod-inv-2}}]$ 
is invariant under the action of $S\times S$. 
Consider the element $g$ of
$\oplus_{j,k\in \Ic} \Hom_{\CcC}(R(A) \otimes (U_j^{\vee}\times U_k),  
U_j^{\vee}\times U_k)$ given by
\be
g \quad =  \bigoplus_{j,k\in \Ic} ~ \sum_\alpha ~~  
  \raisebox{-67pt}{
  \begin{picture}(110,150)
   \put(0,8){\scalebox{.75}{\includegraphics{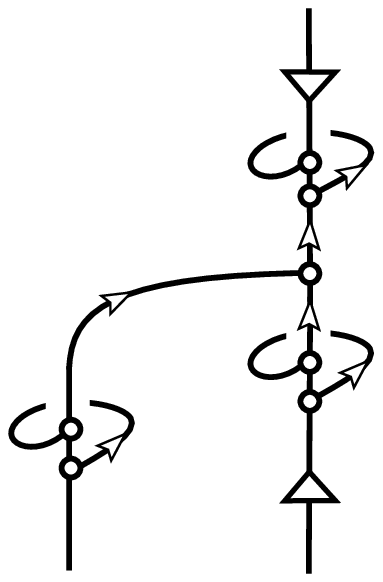}}}
   \put(0,8){
     \setlength{\unitlength}{.75pt}\put(-18,-19){
     \put(25, 10)       {\scriptsize $ R(A) $}
     \put(120,128)    {\scriptsize $ R(A) $}
     \put(120,95)      {\scriptsize $ R(A) $}
     \put(35, 115)     {\scriptsize $ R(A) $}
     \put(97,10)      {\scriptsize $ U_j^{\vee} \times U_k$ }
     \put(97,195)    {\scriptsize $ U_j^{\vee} \times U_k$ }
     \put(116, 163)   {\scriptsize $\alpha$ }
     \put(116, 45)     {\scriptsize $\alpha$ }
     }\setlength{\unitlength}{1pt}}
  \end{picture}} 
  \quad .
\ee
By the above arguments,  
proving S-invariance of \eqref{eq:RA-S-inv} is 
equivalent to proving invariance of $g$
under the action of $S \times S$.

For $i \in \Ic$, we denote by $g_i$ the component of $g$ in 
$$
\big( \oplus_{j\in \Ic} \Hom_{\Cc}(A\otimes U_i^{\vee}\otimes U_j^{\vee}, 
U_j^{\vee}) \big)
\otimes \big( \oplus_{k\in \Ic} \Hom_{\Cc}(U_i\otimes U_k, U_k)\big) ~.
$$
We view the second Hom-space in 
above tensor product as a Hom-space in $\Cc_+$ instead of $\Cc_-$. 
It is enough to show that 
$g_i$ is invariant under the action of $S\times S^{-1}$. 
Note that the action of $S^{-1}$ in $\Cc_+$ is equivalent to that of 
$S$ in $\Cc_-$.   

The morphism $g_i$ can be canonically identified with a bilinear pairing
\be
  (\,\cdot \,,\, \cdot\,)_i ~:~ 
  \big( \bigoplus_{j\in \Ic} \Hom_{\Cc}(U_j^{\vee}, 
  A\otimes U_i^{\vee}\otimes U_j^{\vee})\big) 
  \times \big( \bigoplus_{k\in \Ic} \Hom_{\Cc}(U_k, U_i\otimes U_k)\big)
  ~ \longrightarrow ~
  \Cb 
\ee
as follows. 
For $h_1\in \Hom_{\Cc}(U_j^{\vee}, A\otimes U_i^{\vee}\otimes U_j^{\vee})$ 
and $h_2\in \Hom_{\Cc}(U_k, U_i\otimes U_k)$ we set
\be
  (h_1,h_2)_i ~=~
  (\dim U_j \dim U_k)^{-1} \,
  \mathrm{tr}_{U_j^\vee \times U_k}\big[ g_i \circ (h_1 \times h_2) \big] ~.
\ee
When substituting the explicit form of the product $m_{R(A)}$ of $R(A) = (A \ti \one) \otimes R(\one)$, after a short calculation one finds 
\be  
(h_1, h_2)_i = \sum_{\alpha} \frac{1}{\dim U_j} \quad
  \raisebox{-100pt}{
  \begin{picture}(190,195)
   \put(0,8){\scalebox{.60}{\includegraphics{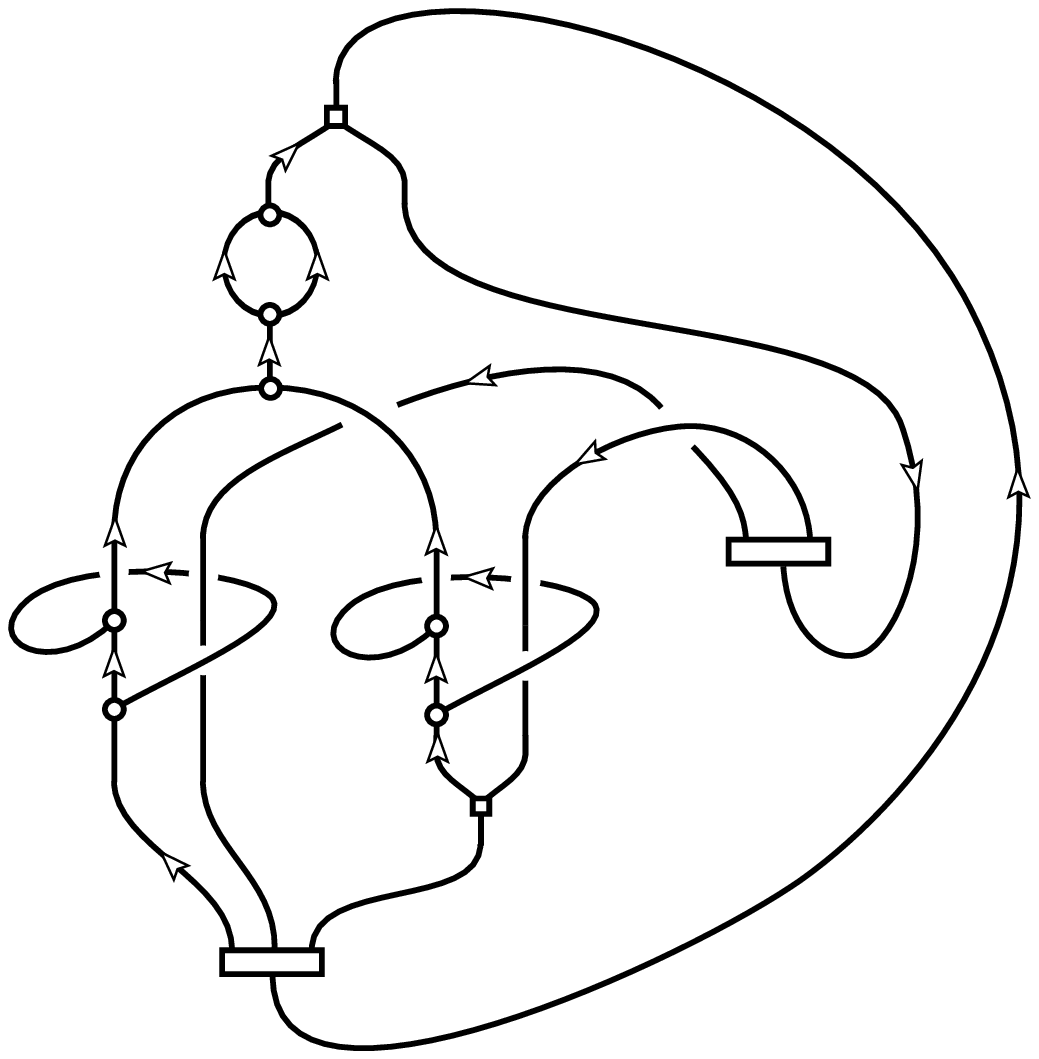}}}
   \put(0,8){
     \setlength{\unitlength}{.60pt}\put(-18,-19){
     \put(53, 65)       {\scriptsize $ A $}
     \put(57, 163)     {\scriptsize $ A $}
     \put(128,102)    {\scriptsize $ A $}
     \put(150,163)    {\scriptsize $ A $}
     \put(64, 245)     {\scriptsize $A$}
     \put(113, 245)   {\scriptsize $A$ }
     \put(82, 280)     {\scriptsize $A$}
     \put(109, 278)   {\scriptsize $\alpha$ }
     \put(151, 99)   {\scriptsize $\alpha$ }
     \put(62,43)     {\scriptsize $h_1$ }
     \put(208,162)   {\scriptsize $h_2$ }
     \put(308,220)   {\scriptsize $U_j$ }
     \put(215, 230)   {\scriptsize $ U_k $ }
     \put(120, 76)     {\scriptsize $ U_j^{\vee} $ }
     \put(148,198)        {\scriptsize $ U_i $ }
     \put(182,176)    {\scriptsize $ U_k $ }
     }\setlength{\unitlength}{1pt}}
  \end{picture}} 
~~.
\labl{eq:biform}
Here the top morphism $P^l_{R(A)}$ has been simplified with the
help of the identity
\be\begin{array}{l}\displaystyle
P_{R(A)} \circ m_{R(A)} \circ (P_{R(A)} \otimes P_{R(A)})
\enl
\quad =
\big( \bigoplus_{i\in \Ic} ((m_A\circ \Delta_A) \otimes \id_{U_i^{\vee}}) \times \id_{U_i} \big)
\circ m_{R(A)} \circ ( P_{R(A)} \otimes P_{R(A)} )~~,
\eear
\labl{eq:PmPP-mPP}
which can be checked by direct calculation along the same lines as in the
proof of \cite[lem\,3.10]{corr}.

The action of the modular transformation $S$ on
$\oplus_{i\in \Ic} (B\otimes U_i, U_i)$
for $B\in \Cc$ 
naturally induces an action on
$\oplus_{i\in \Ic} (U_i,B\otimes U_i)$ \cite[prop.\,5.14]{cardy},
which we denote by $S^*$.
In the present case we get an action of $S^*$ on
$\oplus_{j\in \Ic} \Hom_{\Cc}(U_j^{\vee}, 
A\otimes U_i^{\vee}\otimes U_j^{\vee})$ and 
$\oplus_{k\in \Ic} \Hom_{\Cc}(U_k, U_i\otimes U_k)$.
Then to show $g_i$ is invariant under the action of $S\times S^{-1}$
amounts to showing that 
\be
(h_1, h_2)_i = ( (S^{-1})^*\, h_1, S^* \, h_2)_i ~,
\ee
for all $h_1\in 
\Hom_{\Cc}(U_j^{\vee}, A\otimes U_i^{\vee}\otimes U_j^{\vee})$ 
and $h_2\in 
\Hom_{\Cc}(U_k, U_i\otimes U_k)$. 
We have
\be  \label{eq:biform-S}
( (S^{-1})^*\, h_1, S^* \, h_2)_i
=\sum_{m,n,\alpha} \frac{\dim U_n}{\Dim \Cc} \quad
  \raisebox{-100pt}{
  \begin{picture}(190,200)
   \put(0,8){\scalebox{.60}{\includegraphics{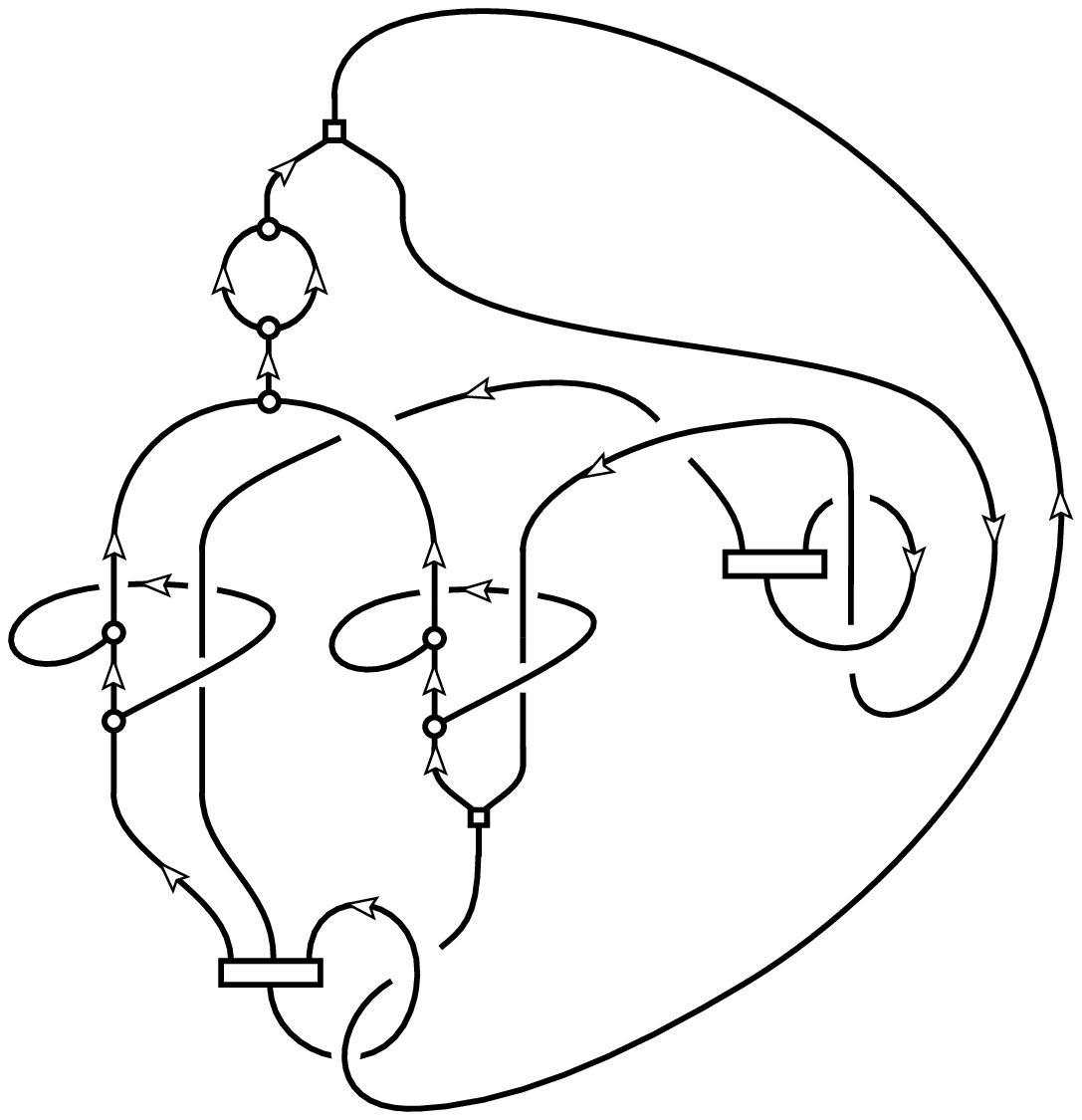}}}
   \put(0,8){
     \setlength{\unitlength}{.60pt}\put(-18,-19){
     \put(47, 83)       {\scriptsize $ A $}
     \put(57, 181)     {\scriptsize $ A $}
     \put(125,115)    {\scriptsize $ A $}
     \put(150,177)    {\scriptsize $ A $}
     \put(64, 258)     {\scriptsize $A$}
     \put(113, 258)   {\scriptsize $A$ }
     \put(82, 293)     {\scriptsize $A$}
     \put(225, 158)   {\scriptsize $ U_k $ }
     \put(111, 293)   {\scriptsize $\alpha$ }
     \put(150, 114)   {\scriptsize $\alpha$ }
     \put(63,59)     {\scriptsize $h_1$ }
     \put(208,177)   {\scriptsize $h_2$ }
     \put(150,213)    {\scriptsize $ U_i $}
     \put(106, 88)     {\scriptsize $ U_j $ }
     \put(328,200)     {\scriptsize $ U_m $}
     \put(185,192)     {\scriptsize $ U_n $}
     }\setlength{\unitlength}{1pt}}
  \end{picture}} 
~~.
\ee
Now drag the upper vertex indexed by $\alpha$ 
in the above graph along its $U_m^{\vee}$-leg until it meets the lower vertex
also indexed by $\alpha$, then sum over $\alpha$ and $m$. 
This gives
\be  \label{eq:biform-S-alpham}
( (S^{-1})^*\, h_1, S^* \, h_2)_i
=\sum_{n} \frac{\dim U_n}{\Dim \Cc }  \quad 
  \raisebox{-100pt}{
  \begin{picture}(190,200)
   \put(0,8){\scalebox{.60}{\includegraphics{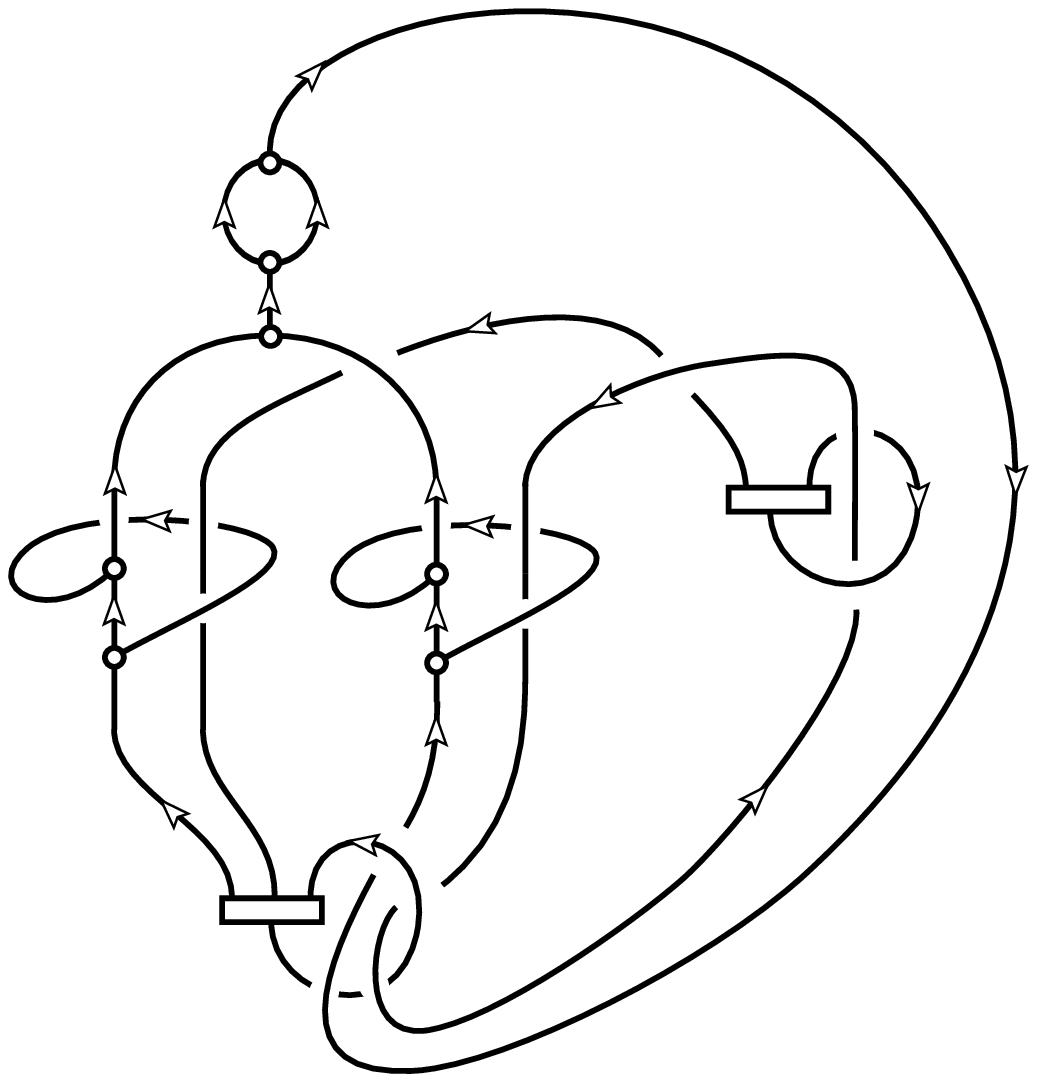}}}
   \put(0,8){
     \setlength{\unitlength}{.60pt}\put(-18,-19){
     \put(47, 91)       {\scriptsize $ A $}
     \put(57, 189)     {\scriptsize $ A $}
     \put(125,115)    {\scriptsize $ A $}
     \put(150,185)    {\scriptsize $ A $}
     \put(64, 266)     {\scriptsize $A$}
     \put(113, 266)   {\scriptsize $A$ }
     \put(82, 299)     {\scriptsize $A$}
     \put(226, 163)   {\scriptsize $ U_k $ }
     \put(215, 105)   {\scriptsize $ U_n $ }
     \put(63,67)     {\scriptsize $h_1$ }
     \put(208,185)   {\scriptsize $h_2$ }
     \put(150,244)    {\scriptsize $ U_i $}
     \put(106, 94)     {\scriptsize $ U_j $ }
     }\setlength{\unitlength}{1pt}}
  \end{picture}} 
~~.
\ee
If we just look at the neighbourhood of the $U_n$-loop in the above graph, we 
see the following subgraphs,
\be  \label{eq:Umloop}
\sum_{n} \frac{\dim U_n}{\Dim \Cc }  \quad 
  \raisebox{-50pt}{
  \begin{picture}(95,100)
   \put(0,8){\scalebox{.75}{\includegraphics{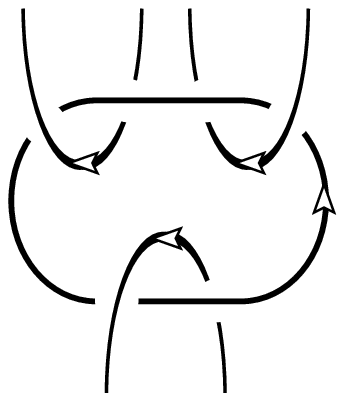}}}
   \put(0,8){
     \setlength{\unitlength}{.75pt}\put(-18,-19){     
     \put(18, 136)     {\scriptsize $A$}
     \put(50, 136)   {\scriptsize $A^{\vee}$ }
     \put(43,8)      {\scriptsize $ U_j^{\vee} $ }
     \put(76, 8)   {\scriptsize $ U_j $ }
     \put(68,136)    {\scriptsize $ U_k $ }
     \put(102, 136)     {\scriptsize $ U_k^{\vee} $ }
     \put(115, 70)   {\scriptsize $ U_n $ }
     }\setlength{\unitlength}{1pt}}
  \end{picture}} 
= \,\, \sum_{\alpha} \frac{1}{\dim U_j} \quad
 \raisebox{-50pt}{
  \begin{picture}(100,100)
   \put(0,8){\scalebox{.75}{\includegraphics{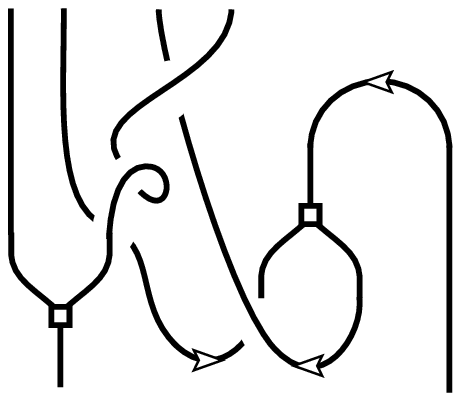}}}
   \put(0,8){
     \setlength{\unitlength}{.75pt}\put(-18,-19){
     \put(13, 136)       {\scriptsize $ A $}
     \put(30, 136)     {\scriptsize $ A^{\vee} $}
     \put(57,136)      {\scriptsize $ U_k $ }
     \put(80, 136)   {\scriptsize $ U_k^{\vee} $ }
     \put(27,8)    {\scriptsize $ U_j^{\vee} $ }
     \put(142, 8)     {\scriptsize $ U_j $ }
     \put(30, 51)    {\scriptsize  $ \alpha $}
     \put(102, 60)    {\scriptsize  $ \alpha $}
     }\setlength{\unitlength}{1pt}}
  \end{picture}} 
~~,
\ee
where we have applied \cite[cor.\,3.1.11]{baki}. 
Substituting this subgraph back to the original graph in 
(\ref{eq:biform-S-alpham}), we obtain 
\be  \label{eq:biformS-loopsum}
( (S^{-1})^*\, h_1, S^* \, h_2)_i
=\sum_{\alpha} \frac{1}{\dim U_j} \quad
  \raisebox{-100pt}{
  \begin{picture}(190,200)
   \put(0,8){\scalebox{.60}{\includegraphics{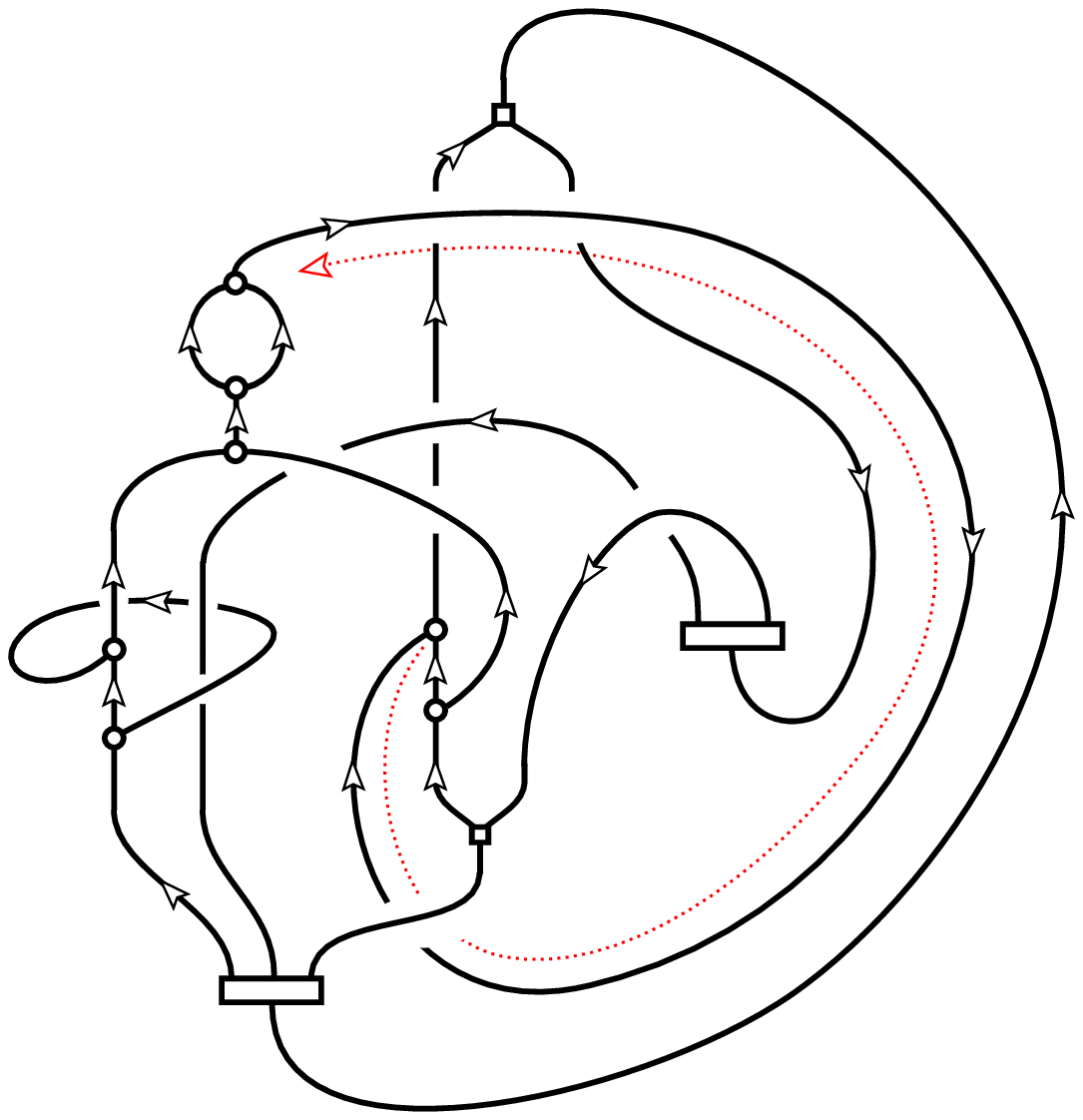}}}
   \put(0,8){
     \setlength{\unitlength}{.60pt}\put(-18,-19){
     \put(48, 75)       {\scriptsize $ A $}
     \put(57, 173)     {\scriptsize $ A $}
     \put(132,100)    {\scriptsize $ A $}
     \put(148,169)    {\scriptsize $ A $}
     \put(53, 240)     {\scriptsize $A$}
     \put(103, 240)   {\scriptsize $A$ }
     \put(107, 281)     {\scriptsize $A$}
     \put(105,117)     {\scriptsize $A$}
     \put(148, 248)    {\scriptsize $A$}
     \put(160, 80)     {\scriptsize $ U_j^{\vee} $ }
     \put(158, 294)   {\scriptsize $\alpha$ }
     \put(151, 107)   {\scriptsize $\alpha$ }
     \put(63,53)     {\scriptsize $h_1$ }
     \put(196,155)   {\scriptsize $h_2$ }
     \put(118,165)   {\scriptsize $m_A$ }
     \put(161,224)    {\scriptsize $ U_i $}
     \put(328,200)     {\scriptsize $ U_j $}
     \put(220,193)     {\scriptsize $ U_k $}
     \put(225,123)     {\scriptsize $ U_k $}
     }\setlength{\unitlength}{1pt}}
  \end{picture}} 
~~.
\ee
The graph in (\ref{eq:biformS-loopsum}) is equal to
that in (\ref{eq:biform}). 
In order to see this, we first drag the ``bubble'' ($m_A\circ \Delta_A$)
along $A$ lines and through the $m_A$ vertex (because $m_A\circ \Delta_A$
is a bimodule map) until it reaches the lower-left leg of the upper vertex 
indexed by $\alpha$. Then drag the $m_A$ vertex along the (red)
dotted line in above graph. 
Finally, we apply the associativity of $A$, (\ref{eq:PmPP-mPP}), 
and \cite[lem\,3.11]{corr}. Then we see that 
the graph in (\ref{eq:biformS-loopsum}) exactly matches with 
the one in (\ref{eq:biform}). 
\epf

\end{document}